\documentstyle[11pt]{article}
\newcommand{\bC}{{\Bbb C}}

\newcommand{\bT}{{\bf T}}
\newcommand{\cM}{{\cal M}}

\newcommand{\cB}{{\cal B}}
\newcommand{\cG}{{\cal G}}
\newcommand{\cN}{{\cal N}}
\newcommand{\bQ}{{\Bbb Q}}
\newcommand{\cK}{\cal K}
\newcommand{\bh}{{\bf h}}
\newcommand{\cQ}{{\cal Q}}
\newcommand{\bZ}{{\Bbb Z}}

\newcommand{\bF}{{\bf F}}
\newcommand{\cF}{{\cal F}}
\newcommand{\cO}{{\cal O}}
\newcommand{\lU}{{}_U}
\newcommand{\bI}{{\bf I}}
\newcommand{\cW}{{\cal W}}
\newcommand{\bW}{{\bf W}}

\newcommand{\bD}{{\bf D}}
\newcommand{\cV}{{\cal V}}
\newcommand{\bV}{{\bf V}}
\newcommand{\bA}{{\bf A}}
\newcommand{\bH}{{\bf H}}
\newcommand{\bG}{{\bf G}}
\newcommand{\cX}{{\cal X}}
\newcommand{\cI}{{\cal I}}
\newcommand{\cT}{{\cal T}}
\newcommand{\cH}{{\cal H}}
\newcommand{\cC}{{\cal C}}

\newcommand{\cA}{{\cal A}}

\font\twelvefrak=eufm10 at 12pt

\font\sevenfrak=eufm7
\font\fivefrak=eufm5
\newfam\frakfam
\textfont\frakfam=\twelvefrak
\scriptfont\frakfam=\sevenfrak
\scriptscriptfont\frakfam=\fivefrak
\def\frak{\fam\frakfam\twelvefrak}

\font\twelveBbb=msbm10 at 12pt

\font\sevenBbb=msbm7
\font\fiveBbb=msbm5
\newfam\Bbbfam
\textfont\Bbbfam=\twelveBbb
\scriptfont\Bbbfam=\sevenBbb
\scriptscriptfont\Bbbfam=\fiveBbb
\def\Bbb{\fam\Bbbfam\twelveBbb}

\def\lam{\lambda}

\def\S{{\frak S}}
\def\a{{\frak a}}
\def\p{{\frak p}}
\def\g{{\frak g}}
\def\t{{\frak t}}

\def\f{{\frak f}}
\def\q{{\frak q}}

\def\n{{\frak n}}
\def\u{{\frak u}}

\def\U{{\cal U}}
\def\cZ{{\cal Z}}
\def\w{{\bf w}}
\def\C{{\Bbb C}}
\def\m{{\frak m}}

\def\bW{\bf W}
\def\N{{\Bbb N}}
\def\A{{\Bbb A}}
\def\Z{{\Bbb Z}}
\def\Zp{{{\Bbb Z}_p}}
\def\Q{{\Bbb Q}}
\def\R{{\Bbb R}}
\def\C{{\Bbb C}}
\def\F{{{\Bbb F}}}
\def\S{{{\Bbb S}}}
\def\barX{{\overline{X}}}
\def\X{{{\Bbb X}}}
\def\G{{{\Bbb G}}}

\def\sp{\mbox{\bf spin}}
 
\def\adots{\mathinner{\mkern2mu\raise1pt\hbox{.}\mkern3mu\raise4pt\hbox{.}\mkern1mu\raise7pt\hbox{.}}}

\newtheorem{thm}{Theorem}
\newtheorem{cor}{Corollary}
\newtheorem{pro}{Proposition}
\newtheorem{lem}{Lemma}
\newtheorem{sublem}{Sublemma}
\newtheorem{de}{Definition}
\title{Cohomology of Siegel varieties with $p$-adic integral coefficients and Applications}
\author{A. Mokrane and J. Tilouine}
\date{October 26 2000}
\begin{document}
\maketitle

\section{Introduction}

{\bf 1.1.} 
Let $G$ be a connected reductive group over $\Q$. Diamond \cite{D} and
Fujiwara \cite{Fu} (independently) have axiomatized the Taylor-Wiles method which allows to study
some local components
$\bT_\m$ of a Hecke algebra
$\bT$ for $G$ of suitable (minimal) level; when it applies, this method shows at the same time that
$\bT_\m$ is complete intersection and that some cohomology module, viewed as $\bT$-module, is locally
free at
$\m$. It has been successfully applied to
$GL(2)_{/\Q}$ \cite{TW}, to some quaternionic Hilbert modular cases \cite{Fu}, and to some inner
forms of unitary groups \cite{HT}. If one tries to treat other cases, one can let
 the Hecke algebra act
faithfully on the middle degree Betti cohomology of an associated
Shimura variety; then,
one of the problems to overcome is the possible presence of torsion in the cohomology modules with $p$-adic integral
coefficients. For $G=GSp(2g)$ ($g\geq 1$), we want to explain in this paper why this torsion is not
supported by maximal ideals of $\bT$ which are ``{non-Eisenstein}'' and ordinary (see below for precise
definitions), provided the residual characteristic
$p$ is prime to the level and greater than a natural bound. A drawback of our method is that it necessitates
to assume that the existence and some local properties of the Galois representations associated to cohomological
cuspidal representations on $G$ are established. For the moment, they are proven for $g\leq 2$ (see below). In his recent
preprint 
\cite{H2}, Hida explains for the same symplectic groups $G$ how by considering only coherent cohomology, one
can let the Hecke algebra act faithfully too on cohomology modules whose torsion-freeness is
built-in (without assuming any conjecture). However for some applications (like the relation, for some groups $G$, between
special values of adjoint
$L$-functions, congruence numbers, and cardinality of adjoint Selmer groups), the use of the Betti cohomology seems
indispensable.

\vskip 1cm   
\noindent {\bf 1.2.}
Let $G=GSp(2g)$ be the group of symplectic similitudes given by the matrix
$J=\left(\begin{array}{cc} 0&s\\-s&0\end{array}\right)$, whose entries are $g\times g$-matrices, and
$s$ is antidiagonal, with non-zero coefficients equal to $1$; the standard Borel $B$, resp. torus
$T$, in $G$ consists in upper triangular matrices, resp. diagonal matrices in $G$. For any dominant
weight
$\lambda$ for
$(G,B,T)$, we write
$\hat{\lambda}$ for its dual (that is, the dominant weight associated to the Weyl representation dual
of that of
$\lambda$). Let $\rho$ be the half-sum of the positive roots. Recall that $\lambda$ is given by
a $(g+1)$-uple $(a_g,\ldots,a_1;c)\in\Z^{g+1}$ with $c\equiv a_1+\ldots+a_g\, mod.2$, that
$\hat{\lambda}=(a_g,\ldots,a_1;-c)$ and
$\rho=(g,\ldots,1;0)$ (see section 3.1 below). Throughout this paper, the
following integer will be of great importance:  
$$\w=\vert\lambda+\rho\vert=\vert\lambda\vert+d=\sum_{i=1}^g(a_i+i)=d+\sum_{i=1}^g a_i$$
where $d=g(g+1)/2$.
It can be viewed as a cohomological weight as follows.

Let $\A=\A_f\times\Q_\infty$ be the ring of rational adles; let $G_f$ resp. $G_\infty$ be the group
of $\A_f$-points resp. $\Q_\infty$-points of $G$. Let $U$
be a ``{good}'' open compact subgroup of
$G(\bA_f)$ (see Introd. of Sect.2); let $S$ resp.$S_U$ be the Shimura variety of infinite level,
 resp. of level $U$
associated to
$G$; then
$d=dim\,{}S_U$ is the middle degree of the Betti cohomology of $S_U$. Let 
$V_\lambda(\C)$ be the
coefficient system over $S$ resp. $S_U$ with highest weight $\lambda$. See Sect.2.1 for precise
definitions. 

Let
$\pi=\pi_f\otimes\pi_\infty$ be a cuspidal automorphic representation of
$G(\A)$ which occurs in $H^d(S_U,V_\lambda(\C))$. This means that 
\begin{itemize}
\item the
$\pi_f$-isotypical component $W_\pi=H^d(\pi_f)$ of the $G_f$-module $H^\bullet (S,V_\lambda
(\C))$ is non-zero, and
\item  $\pi_f^U\neq 0$.
\end{itemize}

It is known (see Sect.2.3.1 below) that the first condition is equivalent to the fact that
$\pi_\infty$ belongs to the
$L$-packet $\Pi_{\hat{\lambda}+\rho}$ of Harish-Chandra's parameter $\hat{\lambda}+\rho$ in the
discrete series.

By a Tate twist, we can restrict ourselves to the case where
$c=a_g+\ldots+a_1$.  We do this in the sequel. 
Then, $\vert \lambda\vert$ is the Deligne weight of the coefficient system
$V_\lambda$ and $\w=\vert\lambda+\rho\vert$ is the cohomological weight of $H^d(S_U,V_\lambda)$,
hence the (hypothetical) motivic weight of $\pi$.  

\bigskip

Let $p$ be a prime. Let us fix an embedding
$\iota_p:\overline{\Q}\hookrightarrow \overline{\Q}_p$. Let $v$ be the valuation of $\overline{\Q}$
induced by 
$\iota_p$ normalized by
$v(p)=1$; let $K$ be the $v$-adic completion of a number field containing the Hecke eigenvalues of
$\pi$. We denote by $\cO$ the valuation ring of $(K,v)$; we fix a local parameter $\varpi\in\cO$.  Let
$N$ be the level of $U$, that is, the smallest positive integer such that the principal congruence subgroup $U(N)$ is
contained in $U$. Let ${\cal H}^{N}$ resp.
${\cal H}_U(\cO)$ be the abstract Hecke algebra outside $N$ generated over $\Z$, resp. over $\cO$ by the standard Hecke
operators for all primes $\ell$ prime to $N$; for any such prime $\ell$, let 
$P_\ell(X)\in {\cal H}^N[X]$ be the minimal polynomial of the Hecke-Frobenius element (it is
monic, of degree $2^g$, see
\cite{CF} page 247).  Let
$\theta_{\pi}:{\cal H}^N(\cO)\rightarrow
\cO$ be the $\cO$-algebra homomorphism associated to $\pi_f$.

Let
$\hat{G}=GSpin_{2g+1}$ be the group of spinorial similitudes for the quadratic form
$$\sum_{i=1}^{g} 2x_ix_{2g+1-i}+x_{g+1}^2;$$ 
it is a split Chevalley group over $\Z[1/2]$ (we won't consider the prime $p=2$ in the sequel); it can be viewed as the dual
reductive group of
$G$ (see Sect.3.2 below); let
$\hat{B}$,
$\hat{N}$,
$\hat{T}$ the standard Borel, its unipotent radical, resp. standard maximal torus therein. The group $\hat{G}$ acts faithfully
irreducibly on a space
$V_{/\Z}$ of dimension $2^g$, via the spinorial representation. Let $B_V$ be the upper triangular
Borel of $GL_V$. Note that $\hat{B}$ is mapped into $B_V$ by the spin representation.

\vskip 1cm
\noindent{\bf 1.3.}
We put $\Gamma=Gal(\overline{\Q}/\Q)$. We assume that 

\vskip 5mm

{\bf (Gal)}  there exists a continuous homomorphism
$$\rho_\pi:\Gamma\rightarrow GL_V(\cO)$$ associated to $\pi$: that is, unramified outside $Np$,
and such that the characteristic polynomial of the (arithmetic) Frobenius element $\phi_q$ at a prime
$q$ not dividing $Np$ is equal to $\theta_\pi(P_q(X))$.

We shall make below an assumption on the reduction of $\rho_\pi$ modulo the maximal ideal of $\cO$ 
which will imply that $\rho_\pi$ act absolutely irreducibly on $V$ for each geometric fiber; hence
the choice of a stable $\cO$-lattice $V_\cO$ in $V\otimes K$ is unique up to homothety. 

Evidences for {\bf (Gal)} : For $g=2$, assuming

\vskip 5mm

{\bf (Hol)} $\pi_\infty$ is in the holomorphic discrete series,

\vskip 5mm

 Weissauer
\cite{We2} (see also
\cite{Har} and
\cite{La})  has shown the existence of a four-dimensional
$p$-adic Galois representation
$$\rho_\pi:\Gamma\rightarrow GL_V({\overline{\Q}_p})$$
 Moreover, his construction, relying on trace formulae, shows actually that
$$L(W_\pi,s)^4=L(\rho_\pi,s)^m$$
From this relation, one sees easily that the irreducibility of $\rho_\pi\otimes Id_{\overline{\Q}_p}$
implies that  the (Galois) semisimplification of 
$W_{\pi,p}$ is isomorphic to $n.\rho_\pi$ ($m=4n$).

Another crucial assumption for us will be that $p$ is prime to $N$ (hence $\pi$ is
unramified at $p$). Recall that under this assumption, Faltings has shown (Th.6.2 (iii) of
\cite{CF} and Th.5.6 of
\cite{Fa}) that for any
$q$, the
$p$-adic representation 
$H^q(S_U\otimes\overline{\Q}_p,V_\lambda(\overline{\Q}_p))$ is crystalline.

Let $D_p$, resp. $I_p$ be a decomposition subgroup, resp. inertia subgroup of $\Gamma$. Via
the identification $X^*(T)=X_*(\hat{T})$, we can view any $\mu\in X^*(T)$ as a cocharacter of
$\hat{T}$, hence as a homomorphism $I_p\rightarrow \Z_p^\times\rightarrow \hat{T}(\Z_p)\rightarrow
\hat{G}({\Z_p})$ where the first map is the cyclotomic character $\chi:I_p\rightarrow \Z_p^\times$.
We denote by $\tilde{\rho}$ the character of $T$ whose semisimple part is that of $\rho$, but whose central parameter is $d$.
it is the highest weight of an irreducible representation of
$G$ given by
$\rho$ on the derived group $G'$.
The character $\lambda+\tilde{\rho}$ has coordinates $(a_g+g,\ldots,a_1+1;\w)$. Let us introduce the assumption of Galois
ordinarity, denoted in the sequel {\bf (GO)} : 

\vskip 5mm

\begin{enumerate}
\item The image $\rho_\pi(D_p)$ of the decomposition group is contained in $\hat{G}$,
\item There exists $\hat{g}\in \hat{G}({\cO})$ such
that
$$\rho_\pi(D_p)\subset \hat{g}\cdot \hat{B}(\cO)\cdot \hat{g}^{-1},$$
\item the restriction of the conjugate $\rho_\pi^{\hat{g}}$ to $I_p$, followed by the quotient by the unipotent radical
$\hat{g}\cdot \hat{N}\cdot \hat{g}^{-1}$ of $\hat{g}\cdot \hat{B}\cdot \hat{g}^{-1}$ factors through
$-(\lambda+\tilde{\rho}):I_p\rightarrow \hat{T}(\Z_p)$.
\end{enumerate}

\noindent{\bf Comments:} 1) Let us introduce the condition of automorphic ordinarity:

{\bf (AO)} For each $r=1,\ldots,g$, 
$$v(\theta_\pi(T_{p,r}))=a_{r+1}+\ldots+a_g$$

where $T_{p,r}$ is the classical Hecke operator associated to the double class of 
$$diag(1_{r},p.1_{2g-2r},p^2.1_r).$$

We conjecture that for any $g$, if $\rho_\pi$ is residually absolutely irreducible, (AO) implies (GO). It is well-known for
$g=1$ (\cite{W1} Th.2.2.2, and
\cite{MT}). Moreover, for $g=2$, it follows from Proposition 7.1 of
\cite{TU} together with a recent result of E. Urban \cite{U2}.

2)  The minus sign
in front of $(\lambda+\tilde{\rho})$ comes from the definition of Hodge-Tate weights (for us: the jumps of the
Hodge filtration): the weight of the Tate representation $\Z_p(n)$ is $-n$.

\vskip 5mm

Let $\overline{\theta}_\pi=\theta_\pi\,{\rm mod.}\,\varpi$ and $\m=Ker\,\overline{\theta}_\pi$. Our
last assumption concerns ``{non-Eisenstein-ness}'' of the maximal ideal
$\m$.  It says that the image
of the residual representation
$\overline{\rho}_\pi$ induced by $\rho_\pi$ on $V_\cO/\varpi V_\cO$ is ``{large enough}''. More
precisely, let
$W_{\hat{G}}$ be the Weyl group of $\hat{G}$, viewed as a subgroup of $\hat{G}$. Recall the standard
description
$W_{\hat{G}}\cong S_g\propto
\{\pm 1\}^g$. Let $W'\subset \hat{G}$ corresponding to $\{\pm 1\}^g$. The ``{residually large image assumption}''
is as follows: 

\bigskip

{\bf (RLI)} : there exists a split (non necessarily connected) reductive
Chevalley subgroup
$H$ of
$\hat{G}_{/\Z}$  with
$W'\propto \hat{T}\subset H$, and a subfield $k'\subset k$, of order say
$\vert k'\vert=q'=p^{f'}$
($f'\geq  1$), so that $H(k')_\nu\subset Im\,\overline{\rho}_\pi$  and
$\overline{\rho}_\pi(I_p)\subset H^0(k')$.

Here, $H(k')_\nu$ denotes the subgroup of $H(k')$ consisting in elements whose
$\nu$ belongs to
$Im\,\nu\circ \bar{\rho}_\pi$.

\bigskip
 It has the consequence 
that $\overline{\rho}_\pi$ and ${\rho}_\pi$ are absolutely irreducible, hence
the uniqueness of the stable lattice
$V_\cO$ up to homothety.

\vskip 1cm
\noindent {\bf 1.4.}
One defines the sheaf $V_\lambda(\cO)$ over $S_U$ using the right action of
$U_p=G(\Z_p)$ (see
\cite{TU} Sect.2.1). We put
$V_\lambda(A)=V_\lambda(\cO)\otimes A$ for any
$\cO$-module
$A$; these are locally constant sheaves on $S_U$. Our main result is
as follows.

\begin{thm}  Let $\pi$ be cuspidal with $\pi_\infty$ in the discrete series and of good
level group
$U$, occuring in
$$H^d(S_U,V_\lambda(\C));$$ 
let
$p$ be a prime not dividing $N=level(U)$, assume {\bf (Gal)}, {\bf (GO)}, {\bf (RLI)}, $p>5$ and
that the weight $\lambda$ is small with respect to $p$:

$$ p-1>\vert \lambda+\rho\vert$$ 

Then, one has:
$$(i)\quad H^\bullet(S_U,V_\lambda(k))_\m=H^d(S_U,V_\lambda(k))_\m$$ 

$$(ii)\quad H^\bullet(S_U,V_\lambda(\cO))_\m=H^d(S_U,V_\lambda(\cO))_\m$$ 
and this $\cO$-module is free of finite rank. Similarly,
$$(iii)\quad H^\bullet(S_U,V_\lambda(K/\cO))_\m=H^d(S_U,V_\lambda(K/\cO))_\m$$
and this $\cO$-module is cofree of
finite rank.

The same statements hold for the cohomology with compact supports.
\end{thm}

{\bf Comments:} 

1) By standard arguments, the whole theorem follows if we show that: 

$$H_*^q(S_U,V_\lambda(k))[\m]=0\quad for \,q<d$$
where $*=c,\emptyset$, and for any Hecke-module $M$, $M[\m]$ stands for its $\m$-torsion.
This is the main result of the text.

2) In several instances in the proof, it is important that the maximal Hodge weights of the cohomology modules
involved are distinct for distinct modules, and are smaller than $p-1$; the condition

$$ p-1>a_1+\ldots+a_g+d$$
 
implies this; at the same time, it is also the condition needed to apply a 
comparison theorem of Faltings (Th.5.3 of
\cite{Fa}). We shall refer to this condition throughout the paper by saying that $\lambda$ is $p$-small. This terminology
has not the same meaning here than in \cite{PT}, but is in fact stronger than what is called $p$-smallness there.
Hence, under the present assumption, we can make use of Theorem D of \cite{PT}. In brief, this assumption is 
 unavoidable in our approach. The condition $p>5$ comes from the theory of modular
representations of reductive groups and has been pointed out to us by P. Polo. It is
necessary for the validity of Lemma 13 of Section 7.1, as there is a counterexample to this Lemma
for
$p=5$ and
$G=GSp(4)$; hence in our approach, the minimal possible $p$ is $7$ (for $g=2$ and $a_1=a_2=0$). 
Observe anyway that our bound on
$p$ depends only on $\lambda$ (not on the level group
$U$). This is crucial for the applications we have in view.   

3) The assumption {\bf (RLI)} is used only in Lemma 13 of Section 7.1, 
but this lemma is crucial for our proof of the Theorem.

4) Note that for $\lambda$ regular and for $g=2$, by
calculations of \cite{RT2}, and results of Schwermer and Franke (see Theorem 3.2 (i) of \cite{TU}),
 one has $H^q(S_U,V_\lambda(\C))=0$ for any $q<3$, while this is not so for the compact
support cohomology: the boundary long exact sequence for Borel-Serre compactification
relates
$H^2_c(S_U,V_\lambda(\C))$ to an $H^1$ of modular curves, which does not vanish. Our
vanishing statement concerns the localization at $\m$ and means that there is no mixing of Hodge
weights between the
$\m$-part of $H^2_c$ and that of $H^3_c$.

5) For $g=2$, E. Urban \cite{U1} has found a completely different proof of the absence of torsion of
$H^2(S_U,V_\lambda({\cal O})_\m$ under mild assumptions (with $\m$ non-Eisenstein). His proof is
much shorter than ours but relies on the fact that the complement in
$S_U$ of the Igusa divisor is affine, which is particular to the Siegel threefold. Whereas our
 theorem seems to carry over
(with the same proof) to various other situations, like the Hilbert (or quaternionic) modular case,
or unitary groups
$U(2,1)_{/\Q}$.
 
{\bf Evidences:} 1) If $g=2$ and $\pi$ is neither CAP nor endoscopic, one can conjecture that for
$p$ sufficiently general, $Im\,\rho_\pi$ contains the derived group $\hat{G}(\Z_p)$. Then {\bf (RLI)} is trivially
satisfied; if moreover
$p$ is also ordinary, the situation is as desired. Such a conjecture is unfortunately presently out of
reach. 

2) A more tractable situation is the following. See the details in Section 7.3. Let $F$ be a real
quadratic field with ${\rm Gal}(F/\Q)=\{1,\sigma\}$. Let $f$ be a holomorphic Hilbert cusp form for
$GL(2)_{/F}$,of weight
$(k_1,k_\sigma)$,
$k_1, k_\sigma\geq 2$, $k_1= k_\sigma+2m$ ($m\geq 1$). One can show (\cite{Y} and \cite{R1}) the
existence of a holomorphic theta lift from
$GL(2)_{/F}$ to
$G=GSp(4)_{/\Q}$ for $f$. Let $\pi$ be the corresponding automorphic representation of $G(\bA)$. It
is cohomological for a suitable coefficient system. Since
$f$ is not a base change from $GL(2)_{/\Q}$, $\pi$ is cuspidal, neither CAP nor endoscopic. We allow that $f$ is CM
of type
$(2,2)$; that is,  is a theta series coming from a CM quadratic extension $M=FE$ of
$F$, where
$E$ is imaginary quadratic. Moreover, $\pi$ is stable at $\infty$ (see \cite{R2}), $\rho_\pi$ exists and is motivic, namely:
$\rho_\pi=Ind_F^\Q\,\rho_f$, and it is absolutely irreducible. Moreover, for
$p$ sufficiently large (and splitting in $E$ in the $(2,2)$-CM case), the image of the associated Galois
representation
$\rho_\pi:\Gamma\rightarrow GL_{K}(V)$ is equal (up to explicit finite index) to the group of
points over a finite extension of $\Z_p$ of either the $L$-group ${}^L (Res^F_\Q GL(2)_{/F})={\rm Gal}(F/\Q)\propto (GL(2)
\times GL(2))^0$ (if $f$ is not CM), or
those of
${}^L Res^M_\Q{M}^\times=Gal(M/\Q)\propto (\G_m^2\times \G_m^2)^0$ if
$f$ is CM of type $(2,2)$. The subgroup
$H$ of
$\hat{G}$ whose image by the spin representation is
${}^L GL(2)_{/F}$ resp. ${}^L{M}^\times$,  does contain
$W'\propto \hat{T}$; that is, the assumption {\bf (RLI)} is satisfied for $H$. If $p$ is ordinary for
$f$ and splits in $F$, $\rho_\pi$ satisfies {\bf (GO)}; assume finally that $p$ satisfies
$p-1>k_1-1$; then, our result applies. See Sect.7.3 for numerical examples.      

\vfill\eject

In Section 8, we obtain a refinement of Theorem 1 as follows :

\begin{thm} Under the assumptions of Theorem 1, 

1) the finite free 
$\cO$-module $H^\bullet(S_U,V_\lambda(\cO))_\m$
coincides with the $\m$-localizations of 
\begin{itemize} 
\item the middle degree interior cohomology
$H_!^d(S_U,V_\lambda(\cO))=Im(H^d_c\rightarrow H^d)$,
\item the middle degree intersection cohomology
$IH^d(S_U,V_\lambda(\cO)).$
\end{itemize}

2) if $\lambda$ is regular, $H_!^d(\lU
S,V_\lambda(K))_\m$ contains only cuspidal eigenclasses, whose infinity type are in the discrete
series of HC parameter $\hat{\lambda}+\rho$.
\end{thm}

 The main tool for the proof of the first assertion is
the solution by Pink of a conjecture of Harder  
\cite{Pink}, together with a repeated use of our Theorem 1 for $GSp(2(g-r))$ for all integers
$r=1,\ldots,g$. To apply this argument, we need a mod. $p$ version of Kostant's formula,
proven in Theorem B of \cite{PT} under the assumption of
$p$-smallness. This is to apply Pink's theorem
in a fashion similar to \cite{Ha2} (who worked in characteristic zero). 
The second assertion follows by using a result of Wallach
\cite{Wall}. 

 We state in Section 9 and 10 several consequences of these results: 
\begin{itemize}
\item Control theorem and existence of $p$-ordinary cuspidal Hida families for $G$, improving upon
\cite{TU}, 
\item Verification of a condition of freeness of a cohomology module occuring in the definition of a
Taylor-Wiles system.
\end{itemize}

\vskip 1cm
\noindent {\bf 1.5.}
Let us briefly discuss the proof of Theorem 1. Let
$V_\lambda(\F_p)$ resp. $V_\lambda(k)$ be the etale sheaf over
$X\otimes\Q$ associated to the representation
$V_{\lambda\,/\F_p}$ of
$G_{\F_p}=G\otimes \F_p$, of highest weight $\lambda$, resp. its extension of scalars to $k$.
As mentioned in Comment 1) to Theorem 1, it is enough to show that 
$$(*)\quad W^j_*=H^j_*(X\otimes \overline{\Q},V_\lambda(k))[\m]=0$$
where
$*=\emptyset\, or \, c$, and for any $j<d$. 

Let $X_{/\Z[1/N]}$ be the moduli scheme classifying $g$-dimensional p.p.a.v. with level $U$
structure over
$\Z[1/N]$. Let $\overline{X}$ be a given toroidal compactification over $\Z[1/N]$ (see Th.6.7 of
Chapt.IV 
\cite{CF}, or Fujiwara \cite{Fu2}). Let $X_0=X\otimes \F_p$, $\overline{X}_0=\overline{X}\otimes \F_p$.

To the representation $V_{\lambda\,/\F_p}$ (with $\vert\lambda+\rho\vert<p-1$), one associates also a filtered log-crystal 
$\overline{\cV}_\lambda$ over $\overline{X}_0$ (see Section 5.2 below); the $F$-filtration  on the dual
$\overline{\cV}_\lambda^\vee$, satisfies $Fil^0=\overline{\cV}_\lambda^\vee$ and $Fil^{\vert\lambda\vert+1}=0$.  
Then, the main tools for proving $(*)$ are
\begin{itemize}
\item  Faltings's Comparison Theorem (\cite{Fa}, Th.5.3, see Sect.6.1). It says that, since
$p-1>\w$, for any $j\geq 0$, the linear dual of
$H^j_*(X\otimes
\overline{\Q}_p,V_\lambda(\F_p))$ is the image by the usual contravariant Fontaine-Lafaille functor
$\bV^*$ of the logarithmic de Rham cohomology
$$M=H^j_{log-dR,*}(\overline{X}\otimes
\F_p,\overline{\cV}_\lambda^\vee)=H^j(\overline{\cV}_\lambda^\vee\otimes
\Omega^\bullet_{\overline{X}_0}(log\,\infty)).$$

\item The mod. $p$ generalized Bernstein-Gelfand-Gelfand dual complex (section 5.4)
$$\kappa :\overline{\cal K}_\lambda^\bullet\hookrightarrow
\overline{\cV}_\lambda^\vee\otimes\Omega^\bullet_{\overline{X}_0}.$$ 
This is the mod. $p$ analogue of a construction carried in
Chapter VI of
\cite{CF}. The main result is that $\kappa$ is a
filtered quasi-isomorphism: it provides an explicit description of the jumps of the Hodge filtration in terms of
group-theoretic data. In particular for $j<d$, $\w$ is not a jump.
\item Lemma 13 in Section 7.1 shows, assuming {\bf (RLI)} and {\bf (GO)}, that if $W_j\neq 0$, its
restriction to the inertia group $I_p$ admits
$k\otimes \mu_p^{-\w}$ as subquotient.
\end{itemize}

Thus if $W_j\neq 0$ we obtain a contradiction since the maximal weight
 $\w$ should not occur in $W_j$. 

\vskip 4mm

 Theorem 2 is equivalent to the fact that the localization at $\m$ of the degree $d$ boundary cohomology of
$V_\lambda(k)$ vanishes. The argument for this is similar to the previous one, but makes use of the minimal compactification
$j:X_\Q\hookrightarrow X_\Q^*$ of $X_\Q=X\otimes \Q$ (instead of the toroidal one). The advantage of this compactification is
that Hecke correspondences extend naturally. We use crucially a theorem of R. Pink (Th.4.2.1 of
\cite{Pink}) which describes the Galois action on the cohomology of each stratum with coefficents
in the Žtale sheaves $R^qj_*V_\lambda(k)$; by the spectral sequence of the stratification it is
enough to show the vanishing of the localization at $\m$ of the degree $d$ cohomology of each
individual stratum. For this, we follow the same lines as for the proof of Theorem 1: the jumps of
the Hodge filtration in the degree $d$ cohomolology with compact support $H^d_c(X_r)$ of the
non-open strata $X_r$ cannot contain both $\w$ and $0$; on the other hand, if the
$\m$-torsion of $H^d_c(X_r)$ is not 0, Lemma 13 does imply that these weights both occur. Hence,
$H^d_c(X_r)_\m=0$. The last two sections contain two applications which were the original motivations for this
work.

\vfill\eject
{\bf 1.6.}
Acknowledgements: To start with, the authors have greatly benefitted of the seminar on toroidal
compactifications and cohomology of Siegel varieties held at the University of Paris-Nord in 97-98.
Without it, this work wouldn't have existed. They wish to thank the participants thereof, and in
particular, A. Abbs. Part of the writing was done by the second author at UCLA, at  MPI (Bonn) and
at Kyoto University; the excellent working conditions of these institutions were appreciated. A series of lectures on this work
at MRI (Allahabad) on the invitation of D. Prasad were also quite helpful. The first author would like to express his gratitude to the
Mathematic Department of Padova University,  and especially B. Chiarellotto for an invitation where a part of this work was exposed. 
During the preparation  of this text, we had
useful discussions or correspondence with many persons, in particular,  A. Abbs, D. Blasius, M. Dimitrov, A. Genestier, G. Harder, 
J.C. Jantzen, K. Khuri-Makdisi, K. Kuennemann, A. Nair,  B.C. Ngo, R. Pink, P. Polo, B. Roberts, J. Wildeshaus, H. Yoshida.  We thank
them heartily for their patience and good will.  
\vfill\eject
\tableofcontents

\vfill\eject
\section{Cohomology of Siegel varieties and automorphic representations}

We keep the notations of the introduction. Let us make precise what we mean by a good open compact subgroup of $G(\hat{\Z})$:
$U$ is good if 

1) it is neat: the subgroup of $\C^\times$ generated by the eigenvalues of elements in $U\cdot
G_\infty\cap G_\Q$ does not contain any root of unity other than $1$, and

2) 
$\nu(U)=\hat{\Z}^\times$. 

Let us now recall some properties of the cohomology groups
$H^\bullet_* (S_U,V_\lambda (K))$, for $K$ a field of characteristic zero and $*=\emptyset,c$ or
$!$ (as usual, $H^\bullet_!$ denotes the image of $H^\bullet_c$ in $H^\bullet$). In this section, $\g=Lie(G)$
will denote the real Lie algebra.

\subsection{Generalities over $\C$}
Let $U_\infty$ be the stabilizer in $G_\infty$ of the map
$$h:\C^\times\rightarrow G_\infty,\quad z=x+iy\mapsto \left(\begin{array}{cc} x\cdot1_g&y\cdot
s_g\\-y\cdot s_g&x\cdot 1_g\end{array}\right)$$
with $s_g$ the $g\times g$ antidiagonal matrix, with non-zero entries equal to $1$. 
For any good compact open subgroup $U\subset G(\hat{\Z})$, let
$$S_U=G(\Q)\backslash G(\A)/UU_\infty\quad {\rm and}\,  S=G(\Q)\backslash G(\A)/U_\infty$$
be the Siegel varieties of level $U$, resp. infinite level. Since $U$ has no torsion, $S_U$ is a smooth quasi-projective
algebraic variety of dimension $d={g(g+1)\over 2}$. $S$ is a pro-variety. For any (rational) 
irreducible representation $V_\lambda$ of $G$ of highest weight $\lambda$,
, we define the local system
$V_\lambda(\C)$ on $S_U$ as the locally constant sheaf of sections of
$$pr_1: G(\Q)\backslash \left( G(\A)\times V_\lambda\otimes \C\right)/UU_\infty\rightarrow S_U
$$
 By Prop.2.7 of
\cite{BW}(which does not require cocompactness), one has
$$H^\bullet(S_U,V_\lambda(\C))= H^\bullet(\g,U_\infty,\cC^\infty(G_\Q\backslash
G_\A,\C)\otimes V_\lambda(\C)).$$

The inclusions of spaces 
$$\cC_{cusp}^\infty(G_\Q\backslash
G_\A,\C)\subset \cC_c^\infty(G_\Q\backslash
G_\A,\C)\subset \cC^\infty_{(2)}(G_\Q\backslash
G_\A,\C) \subset \cC^\infty (G_\Q\backslash
G_\A,\C) $$
(where  $\cC^\infty_{cusp}=\cC^\infty_c \cap L^2_0$ and $\cC^\infty_{(2)}=\cC^\infty \cap L^2$)
 give rise to maps
$$ H^\bullet_{cusp}(S,V_\lambda(\C)) \rightarrow
H^\bullet_{c}(S,V_\lambda(\C))\rightarrow
H^\bullet_{(2)}(S,V_\lambda(\C))\rightarrow H^\bullet (S,V_\lambda(\C))$$ 
and a well-known theorem of Borel \cite{Bo1} asserts that their composition is injective:
$$H^\bullet_{cusp}(S,V_\lambda(\C)) \hookrightarrow H^\bullet_!(S,V_\lambda(\C)).$$
Moreover, as in the proof of Th.3.2 (or Th.5.2) of \cite{BW}, one has a $G_f$-equivariant
decomposition
$$H^\bullet_{cusp}(S,V_\lambda(\C))=H^\bullet(\g,U_\infty,\cC^\infty_{cusp}(G_\Q\backslash
G_\A,\C)\otimes
V_\lambda(\C))=$$
$$=\bigoplus_{\pi}\pi_f\otimes H^\bullet(\g,U_\infty,\pi_\infty^{U_\infty}\otimes
V_\lambda(\C))$$ 
where
$\pi=\pi_f\otimes \pi_\infty$ runs over the set of isomorphism classes of cuspidal
representations and $\pi_\infty^{U_\infty}$ is the Harish-Chandra module of
$\pi_\infty$.

\begin{pro} If $\lambda$ is regular dominant, the
interior, $L^2$ and cuspidal cohomology groups coincide and are concentrated in middle degree:
$$H^\bullet_{cusp}(S,V_\lambda(\C))=H^\bullet_{(2)}(S,V_\lambda(\C))=H^\bullet_!(S,V_\lambda(\C))=
$$
$$=H^d_!(S,V_\lambda(\C)).$$
\end{pro}

\noindent {\bf Proof:} Recall first that $H^\bullet_{cusp}=H^\bullet_{(2)}$ implies 
$H^\bullet_{cusp}=H^\bullet_{(2)}=H^\bullet_!(S,V_\lambda(\C))$ (see also Cor. to Th. 9
of \cite{Fa0}). 

By Th.4 of \cite{Bo2} (which applies here since
$rk\,G=rk\,U_\infty$): 
$$H^\bullet_{(2)}(S,V_\lambda(\C))=H^\bullet(\g,U_\infty,\cC^\infty_{(2)}(G_\Q\backslash
G_\A,\C)\otimes
V_\lambda(\C))=
$$
$$\bigoplus_{\pi}\pi_f\otimes H^\bullet(\g,U_\infty,\pi_\infty^{U_\infty}\otimes
V_\lambda(\C))$$ 
where $\pi$ runs over the discrete spectrum of $L^2(Z_\A G_\Q\backslash G_\A,\omega)$ where
 $\omega$ is the central
character of $V_\lambda^\vee$.

Let $\pi=\pi_f\otimes \pi_\infty$ be such an automorphic representation; its local components are
unitary. Moreover, one must have
$H^\bullet(\g,U_\infty,\pi_\infty^{U_\infty}\otimes V_\lambda(\C))\neq 0$ . 
By \cite{VZ} Th.5.6, the assumption that $\lambda$ is regular implies that  $\pi_\infty=A_\q(\lambda)$,
is a cohomological induction from
a parabolic subalgebra $\q$ which must be that of the Borel. In that case, this induction 
provides the discrete series. So, $\pi_\infty$ is one of the
$2^{g-1}$ unitary representations of $G_\infty$ in the discrete series of HC parameter
$\hat{\lambda}+\rho$. By
\cite{BW} Chapt.III, Cor.5.2 (iii), the tempered unitary $\pi_\infty$'s contribute only in
middle degree;  Moreover, since the automorphic representation $\pi=\pi_f\pi_\infty$ occurs in the global
discrete spectrum and admits at least one local component which is tempered,
it must be cuspidal; indeed, a theorem of Wallach (\cite{Wall}, Th.4.3) asserts that
if 
$\pi_\infty$ is tempered, the multiplicity of $\pi$ in 
$L^2_{disc}$ is equal to that in $L^2_0$.

\noindent{\bf Remark:} If $\lambda$ is not regular, there may also be non-tempered representations $\pi_\infty$
which occur as infinity type of $\pi$. However, by Langlands
classification (\cite{BW}, Sect.4.8, Th.4.11) and  Th.6.1 of \cite{BW},
 it implies that $H_{(2)}^q(S,V_\lambda)(\pi_f)\neq 0$ for some
$q<d$.  Franke's
spectral sequence (below) seems to suggest then that $H^q(S,V_\lambda)(\pi_f)\neq 0$ (we leave this as a question). 

This proposition will be used in the proof
of Theorem 2 (in Section 8 below) to rule out the occurence of non-cuspidal
representations in the localization
of the middle degree $L^2$-cohomology  $H^\bullet_{(2)}(S_U,V_\lambda)$, at a ``{
non-Eisenstein}'' maximal ideal of the Hecke algebra (that is, satisfying {\bf (RLI)}).  

\subsection{Franke's spectral sequence}

This section is not used in the sequel, but it provides a motivation for Section 8. 
By \cite{BW} Chap.VII
Cor.2.7, we have

$$H^\bullet(S,V_\lambda(\C))=H^\bullet(\g,U_\infty;
\cC^\infty(G(\Q)\backslash G(\A))\otimes V_\lambda(\C))$$

By \cite{Bo4}, one can replace the space of $\cC^\infty$-functions by those of uniformly moderate
growth. Franke has shown (\cite{Fr}, Th.13, or \cite{Walds} 2.2) that one can even replace this
space by the space $\cA(G)$ of automorphic forms on $G$. He has moreover defined a filtration
on$\cA(G)$, called the Franke filtration (see \cite{Walds} 4.7) whose graded pieces interpret as
$L^2$-cohomology. This yields an hypercohomology spectral sequence associated to
a filtered complex; more precisely:
  
 Let $\Phi^+$, resp. $\Phi_L^+$, be the positive
root system of $G$, resp. of a standard Levi $L$ of $G$, given by
$(G,B,T)$, resp. $(L,B\cap L,T)$. The corresponding simple roots are denoted by
$\Delta$, resp.
$\Delta_L$. For each standard parabolic
$P=L.U$, let $\a_P$ is the Lie algebra of the center of $L$.
Recall then Franke's spectral sequence (\cite{Fr} Th.19 or \cite{Walds} Corollaire 4.8)
$$E_1^{p,q}=H^{p+q}_{(2)}(S,V_\lambda(\C))\oplus_{{P}}\oplus_{w\in W^P(\lambda,p)}
{\rm
Ind}_{P_f}^{G_f}H^{p+q-\ell(w)}_{(2)}(S(L),V(L;w\cdot(\lambda+\rho))$$
$$\Rightarrow
H^{p+q}(S,V_\lambda(\C)) 
$$
where 
\begin{itemize}
\item $P=L\cdot U_P$ runs over the set of proper standard parabolic subgroups, 
\item $P_f$, resp. $G_f$ denotes the group of $\A_f$-points of $P$, resp. $G$, 
\item for each
$p$, $W^P(\lambda,p)$ is a certain subset of $W^L=\{w\in W;w^{-1}(\alpha)>0, \mbox{\rm for all}\,
\alpha\in
\Phi_L\}$, so that 
$W^L=\coprod_{p}W^P(\lambda,p)$,
\item the locally constant sheaf $V(L;w\cdot (\lambda+\rho))$ on the provariety 
$S(L)$
is attached to the representation of $L$
of highest weight $w\cdot (\lambda+\rho)=w(\lambda+\rho)-\rho$ (dominant for the order given
by $(L,B\cap L, T)$),  
twisted by $-w(\lambda+\rho)\vert_L$, that is, by the
one-dimensional representation of
$L$ attached to the (exponential of the) restriction of $-w(\lambda+\rho)$ to its (co-)center $\a_P$.
\end{itemize}
 This
spectral sequence is
$G_f$-equivariant. It allows one to represent any
$G_f$-irreducible constituent of
$H^{p+q}(S,V_\lambda(\C))$ as  ${\rm Ind}_{P_f}^{G_f}\pi_f$ where $\pi_f$ is an irreducible
admissible representation of $L_f$ such that $\pi=\pi_f\otimes\pi_\infty$ is automorphic, in the
discrete spectrum of $L^2(L_\Q Z_\A\backslash L_\A,\phi)$ with $P$ a rational parabolic in $G$,
$L$ its Levi quotient, and $\phi$ some unitary Hecke character.

Moreover, by Th.19 (ii) of \cite{Fr}, if $\lambda$ is regular, Franke's spectral sequence
 degenerates at
$E_1^{p,q}$. So, we have a Hecke-equivariant decomposition for each degree $q\in[0,2d]$:
$$\qquad H^{q}(S,V_\lambda(\C))=IH^{q}(S_U,V_\lambda(\C))\oplus
$$
$$\bigoplus_{{P}}\bigoplus_{p=0}^q \bigoplus_{w\in W^P(\lambda,p)}
IH^{q-\ell(w)}(S^L,V^L_{w(\lambda+\rho)-\rho}(\C)(-w(\lambda+\rho)_L)) 
$$

However, unlike the $GL_n$-case, the question of the rationality of this splitting for the group $G$
is open (with a possibly
negative answer). We nevertheless expect that it should yield, after localization at a ``{ non-Eisenstein}'' maximal prime
ideal of the Hecke algebra, an equality of the form
$$IH^{q}(S_U,V_\lambda(\C))_\m=H^{q}(S_U,V_\lambda(\C))_\m$$
for $\lambda$ regular. We establish this in Section 8 below for a suitable $\m$, by a Galois-theoretic argument which in some
sense replaces the lacking Jacquet-Shalika theorem.

\subsection{Hodge filtration in characteristic zero}

Recall we assumed that $U$ is good, so that its projection to any Levi quotient of $G$ is
torsion-free and 
$\nu(U)=\hat{\Z}^\times$. By the first condition, $S_U$ is smooth; the second condition implies that
$S_U$ admits a geometrically connected canonical model over
$\Q$. Let $X$ be this canonical model; it is a geometrically connected smooth quasi-projective scheme
over $\Q$. Let
$\overline{X}$ a toroidal compactification of $X$ defined by an admissible polyhedral cone
decomposition of $Sym^2 X^*(T)$ (\cite{AMRT} Chapt.3 and \cite{CF} Chapt.IV, Th. 5.7). Let
$\infty_X=\overline{X}-X$ be the divisor with normal crossing at infinity. Let $f:A\rightarrow X$ be
the universal principally polarized abelian variety with level $U$-structure over $X$ (it exists
over $\Q$).  Let
$Q$ be the Siegel parabolic of $G$, that is, the maximal parabolic associated to the longest simple
root for
$(G,B,T)$; let $M$ its Levi subgroup. For any
$B_M$-dominant weight
$\mu$, let
$\cW(\mu)$ resp.
$\overline{\cW}(\mu)$, be the corresponding automorphic vector bundle on
$X$, resp. its canonical Mumford extension to $\overline{X}$ 
 (see Th.4.2, Chap.VI). These are coherent sheaves. As observed by Harris \cite{Ha}, the coherent
cohomology 
$H^\bullet (\overline{X},\overline{\cW}(\mu))$ has a natural action of the Hecke algebra. Let
$\lambda=(a_g,\ldots,a_1;c)$ as above (recall that for simplicity we assume
$c=a_g+\ldots+a_1$). Let $H=diag(0,\ldots,0,-1,\ldots,-1)\in \g$.

\subsubsection{Complex Hodge Filtration}

It results from Deligne's mixed Hodge theory that the complex cohomology $H^m(X,V_\lambda(\C))$
carries a mixed Hodge structure with Hodge weights greater than, or equal to
$m+|\lambda|$ and that the interior cohomology (image of
$H^m_c\mapsto H^m$) is pure of Hodge weight $m+|\lambda|$. It is studied in greater details in
Sect.6.5 of
\cite{CF}. We won't need any information about its $W$-filtration, so we concentrate on its
$F$-filtration (Hodge filtration). With the notation of 6.4 of \cite{CF}, de Rham comparison theorem 
reads:
$$H^m(X(\C),V_\lambda(\C))=H^m(X(\C),\cV_\lambda^\vee)$$
where $\cV_\lambda$ denotes the coherent sheaf associated to the $Q$-representation restriction to the
Siegel parabolic $Q$ of the $G$-representation of highest weight
$\lambda$. The reason for the dual (denoted
${}^\vee$) is the following. The de Rham comparison theorem sends the local system $R^1f_*\C$ on
$R^1f_*\Omega_{A/X}^\bullet$;
however, as  explained on top
of page 224 of
\cite{CF}, the construction of coherent sheaves from $Q$-representations associates to the standard
representation the dual of
$R^1f_*\Omega_{A/X}^\bullet$, while the locally constant sheaf associated to the standard
representation is $R^1f_*\C$. 

  Let
$\g$, resp. $\t$, be the Lie algebra of $G$, resp. $T$.  Let 
$$H=diag(0,\ldots,0,-1,\ldots,-1)\in \t$$
Let 
$W^M$ be the set of Kostant representatives of the quotient
$W_M\backslash W_G$ of the Weyl groups; for each $w\in W^M$, let
$p(w)=-\left(w(\lambda+\rho)-\rho\right)(H)$; it is a non-negative integer. The main result of
Sect.6.5 (Theorem 5.5 (i), Chapt.VI) of
\cite{CF} gives a Hecke-equivariant description of the graded pieces of the $F$-filtration in terms
of coherent cohomology of automorphic vector bundles extended to a toroidal compactification
$\overline{X}$ of $X$, as follows:   
$$ (BGG)\quad gr_F^p H^\bullet(X,V_\lambda(\C))=\bigoplus_{w\in W^M,
p(w)=p}H^{\bullet-\ell(w)}(\overline{X},\overline{\cW} (w(\lambda+\rho)-\rho)^\vee)
$$

Because of our comment on de Rham comparison theorem, we see that contrary to what is mentioned in  R.
Taylor's paper (\cite{RT2} p.295, l.14 from bottom), the
statement of Th.5.5, l.6 in
\cite{CF} is correct, because the local system denoted $V_\lambda$
 in Faltings-Chai is actually dual to the one denoted $V_\lambda$ in Taylor's and in the present
paper. Our statement, in accordance to Faltings', is that the sum runs over the
$w$ such that $w(\lambda+\rho)(H)+p=\rho(H)$. We think therefore that Taylor's statement cited above
is incorrect (but correct after a Tate twist, anyway).

For any subset
$B$ of
$A=\{1,\ldots,g\}$, let
$(B,\overline{B})$ the corresponding partition of
$A$ and 
$w_B\in W_G$ such that for $(t;\nu)\in T$,
$w_B(t;\nu)=(t_B^{-1},t_{\overline{B}};\nu)$. An easy calculation
shows that for any $w\in W_G$, if $w=(\sigma,w_B)$ for some permutation $\sigma$ of $A$ and $B$ some
subset of $A$, one has:
$$-(w(\lambda+\rho)-\rho)(H)=-(w_B(\lambda+\rho)-\rho)(H)= \sum_{i\in B}(a_i+i)$$

We put $j_B=\sum_{i\in B}(a_i+i)$, so $j_A=\w$ is the motivic weight defined in the introduction. The $j_B$'s belong to
the closed interval $[0,\w]$. They are indexed by a set of cardinality $2^g$, but need not be mutually
distinct, from $g=3$ on. Note that for any degree $m$ of the cohomology, the jumps of the Hodge filtration 
occuring in
$H^m$ always form a subset of $\{j_B;B\subset A\}$.

Let $\pi=\pi_f\otimes\pi_\infty$ be a cuspidal representation of $G(\A)$, with $\pi_\infty$
holomorphic in the discrete series of HC parameter $\hat{\lambda}+\rho$; let
$\theta_\pi:\cH^N\rightarrow \C$ be the character of the (prime-to-$N$)  Hecke algebra, 
associated to $\pi$ and $\p_\pi={\rm Ker}\,\theta_\pi$. By \cite{BW}
Chapt.III Th.3.3 (ii), the $(\g,U_\infty)$-cohomology of $\pi_\infty\otimes V_\lambda$ 
is concentrated
in degree $d$. we put 
$$W_\pi=H^d(X,V_\lambda(\C))[\p_\pi]$$
By cuspidality of $\pi$, $W_\pi$ has a Hodge structure pure of weight $\w=d+|\lambda|$:
$$W_\pi=\bigoplus_{p+q=\w}W_\pi^{p,q}$$
Let us show that $ W_\pi^{\w,0}$ and $ W_\pi^{0,\w}$ are both non-zero. 
More precisely, let $w'\in W^M$ be the Kostant representative of largest length, namely $d$ (it is
unique, and if $w''\in W_M$ is the unique element of largest
length, then $w'w''$ is the unique element of largest length in $W_G$). Then, 

\begin{pro} There is a
$\cH^N$-linear embedding

$$\pi_f^U\subset H^{\w,0}=H^0(\overline{X},\overline{\cW}_{w'(\lambda+\rho)-\rho}),\quad
\pi_f^U\subset H^{0,\w}=H^d(\overline{X},\overline{\cW}_\lambda)$$
\end{pro}

{\bf Proof:}  Let $\q$ be the Lie algebra of the Siegel
 parabolic. Since
$\pi$ is cuspidal, a calculation of M. Harris, Prop.3.6 of \cite{Ha} shows that for
any $q$ and $\mu$ $M$-dominant, $\pi_f^U\otimes H^q(\q,U_\infty,\pi_\infty\otimes W_\mu)$ embeds
$\cH^N$-linearly into
$H^q(\overline{X},\overline{\cW}_\mu)$. Moreover by Theorem 3.2.1 of \cite{BHR}, 
$H^q(\q,U_\infty,\pi_\infty\otimes W_\mu)$ does not vanish in only two cases:
$\mu=\lambda$ and $q=d$, or $\mu=w'(\lambda+\rho)-\rho$ and $q=0$. 

{\bf Remark:} If $\pi$ is stable at infinity, that is, if all the possible infinity types $\pi'_\infty$ in the discrete series
of HC parameter $\hat{\lambda}+\rho$ give rise to automorphic cuspidal representations $\pi'=\pi_f\otimes \pi'_\infty$, then
all the possible Hodge weights do occur in $W_\pi$:

$$\mbox{\rm For any}\quad j_B, B\subset A,\, A=B\coprod \overline{B}\quad W_\pi^{j_B,j_{\overline{B}}}
\neq 0.
$$

\subsubsection{$p$-adic Hodge filtration}

We view now $\theta_\pi$ as taking values in a $p$-adic field $K\subset \overline{\Q}_p$. Let
$\p_\pi$ be the prime ideal of $\cH^N(K)=\cH^N\otimes K$ defined as kernel
of $\theta_\pi\otimes\,Id_K$.

The Hodge-to-de Rham spectral sequence

$$(BGG)_\Q\qquad E_1^{p,q}=\bigoplus_{w\in
W^M,p(w)=p}H^{p+q-\ell(w)}(\overline{X},\overline{\cW} (w(\lambda+\rho)-\rho))$$
$$
\Rightarrow H^{p+q}(\overline{X},\overline{\cV}_\lambda\otimes \Omega^\bullet_{\overline{X}/\Q}(log\,\infty_X))$$

\noindent makes sense over $\Q$ and degenerates in $E_1^{p,q}$ (\cite{CF} Sect.VI.6, middle of page 238). Here, $\cV_\lambda$
denotes the flat vector bundle defined over $\Q$  associated to the rational representation
$V_\lambda$ of $G$. More explanations on the rational structures involved, as well as integral
versions thereof will be given in Section 4.2.

Actually, let $C$ be the completion of an algebraic closure of $\Q_p$; by Th. 6.2 of \cite{CF}, there is a Hodge-Tate
decomposition theorem inducing the splitting of
$(BGG)_{\C}$;  More precisely:

$$(BGG)_{HT}\qquad H^{p+q}(X,V_\lambda(\Q_p))\otimes C\cong
$$
$$\bigoplus_{w\in
W^M,p(w)=p}H^{p+q-\ell(w)}(\overline{X},\overline{\cW} (w(\lambda+\rho)-\rho))\otimes C(p(w))
 $$

 By a theorem of Harris \cite{BHR}, the Hecke algebra $\cH^N$ acts naturally on each
summand of the LHS of this splitting. Now, the 
main feature of the above splitting is its naturality for algebraic correspondences on
$\overline{X}$.
It implies the compatibility of the decomposition $(BGG)_{HT}$ with the action of $\cH^N$.  
Let $W_{\pi,p}=W_\pi\otimes K$; it is cut by algebraic correspondences of the Hecke algebra $\cH^N(K)=\cH^N\otimes K$ with
coefficients in $K$. The restriction of $(BGG)_{HT}$ to the part killed by
$\p_\pi$ is still a $\cH^N$-equivariant decomposition of
$W_{\pi,p}\otimes_K C$. If we assume {\bf (Hol)}, we see
from Prop.1 above that the Hodge-Tate weights $\w$ and $0$ do occur; indeed,

$$W_{\pi,p}^{\w,0}=H^0(\overline{X},\overline{\cW}_{w'(\lambda+\rho)-\rho})[\p_\pi]$$
and 
$$W_{\pi,p}^{0,\w}=H^0(\overline{X},\overline{\cW}_\lambda)[\p_\pi]$$
by  comparing to complex cohomology, we see from Prop.1 that these two spaces are non-zero.

Let us remark that if $\pi$ is stable at infinity, the analogue of Prop.2 for all possible infinity
types in the discrete series of HC parameter $\hat{\lambda}+\rho$ (in number $2^g$, but isomorphic
two by two) implies that all the possible Hodge-Tate weights $j_B$ ($B\subset A$) do occur in the
Hodge-Tate decomposition of $W_{\pi,p}$.

\vfill\eject

\section{Galois representations}

\subsection{Relation between $\rho_\pi$ and $W_{\pi,p}$}

The absolute Galois group $\Gamma$ acts on  $W_{\pi,p}$. Let us first recall, for later use, the following well-known fact.

\begin{lem} $W_{\pi,p}$ is pure of weight $\w$. That is, for any $\ell$ prime to $Np$, all the eigenvalues ofthe geometric
Frobenius at $\ell$ have archimedean absolute value $\ell^{\w/2}$.
\end{lem}

{\bf Proof:} Since $\pi$ is cuspidal, we know by a theorem of Borel (see Sect.2.1) that
$W_{\pi,p}$ is contained in the interior cohomology  $H^d_!(X,V_\lambda)$. By Th.1.1 of Chap.VI of \cite{CF}, there is a
toroidal compactification $Y\subset \overline{Y}$ of the $\vert\lambda\vert$-times fiber product $Y=A^{\vert\lambda\vert}$  
of the universal abelian variety
$A$ above a toroidal compactification of the Siegel variety $X\subset \overline{X}$, all these schemes being flat over
$\Z[1/N]$; over this base,
$\overline{Y}$ is smooth and $\overline{Y}-Y$ is a divisor with normal crossings. One can interpret the etale sheaf as cut 
by algebraic correspondences in $(R^1\pi_*\Q_p)^{\otimes d}$ (see \cite{CF} p.235, and 238, or this
text, Sect.5.2), hence
$H^d_*(X,V_\lambda)\subset H^{\w}_*(Y,\Q_p)$ ($*=\emptyset, c$).  By the classical commutative diagram
$$\begin{array}{ccc} H^w_c(Y,\Q_p)&\rightarrow H^w(\overline{Y},\Q_p)\rightarrow&
H^w(Y,\Q_p)\\\bigcup&&\bigcup\\H^d_c(X,V_\lambda)&\longrightarrow &H^d(X,V_\lambda)\end{array}$$
We conclude that $H^d_!(X,V_\lambda)$ is pure of weight $\w$; recall that this can be interpreted either
in the sense of Deligne (take $\ell$ unramified
and different from $p$, then the eigenvalues of geometric $Fr_\ell$ have archimedean absolute values $\ell^{\w/2}$)
or in a $p$-adic sense (in the crystalline case, say: that the eigenvalues of the crystalline Frobenius
have archimedean absolute values $p^{\w/2}$).

Assume now
that
$\pi$ admits an associated
$p$-adic Galois representation
$\rho_\pi:\Gamma\rightarrow GL_V(\overline{\Q}_p)$; we assume that $\rho_\pi$ is irreducible.  We don't know
a priori whether
$\rho_\pi$ is a Galois constituent of
$W_{\pi,p}$ although,
by \cite{CF} Chapter VII Th.6.2, we know that the characteristic polynomial of
$\rho_\pi$ annihilates the global $p$-adic representation $W_{\pi,p}$. If moreover $p$ does not
divide $N$, we know by Faltings \cite{Fa} Th.5.2 that $W_{\pi,p}$ is
crystalline but we don't know this a priori for
$\rho_\pi$. However, for $g\leq 2$, if $\rho_\pi$ is absolutely irreducible, we do know that
it is a constituent of
$W_{\pi,p}$  (by \cite{RT2} or
\cite{We2}). Indeed,  for $g=2$, Weissauer
\cite{We2} ( completing works of
\cite{Har}, \cite{RT2} and
\cite{La})  has shown the existence of a four-dimensional
$p$-adic Galois representation
$$\rho_\pi:\Gamma\rightarrow GL_V({\overline{\Q}_p})$$
such that
$$L(W_\pi,s)^4=L(\rho_\pi,s)^m$$
thus, the assumption of irreducibility for $\rho_\pi$ implies that the Galois
semisimplification $W_{\pi,p}^{s.s.}$ of $W_{\pi,p}$ is isomorphic to $n.\rho_\pi$, for
$m=4n$. In turn, it also implies that $\rho_\pi$ is pure of weight $\w$ and is crystalline at $p$ if $p$ is prime to $N$.

There are other situations, namely when $\pi$ is a base change of a Hilbert modular eigenform, where
one knows that
$\rho_\pi$ is crystalline, although one may not know that it is a constituent of $W_{\pi,p}$; see
Sect.7.3 below.  One of
the uses of our assumption {\bf (RLI)} will be to relate (residually only)
$W_{\pi,p}$ and
$\rho_\pi$ (see Sect. 7.1).

\subsection{Spin groups and duality}

\subsubsection{description}

For the general definitions on spinors, we follow \cite{FH} Sect.20.2, and \cite{DeH} VIII.8 and IX.2; however 
by lack of references for our precise need, we give some details in this section.  
Although these groups exist over $\Z$, we'll restrict ourselves to $\Z[1/2]$, as $p=2$ is excluded of our study. Let
$\tilde{V}=\A_{\Z[1/2]}^{2g+1}$  endowed with the quadratic form
$q(x)=\sum_{i=1}^{g} 2x_ix'_i+x_0^2$ for $x=\sum_{i=1}^gx_ie_i+x_0e_0+\sum_{i=1}^g x'_ie'_i$; the
scalar product is denoted by $\langle x, y\rangle$. The canonical basis is ordered as
$(e_g,\ldots,e_1,e_0,e'_1,\ldots,e'_g)$, so that $\langle e_i,e'_j\rangle=\delta_{i,j}$, $e_0$ is
unitary, $W=\langle e_g,\ldots,e_1\rangle$ and $W'=\langle e'_1,\ldots,e'_g\rangle$ are totally
isotropic, and the sum
$\tilde{V}=W\oplus W'\oplus
\langle e_0\rangle$ is orthogonal. The Clifford algebra $C(\tilde{V},q)$ over $\Z[1/2]$ is the
quotient of the tensor algebra by the two-sided ideal generated by $x\otimes x-q(x)$, ($x\in
\tilde{V})$; it is $\Z/2\Z$-graded: $C(\tilde{V},q)=C^+\oplus C^-$.  The main
involutive automorphism $\Pi$ is defined as $Id$ on $C^+$ and $-Id$ on $C^-$; the main antiinvolution
$x\mapsto x^*$ is defined by
$v_1\cdot \ldots\cdot v_r\mapsto (-1)^r v_r\cdot
\ldots\cdot v_1$. We write $N(x)=x\cdot x^*=x^*\cdot x$ for the spinor norm. The
$\Z[1/2]$-group scheme
$GSpin_{\tilde{V}}=GSpin_{2g+1}$ (called the regular Clifford group in \cite{DeH} IX.2) is defined as
the group of invertible elements
$g$ of
$C(V,q)$ such that
$g\cdot \tilde{V}\cdot g^{-1}= \tilde{V}$. The group of orthogonal similitudes
$GO_{\tilde{V}}=GO_{2g+1}$ is defined as the group of $h\in GL_{\tilde{V}}$ such that
$q\circ h=c(h)\cdot q$. Consider the group-scheme morphism  
$$\nu:GO_{2g+1}\rightarrow \G_m,\quad, h\mapsto det\,h\cdot c(h)^{-g}.$$ 
One has $c(h)=\nu^2(h)$.  Moreover,
the homomorphism of
$\Z[1/2]$-group schemes 

$$\psi: GSpin_{\tilde{V}}\rightarrow GO_{\tilde{V}}, \quad g\mapsto (x\mapsto \Pi(g)\cdot x\cdot
g^*)$$ 
is an isogeny of degree two (using \cite{DeH} VIII.8) which satisfies $\nu\circ \psi=N$.
The spin representation $\sp$ is a representation of $GSpin_{\tilde{V}}$ on $V=\bigwedge W$; it can be defined via
the universal property of the Clifford algebra, as in
\cite{FH} Lemmata 20.9 and 20.16. We have ${\rm dim}\,V=2^g$. 
We write $\hat{G}$ for $GSpin_{\tilde{V}}$. It is a Chevalley group over
$\Z[1/2]$; the standard maximal torus $\hat{T}$, resp. Borel $\hat{B}$, of $\hat{G}$ is the inverse image by $\psi$ of the
diagonal torus, resp. upper triangular subgroup in $GO_{2g+1}$.

\subsubsection{Dual root data}
 
We want to recall first the notion of a (reduced) based root
datum 
$$(M,R,\Delta,M^*,R^\vee,\Delta^\vee),$$
 consisting of two free $\Z$-modules $M$, $M^*$ of rank, say,
$n$ with a perfect pairing $M\times M^*\rightarrow \Z$ and finite subsets $R\supset \Delta$ in $M$,
resp. $R^\vee\supset \Delta^\vee$ of $M^*$, together with a bijection $R\rightarrow R^\vee$; $R$ is
the set of roots, and $\Delta$ the simple roots; these data should satisfy two conditions RD I and RD
II: cf
\cite{Sp} 1.9 or rather, for the degree of generality that we need, Exp.XXI Sect.1.1 and 2.1.3; here,  ``{reduced}'' means
that in the set of roots $R$, we allow no multiple of any given root except its opposite. 

In order to make some calculations, let us
recall briefly the classification given by these data. The main reference is \cite{DGr}, whose ExposŽs are quoted by their
roman numbering. 

\begin{de} For any scheme
$S\neq
\emptyset$, a split reductive group with ``{Žpinglage}'' over $S$, is a t-uple $(G,B,T, (X_\alpha)_{\alpha\in \Delta})_S$
consisting in a connected reductive group scheme
$G_S$ of rank $n$, together with a Borel $B_S$ and split maximal torus $T_S\subset B_S$ : $T\cong \G_m^n$.
Let $R$, resp.$\Delta\subset R$, be the root system, resp. set of simple roots, attached to $(G,B,T)$  (Exp.XIX Sect.3). The
``{Žpinglage}'' 
$(X_\alpha)_{\alpha\in
\Delta}$ is the datum for each $\alpha\in \Delta$, of a section $X_\alpha\in \Gamma(S,\g_\alpha)$
which is a basis of $\g_\alpha$ at each point $s\in S$.
\end{de}
 
For details on ``{Žpinglages}'', see 
\cite{DGr} XXII 1.13 and XXIII 1.1. 
Any such split reductive group defines
a reduced based root datum
$$(M,R,\Delta,M^*,R^\vee\,\Delta^\vee)$$
Note that the ``{Žpinglage}'' is not needed in the construction, it comes 
in only for the fidelity of the functor. The definition runs as follows
 put $M=X^*(T)$, $M^*=X_*(T)$; the duality $\langle\,,\,\rangle$ between these modules is the
composition $(\lambda,\mu)\mapsto \lambda\circ\mu$, 
$R$, resp.$\Delta$ is the set of roots, resp. simple roots attached to
$(G,B,T)$, and $\alpha^\vee$ is defined for each $\alpha\in \Delta$ as follows: let $T_\alpha$ be the connected component of
$Ker\,\alpha$, let $Z_\alpha$ be its centralizer in $G$. It is reductive of semisimple rank one, hence its derived group
$Z'_\alpha$ is isomorphic to
$SL(2)$ or
$PGL(2)$, and its character group is generated by $\alpha$; then, $\alpha^\vee:\G_m\rightarrow Z'_\alpha\cap T$ is defined as
the unique cocharacter of $Z'_\alpha$ such that $\alpha\circ\alpha^\vee=2$. For details, see Exp.XX, Th.2.1.   
As checked in Exp.XXII 1.13, these data satisfy the two conditions (DR I) and (DR II) of Exp.XXI 1.1, hence do form a
based root datum (donnŽes radicielles ŽpinglŽes). The system thus obtained is reduced.

\begin{thm}
There is an equivalence of categories between reduced based root data and split reductive groups with ``{Žpinglage}''.
\end{thm}
This is the main theorem of \cite{DGr}, it consists in 4.1 of Exp.XXIII Sect.4 and Th.1.1 of Exp.XXV Sect.1.

Now, given a reduced based root datum, one can form its dual by exchanging $(M,R,\Delta)$ and $(M^*,R^\vee,\Delta^\vee)$. This
induces a duality of split reductive group schemes with Žpinglages, over a base $S$.  Let us apply this to our situation. We
take
$G=GSp_{2g}$, $(G,B,T)_{/\Z[1/2]}$;
$M=X^*(T)$ and
$M^*=X_*(T)$, naturally paired by the composition. By using the standard basis of
$X^*(T)$, one identifies $M$ to the subgroup of $\Z^g\times \Z$, consisting in $\mu=(\mu_{ss};\mu_c)$ such that
$\vert\mu\vert\equiv\mu_c\,mod.2$. This lattice is endowed with the standard scalar product; here
$\Z^g$ corresponds to the characters of the semisimple part of $T$, and the last component to the central variable. In this
identification, $R\subset \Z^g\times\{0\}$ and one can write
$\alpha^\vee=2\cdot{\alpha\over\langle\alpha,\alpha\rangle}$ in the space $\Q^g\times\{0\}$. The simple roots of
$G$ are
$\alpha_g=t_g/t_{g-1},\ldots,\alpha_1=t_1^2\nu^{-1}$, for $t=diag(t_g,\ldots,t_1,t_1\nu^{-1},\ldots,t_g\nu^{-1})\in T$; hence
their coordinates in $M=\Z^g\times\Z$ are $(1,-1,0,\ldots;0)$,...,$(0,\ldots,2;0)$. The corresponding coroots have therefore
coordinates
$\alpha_g^\vee=(1,-1,\ldots;0)$,..., $\alpha_1^\vee=(0,\ldots,1;0)$.
 Then, $X_*(T)$ is identified to $\Z^g\times\Z+{1\over 2}\cdot diag(\Z^{g+1})$. 

 The resulting dual of $(G,B,T)_{\Z[1/2]}$ is precisely $(\hat{G},\hat{B},\hat{T})_{\Z[1/2]}$ 
 (it is true as well over
$\Z$, but we don't need, and don't want to consider characteristic $2$ spin groups).

Let $\hat{\varpi}$ be the minuscule weight of $\hat{G}$; it belongs to $X^*(\hat{T})=X_*(T)$. It satisfies the formulae:
$\langle\hat{\varpi},\alpha_i^{\vee\vee}\rangle=\delta_{1,i}$ for $i=1,\ldots, g$. Hence, in the basis we have fixed, its
coordinates are $(1/2,\ldots,1/2;x)$. 
The central parameter $x$ must equal $1/2$ as well, because the homomorphism $\psi$ is etale of degree two, and induces the
standard representation, whose highest weight is therefore $2\hat{\varpi}$, but whose central character is $z\mapsto z$.
Now, any character $\mu\in X^*(T)$ is identified to a cocharacter of $\hat{T}$. Then,

\begin{lem} In $X^*(\G_m)=\Z$, for any $\mu=(\mu_{ss};\mu_c)\in X^*(T)$, one has: 
$$(3.2.2.1)\quad \hat{\varpi}\circ \mu={\vert \mu_{ss}\vert\over 2}+{\mu_c\over 2}$$
\end{lem}

Note that the right-hand side is an integer.

\noindent{\bf Proof:} Clear. 

\noindent Let us make simple remarks:

1) Let $B_V$ be the upper triangular
Borel of $GL_V$. Then $\hat{B}$ is mapped into $B_V$ by the spin representation.

2) In the identification $X_*(T)=X^*(\hat{T})$, the central cocharacter
$z:\G_m\rightarrow T$ becomes the multiplier $N:\hat{T}\rightarrow \G_m$ of our regular Clifford group $\hat{G}$; it is clear
on the level of tangent maps. Dually, via  $X_*(T)=X^*(\hat{T})$, the multiplier $\nu:G\rightarrow \G_m$ becomes the central
cocharacter $\G_m\rightarrow \hat{T}$.

3) If we describe $T_{GO_{\tilde{V}}}(\C)$ as the set $\G_m\times T_{O_{\tilde{V}}}$ of matrices 
$$diag(z\cdot t_g,\ldots,z\cdot t_1, z,z\cdot
t_1^{-1},\ldots, z\cdot t_g^{-1})$$
 then, $\hat{T}(\C)$ can be described as the set of t-uples
$(t_g,\ldots,t_1,[u,\zeta])$ where
$u^2=t_g\cdot\ldots\cdot t_1$ and $\zeta^2=z$, the couple $(u,\zeta)$ being taken modulo the group generated by $(-1,-1)$.
The map $\psi: \hat{T}(\C)\rightarrow T_{GO}(\C)$ is then given by $t_i\mapsto t_i$, $[u,\zeta]\mapsto \zeta^2$. All this
follows easily from the fact that
$\psi$ is dual of the degree two isogeny $T_{ss}\times Z_G\rightarrow T$ given by $(t_{ss},z)\mapsto t_{ss}\cdot
z$.

 Let us apply these considerations to compute the local
Langlands correspondence for a representation
$\pi_p$ of
$G(\Q_p)$ in the principal series. Let us assume $\pi_p=Ind_{B(\Q_p)}^{G(\Q_p)}\phi$ (unitary induction). If $\phi$ is
unramified, it can be viewed as
$(3.2.2.2)\quad \phi=(\alpha_g,\ldots,\alpha_1;\gamma)\in
\C^g\times\C$, the parametrization being given by: 
$$diag(t_g,\ldots,t_1,\nu\cdot t_1^{-1},\ldots,\nu\cdot t_g^{-1})\mapsto \vert
t_g\vert_p^{\alpha_g}\ldots\vert t_1\vert_p^{\alpha_1}\vert \nu\vert_p^{(\gamma-\alpha_g-\ldots \alpha_1)/2}$$
 
Even if it is ramified, we can make the following identifications

$$(3.2.2.3)\quad \underline{\rm Hom}(T(\Q_p),\C^\times)=\underline{\rm Hom}(X_*(T)\otimes
\Q_p^\times,\C^\times)=$$

$${\rm Hom}(X_*(T), \underline{\rm Hom}(\Q_p^\times,\C^\times))=
X^*(T)\otimes \underline{\rm Hom}(\Q_p^\times,\C^\times)=
$$
$$\underline{\rm
Hom}(\Q_p^\times,\C^\times\otimes  X^*(T))=\underline{\rm Hom}(\Q_p^\times,\hat{T}(\C)).$$

So that we can view $\phi$ as a cocharacter $\Q_p^\times\rightarrow \hat{T}(\C)$. We introduce a twist 
of this character by $d$ on the central component ($\gamma\mapsto \gamma-d$), in order to get rid of the irrationality
inherent to Langlands parameters:
$\tilde{\phi}=\phi\cdot \vert\nu\vert_p^{-d}$, it corresponds to the cocharacter $\tilde{\phi}$ obtained by twisting $\phi$
by the unramified cocharacter
$\G_m\rightarrow Z_{\hat{G}}(\C),t\mapsto \vert t\vert_p^{-d}$.In the unramified case, $\tilde{\phi}$ is given by the
formula

$$(3.2.2.4)\quad t\mapsto (\vert t\vert_p^{\alpha_g},\ldots,\vert t\vert_p^{\alpha_1},[\vert t
\vert_p^{{\alpha_g+\ldots+\alpha_1\over 2}},\vert t\vert_p^{(\gamma-d)/2}]).$$

 Consider the canonical map
$a:W_{\Q_p}\rightarrow
\Q_p^\times$ given by class-field theory (sending arithmetic Frobenius to
$p$). The composition $\tilde{\phi}\circ a$ is denoted $\sigma(\pi_p)$ and is called the image by Langlands local
correspondence of
$\pi_p$.

\vskip 5mm

Let us return now to our Galois representations. Note first that
the question whether $\rho_\pi$, if absolutely irreducible, factors through the spin representation 
$$\hat{G}(\overline{\Q}_p)\hookrightarrow GL_V(\overline{\Q}_p)$$
is open.

However, for $g=2$, if $\pi$ is stable at $\infty$ and if $\pi$ satisfies multiplicity
one: 
$m(\pi)=1$, then it can be shown that
$\rho_\pi$ takes values in $\hat{G}$ (see \cite{RT2} p.295-296). This remark, due to E. Urban (to appear) results from
Poincar\'e duality and the autoduality of
$\pi$ (which is well known, at least, at almost all places).

\subsection{Ordinarity}

Let $D_p$, resp. $I_p$ be a decomposition subgroup, resp. inertia subgroup of $\Gamma$. Via
the identification $X^*(T)=X_*(\hat{T})$, we can view any $\mu\in X^*(T)$ as a cocharacter of
$\hat{T}$, hence as a homomorphism $I_p\rightarrow \Z_p^\times\rightarrow \hat{T}(\Z_p)\rightarrow
GL_{\Z_p}(V)$ where the first map is the cyclotomic character $\chi:I_p\rightarrow \Z_p^\times$. Let
$\tilde{\rho}=(g,\ldots,1;d)$. Thus, $\tilde{\rho}$ is the sum of the fundamental weights of $G$; it is the highest weight of
an irreducible representation of
$G$ contained in $St^{\otimes d}$. The assumption of Galois ordinarity, denoted {\bf (GO)} in the sequel, is:

 \begin{itemize}
\item The image $\rho_\pi(D_p)$ of the decomposition group is contained in $\hat{G}$,
\item There exists $\hat{g}\in \hat{G}({\cO})$ such
that
$$\rho_\pi(D_p)\subset \hat{g}\cdot \hat{B}(\cO)\cdot \hat{g}^{-1},$$
\item the restriction of the conjugate $\rho_\pi^{\hat{g}}$ to $I_p$, followed by the quotient by the unipotent radical
$\hat{g}\cdot \hat{N}\cdot \hat{g}^{-1}$ of $\hat{g}\cdot \hat{B}\cdot \hat{g}^{-1}$ factors through
$-(\lambda+\tilde{\rho})\circ \chi:I_p\rightarrow \hat{T}(\Z_p)$.
\end{itemize}

{\bf Example:}

For $g=1$, $\lambda=(n;n)$ corresponds to the representation $Sym^n\,(St)$ of $GL(2)$, and $\tilde{\rho}=(1;1)$ corresponds to
 $St$. Then the weights of the (2-dim.) spin representation of $GSpin_3$ are $\hat{\varpi}=({1\over 2};{1\over 2})$ and
$\hat{\varpi}^{w_0}=(-{1\over 2};{1\over 2})$; hence the composition of $\chi$, $-(\lambda+\tilde{\rho})$ and the spin
representation (modulo unipotent radical) gives the diagonal matrix $diag(\chi^{-(n+1)},1)$ (modulo Weyl group), which is
the usual formula for an ordinary representation coming from an ordinary cusp form of weight $k=n+2$:
$$\rho_f|_{D_p}\,\cong \,\left(\begin{array}{cc}1&*\\0&\chi^{-n-1}\end{array}\right)$$

\noindent{\bf Convention:}
In
the rest of the paper, we make the abuse of notation to write $\hat{B}$, resp.
$\hat{N}$,
$\hat{T}$, instead of their respective conjugates by $\hat{g}$: $\hat{g}\cdot \hat{B}\cdot \hat{g}^{-1}$ and so on.
With this convention, we have $\overline{\rho}_\pi(I_p)\subset  \hat{B}(k)$.

Relative to the triple $(\hat{G},\hat{B},\hat{T})$, we have the notion of dominant characters $\mu\in X^*(\hat{T})$ and Weyl
classification of highest weight $\cO$-representations of $\hat{G}$, provided $p-1>\vert \mu+\rho\vert$ (see Polo-T.
\cite{PT}).
Let $\hat{\varpi}$ be the minuscule weight of $\hat{G}$. As already calculated, its
coordinates are:
$$({1\over 2},\ldots,{1\over 2};{1\over 2})$$ 

\begin{lem} For any $\sigma\in I_p$, 
$$(3.3.1)\quad\hat{\varpi}(\overline{\rho}_\pi(\sigma)\,{\rm mod.}\,\hat{N}(k))=\omega^{-\w}(\sigma).$$
and similarly, for the lowest weight $\hat{\varpi}^{w_0}$
$$(3.3.2)\quad\hat{\varpi}^{w_0}(\overline{\rho}_\pi(\sigma)\,{\rm mod.}\,\hat{N}(k))=1.$$
\end{lem}

{\bf Proof:} By {\bf (GO)}, the left-hand side is given by $\hat{\varpi}\circ [-(\lambda+\tilde{\rho})]\circ \omega
 (\sigma) $; therefore, the desired relation follows from the Lemma above, with $\mu=\lambda+\tilde{\rho}$. 
Indeed, the coordinates of $\lambda+\tilde{\rho}$ in $\Z^g\times\Z$ are $(a_g+g,\ldots,a_1+1;a_g+\ldots+a_1+d)$, hence
the scalar product $\langle \hat{\varpi},\lambda+\tilde{\rho}\rangle$ is equal to $\sum_i{a_i+i\over 2}+{(\sum_i a_i)+d\over
2}$, that is,
${\w\over 2}+{\w\over 2}$ {\it i.e.} $\w$.
Similarly for (3.3.2).

\vskip 5mm

\noindent{\bf Comments:}

\noindent 1)  Let us introduce the condition of automorphic ordinarity:

\vskip 5mm

{\bf (AO)} For each $r=1,\ldots,g$, $v(\theta_\pi(T_{p,r}))=a_{r+1}+\ldots+a_1$,

\vskip 5mm

\noindent where $T_{p,r}$ is the classical Hecke operator associated to the double class of 
$$diag(1_{r},p.1_{2g-2r},p^2.1_r).$$ 
We conjecture that for any $g$, if $\rho_\pi$ is 
a subquotient of $W_{\pi,p}$, then {\bf (AO)}
implies {\bf (GO)}. It is well-known for
$g=1$ (\cite{W1} Th.2.2.2,
\cite{H} and \cite{MT}).

 Consider the statement

${\bf KM_g}(\pi_f,p)$: {\it If $p$ is prime to
$N$, the slopes of the crystalline Frobenius on the isotypical component
$\bD_{crys}(W_{\pi,p})$ are the
$p$-adic valuations of the roots of the polynomial $\theta_\pi(P_p(X))$, reciprocal of the
$p$-Euler factor of the automorphic $L$-function of $\pi$.}

For $g=2$, we have seen in 3.1 that $W_{\pi,p}^{s.s.}$ is $\rho_{\pi}$-isotypical (assuming its absolute irreducibility).
We have observed (Proposition 7.1 of
\cite{TU}) that if ${\bf KM_2}(\pi,p)$ holds and if $\pi$ is stable at infinity, the
condition {\bf (AO)} for $\pi$  implies {\bf (GO)}. In a recent preprint,  E. Urban \cite{U2} has proven
${\bf KM_2}(\pi,p)$; thus, for $g=2$, if $\pi$ is stable at $\infty$, {\bf (AO)} implies {\bf (GO)}.

\vskip 5mm

2) If $\pi_p$ is in the principal series (for instance, if $\pi$ is unramified at $p$), and
if the
$p$-adic representation
$\rho_\pi$ is, say, potentially crystalline at
$p$ (for instance, crystalline), one can ask in general the following question. 

On one hand, the local component
$\pi_p$ of $\pi$ at $p$ is
unitarily induced from $\phi$ for a character
$\phi:T(\Q_p)\rightarrow \C^\times$; we defined in Sect.3.2.2 the local Galois representation
$\sigma(\pi_p)$ of the Weil group
$W_{\Q_p}$ given by 
$$W_{\Q_p}\rightarrow \Q_p^\times\rightarrow \hat{T}(\C)\subset \hat{G}(\C)$$
where $\Q_p^\times\rightarrow \hat{T}(\C)$ is given by the twist $\tilde{\phi}$ through the identification (3.2.2.2). This
representation is rational (the traces belong to some number field). 

Let us consider on the other hand the restriction to
$D_p$ of $\rho_\pi$. By applying the (covariant) Fontaine's functor $D_{pcris}$ (cf. Fontaine, ExposŽ III AstŽrisque 223), we
obtain a representation
$'\rho_{\pi,p}$ of the Weil group
$W_{\Q_p}$: 

$$'\rho_{\pi,p}:W_{\Q_p}\rightarrow GL_V$$

One can conjecture a compatibility at $(p,p)$ between the local and global Langlands correspondences, namely that the
$F$-semisimplification of the two rational representations $'\rho_{\pi,p}$ and $\sigma(\pi_p)$ are isomorphic (where
$a:W_{\Q_p}\rightarrow \Q_p^\times$ is the map induced by class-field theory, sending arithmetic Frobenius to $p$, and the twist is to pass from Langlands parameters to ``{Hecke}'' parameters). This
fact is known in the following cases:
\begin{itemize}
\item for $g=1$, by well-known theorems of Scholl and Katz-Messing,
\item for $g=2$, for a representation $\pi$ on
$GSp(4)$ which is the base change from $GL(2,F)$ ($F$ real quadratic) of a Hilbert modular form which is in the
discrete series at some finite place, and which is unramified at places above $p$ (in which case $\rho_{f,p}$, hence
$\rho_{\pi,p}$ is crystalline at $p$ by Breuil's theorem \cite{Br}). This is a particular case of a theorem of T. Saito
\cite{Sa}.

\end{itemize}  
 
Note however that this statement does not allow one to recover the representation
$\rho_{\pi,p}=\rho_\pi|_{D_p}$ (because it says nothing about the Hodge filtration) unless we assume
it is ordinary (in the usual geometric sense). More precisely, we have two parallel observations: 
\begin{itemize}
\item Let us assume that $\rho_{\pi,p}$ is crystalline; then the assumption of geometric ordinarity
means that the eigenvalues $(\xi_B^{-1})_{B\subset \{1,\ldots,g\}}$ of the crystalline Frobenius are
such that the $ord_p(\xi_B)$ ($B\subset\{1,\ldots,g\}$) coincide (with multiplicities) with the
Hodge-Tate weights; these numbers, if
$\pi$ is stable at infinity,  should be (as mentioned at the end of Sect.2.3.2)
$j_B=\sum_{i\in B}(a_+i)$ ($B\subset A=\{1,\ldots,g\}$). These quantities can also be written
$$\langle
\hat{\varpi}^{w_B},(\lambda+\tilde{\rho})\rangle=\hat{\varpi}^{w_B}\circ(\lambda+\tilde{\rho})$$ where
$w_B\in W_{\hat{G}}$ is the element of the Weyl group such that for
$\hat{t}=(t_g,\ldots,t_1,[u,\zeta])\in \hat{T}$ and $w_B(\hat{t})=\hat{\theta}$, $\theta_i=t_i^{-1}$
if and only if
$i\in B$ and all its other components are those of $\hat{t}$. Therefore, it implies by
Fontaine-Laffaille theory that $\rho_\pi$ is ordinary at
$p$ in the precise sense of ${\bf (GO)}$. Thus the conjunction of geometric ordinarity, and of
stability of $\pi$ at $\infty$ (together with the complete determination of Hodge-Tate weights of
$\rho_\pi$) implies {\bf (GO)}. 

\item Let us assume $\pi$ is unramified at $p$; let us introduce complex numbers $\theta_i$'s and
$\zeta$, such that for any
$t\in
\hat{T}(\C)$ mod.
$W_{\hat{G}}$,

 $$|t_i|_p^{\alpha_i}=\theta_i^{-ord_p(t_i)},\quad {\rm and}\, |z|_p^{\gamma}=\zeta^{-ord_p(z)},$$ 

we can rewrite $(3.2.2.4)$ as
$$\tilde{\phi}(p)=(\theta_g^{-1},\ldots,\theta_1^{-1},[(\theta_g\cdot\ldots\cdot\theta_1)^{-1/2},p^{d/2}\cdot\zeta^{-1}])$$
The composition with $\sp$ gives a complex diagonal matrix whose entries are inverse to the $2^g$
algebraic integers
$$\xi_J=(\prod_{i\in J}\theta_i^{-1}\cdot\prod_{i\notin J}\theta_i)^{1/2})\cdot \zeta.$$

 The Automorphic Ordinarity Conjecture for the $p$-adic embedding $\iota_p$ states
$$ord_p(\iota_p(\xi_J))=\sum_{i\in J}(a_i+i),\quad\mbox{\rm for any}\, J.$$

Therefore,  the quantities $x_i=-ord_p(\iota_p(\theta_i))$ and $y=ord_p(\iota_p(\zeta))$ satisfy  the
linear system in $(x_g,\ldots,x_1;y)\in\Z^{g+1}$: 
$$ -{y+d+\sum_{i\in J}x_i-\sum_{i\notin J}x_i\over 2}=\sum_{i\in J}(a_i+i).$$ It contains a Cramer
system. Therefore, assumption {\bf (AO)} implies 
$$ord_p\theta_i=-(a_i+i),\quad ord_p(\zeta)=a_g+\ldots+a_1$$
 up to permutation of the coordinates. This can be rewritten as an equality in 
$Hom(\Q_p^\times,\hat{T}(K)/\hat{T}(\cO))$:
$$\iota_p\circ \tilde{\phi}=-(\lambda+\tilde{\rho}).$$ We conclude that {\bf (AO)} together with
${\bf KM_g}(\pi,p)$ implies (part of) the compatibility conjecture at
$(p,p)$: the ($p$-adic orders of) the eigenvalues counted with multiplicities of
$D_{cris}(\rho_\pi)(Frob_p)$ coincide with those of
$\sigma(\pi_p)(Frob_p)$.
\end{itemize}
     
\vfill\eject

\section{Crystals and connections}
\subsection{de Rham and crystalline cohomology of open varieties} Let $f: \overline{X}\rightarrow S$
be a smooth proper morphism of schemes; $X\subset \overline{X}$ be an open immersion above $S$,
with complement a relative Cartier divisor $D\rightarrow S$ with normal crossings and smooth
irreducible components. Let $\overline{\cV}$ be a coherent sheaf over $\overline{X}$ endowed with an
integrable connection
$\nabla$ with logarithmic poles along $D$; let $\cV$ its restriction to $X$. Let
$\cI(D)$ be the sheaf of ideals defining $D$. Then the relative de Rham cohomology sheaves
${\cal H}^j_{dR}(X/S,\cV)$ are defined as
$$(2.1)_\emptyset\quad {\bf
R}^jf_*(\overline{\cV}\otimes_{\cO_{\overline{X}}}\Omega^\bullet_{\overline{X}/S}(log\,D))$$ 

Let us now introduce a complex
$$\Omega^\bullet_{\overline{X}/S}(-log\,D)=\Omega^\bullet_{\overline{X}/S}(log\,D)
\otimes_{\cO_{\overline{X}}}\cI(D)$$ We define the cohomology sheaves with compact support
${\cal H}^j_{dR,c}(X/S,\cV)$ by:    
 
$$(2.1)_c\quad {\bf
R}^jf_*(\overline{\cV}\otimes_{\cO_{\overline{X}}}\Omega^\bullet_{\overline{X}/S}(-log\,D))$$

If $S=Spec\,k$, we write $H^j_{dR}$ instead of ${\cal H}^j_{dR}$.

A priori, these definitions depend on the compactification $\overline{X}$ of $X$. One can show for
$S=Spec\,K$ and $\cV$ trivial that the resolution of singularities implies the independence of the
compactification (Th\'eor\`eme 2.11 of \cite{Mo}).  

For the crystalline cohomology there is a similar definition. Our reference is \cite{Ka}, section 5,6.
We use the language of logarithmic schemes; as noted by Kato in Complement 1 of his paper, his results
are compatible with Faltings theory of crystalline cohomology of open varieties: in Faltings approach,
a logarithmic structure on $\overline{X}$ is a family $({\cal L}_i,x_i)_{1\leq i\leq r}$ where
${\cal L}_i$ is an invertible sheaf and $x_i$ a global section thereof, these data always define a
logarithmic scheme in Kato's sense (while the converse is false).  
 Let
$(S,I,\gamma)$ a triple where
$S$ is a scheme, $\cI$ is a quasi-coherent nilpotent ideal of $\cO_S$
 and $\gamma$ is a divided power structure on
$\cI$ (PD-structure, for short). Let $S_0$ the closed subscheme defined by
$\cI$; we consider  a smooth morphism $\overline{X_0}\rightarrow S_0$ and $D_0$ a relative Cartier
divisor with normal crossings. It defines a logarithmic structure $M=\{g\in \cO_{\overline{X}_0}; g\,
\hbox{\rm is invertible outside}\,D\}\subset \cO_{\overline{X}_0}$. One defines the logarithmic
crystalline site of
$(\overline{X}/S_0)^{log}_{crys}$ as in Kato \cite{Ka} Sect.5.2. The objects are $5$-uples
$(U,T,M_T,i,\delta)$ where $U\rightarrow \overline{X_0}$ is \'etale, $(T,M_T)$ is a scheme with fine
logarithmic structure over
$S$, $i:(U,M)\rightarrow (T,M_T)$ is an exact closed immersion over $S$ and $\delta$ is a divided
power structure compatible with $\gamma$. Morphisms are the natural ones. On this site, the structural
sheaf
$\cO_{\overline{X}_0/S}$ is defined by
$$\cO_{\overline{X}_0/S}(U,T,M_T,i,\delta)=\Gamma(T,\cO_T).$$

\begin{de} A crystal on
$(\overline{X}_0/S)^{log}_{crys}$ is a sheaf $\cV$ in $\cO_{\overline{X}_0/S}$-modules satisfying the
following condition: for any morphism $g:T'\rightarrow T$ in $(\overline{X}_0/S)^{log}_{crys}$,
$g^*\cV_T\rightarrow
\cV_{T'}$ is an isomorphism. Here $\cV_T$ and $\cV_{T'}$ denote the sheaves on $T_{et}$ and
$T'_{et}$ defined by $\cV$.  
\end{de}

Let $(\overline{X},D)$ be a lifting of $(\overline{X}_0,D_0)$ to $S$.
Then, by Th.6.2 of \cite{Ka} (see sect.4.2 for more details), the data of a crystal on
$(\overline{X}_0/S)^{log}_{crys}$ is equivalent to that of an $\cO_{\overline{X}}$-module 
${\cal M}$ endowed with a quasi-nilpotent integrable connection with logarithmic singularities 
$$\nabla:{\cal M}\rightarrow {\cal M}\otimes_{\cO_{\overline{X}}}\Omega^1_{\overline{X}/S}(log\, D)$$

To compute the cohomology sheaves of a crystal, we apply the spectral sequence 
$$Rf_{crys,*}\cV=Rf_{et,*}(Ru_*\cV)$$  where $u$ is the canonical projection from the site
$(\overline{X}_0/S)^{log}_{crys}$ to $\overline{X}_{0\, et}$. It is defined, for a sheaf $\cV$ on
$(\overline{X}_0/S)^{log}_{crys}$, and for any \'etale morphism $U\rightarrow \overline{X}_0$, by
$$(u_*\cV)(U)=\Gamma(U,\cV)$$

 Moreover, 
$$Ru_*\cV\cong {\cal M}\otimes_{\cO_{\overline{ X}}}\Omega^\bullet_{\overline{ X}/S}(log\,D).$$

Again, by Th.2.11 of \cite{Mo}, one can show, assuming the resolution of singularities that for
$S=\Z/p^n\Z$, $S_0=\Z/p\Z$ this definition does not depend on the lifting. 

{\bf Remark:} In our case, one even does not need the resolution of singularities. It will be a
consequence of the comparison theorem!

These definitions transfer to the compact support case; it is explained in \cite{Fa} p.58. The
cohomology sheaves $Rf_{crys,*,c}\cV$ are computed by a similar spectral sequence
$$Rf_{crys,*,c}\cV=Rf_{et,*}(Ru_{*,c}\cV)$$ where $u_{*,c}$ is defined, for a sheaf $\cV$ on
$(\overline{X}_0/S)^{log}_{crys}$ and a \'etale morphism
$g:U\rightarrow \overline{X}_0$, by
$$(u_{*,c}(\cV)(U)=\Gamma(U,\cV\otimes_{{\cal O}_U} g^*\cI(D_0))$$

 One has also:
$$Rf_{crys,*,c}\cV=Rf_{et,*}{\cal M}\otimes_{\cO_{\overline{ X}}}\Omega^\bullet_{\overline{
X}/S}(-log\,D)$$  

This result can be proven as in the case without support; it will be explained in the next section.

\subsection{$L$-construction}

In the proof of Theorem 2 below, we will apply the crystalline $L$-construction in the logarithmic
setting (in the classical crystalline setting, cf. Chapt.6 of
\cite{BeO}); we want to explain the definitions and results here.

Let $(S,I,\gamma)$ a triple where
$S$ is a scheme, $\cI$ is a quasi-coherent ideal of $\cO_S$
 and $\gamma$ is a PD-structure on
$\cI$. Let $S_0$ the closed subscheme defined by
$\cI$; we consider a smooth morphism $\overline{X}_0\rightarrow S_0$ and $Y_0$  a relative Cartier
divisor with normal crossings. Let $(\overline{X},Y)$ be a lifting of $(\overline{X}_0,Y_0)$ to $S$;
we suppose that there exists an integer $m>0$ such that $p^{m}\cO_{\overline{X}}=0$. Let
$Z_1,\ldots,Z_a$ be the irreducible components of $Y$. Let $\Xi$ be the blowing-up of
$\overline{X}\times_{S}\overline{X}$ along the subscheme $\sum_{i}( Z_i\times_{S} Z_i)$. Let
$\overline{ X}{\hat\times}_S \overline{ X}$ be the complement in $\Xi$ of the strict transforms of
$\overline{X}\times{Z}_i$ and ${Z}_i\times \overline{X}$, $1\leq i\leq r$ and let
$\tilde{Y}$ be the exceptional divisor in $\overline{X}{\hat\times}_S \overline{ X}$; it  is a divisor
with normal crossings. The couple $(\overline{ X}{\hat\times}_S
\overline{X},\tilde{Y})$ is the categorical fiber product of 
$(\overline{X}, Y)$ by itself over $S$, in the category of logarithmic schemes (cf.
\cite{Fa} IV,c). Locally, if $x_1,\ldots,x_d$ are local coordinates of
$\overline{X}$ over $S$ such that
$Y$ is defined by the equation $x_1\ldots x_a=0$, then $\overline{X}{\hat\times}_S
\overline{X}$ is the relative affine scheme given as spectrum of
 
$$S[x_i{\scriptstyle\otimes} 1,1{\scriptstyle\otimes} x_i]_{1\leq i\leq d}[u^{\pm 1}_j]_{1\leq
j\leq a}/(x_j{\scriptstyle\otimes} 1. u_j - 1{\scriptstyle\otimes} x_j)_{1\leq j\leq a}$$

 and $\tilde{Y}$ is defined by the equation
 $x_1{\scriptstyle\otimes} 1\ldots x_a{\scriptstyle\otimes} 1=0$ (or
 $1{\scriptstyle\otimes} x_1\ldots 1{\scriptstyle\otimes} x_a=0$).

 We endow this product with a PD-structure as follows. Let
${\cal J}$ be the PD-envelope of the diagonal immersion $\overline{
X}\to
\overline{X}{\hat\times}_S\overline{X}$. 
In the local coordinates above, ${\cal J}$ is the
$PD$-polynomial algebra ${\cal O}_{\overline{X}} <v_1,\ldots,v_a,\xi_{a+1},\ldots,\xi_{d}>$ where
$v_i=u_i-1$ and
 $\xi_i=x_i{\scriptstyle\otimes} 1 -1{\scriptstyle\otimes} x_i$. 

We denote by ${\cal J}^{n}$ the $n^{th}$ order divided power neighborhood :
 ${\cal J}^{n}={\cal J}/{\cal I}_\Delta^{[n+1]}$ where ${\cal I}_\Delta$ is the ideal of the diagonal
immersion.

 Let $\cM$ be a sheaf of ${\cal O}_{\overline{X}}$-modules. We recall the interpretation of a
connection on $\cM$ in terms of an HPD-stratification in our context.  For us, the notion of an HPD
stratification on
$\cM$ is defined word for word as in
\cite{BeO} Sect.4.3 (which treats the crystalline situation on $\overline{X}_0$, without the divisor
$Y_0$). It consists namely in the datum of a
$\cal J$-linear isomorphism

$$\epsilon : {\cal J}\otimes_{{\cal O}_{\overline{X}}}{\cM}\rightarrow {\cM}\otimes _{{\cal
O}_{\overline{X}}}{\cal J}$$
 
such that $\epsilon$ reduces to identity modulo $\cI_\Delta$ and the natural
cocycle condition on 
$\overline{X}{\hat\times}_S \overline{ X}{\hat\times}_S \overline{X}$ holds (\cite{BeO} def.2.10). 
In the case
$\cM =\cal J$, we have two canonical HPD stratifications. The first is induced by
extending by (left) $\cI$-linearity the map  
$\theta : c{\scriptstyle\otimes}d \mapsto
((1{\scriptstyle\otimes}d){\scriptstyle\otimes}(1{\scriptstyle\otimes}c)$; it makes use of the right
module structure of $\cI$ over ${\cal O}_{\overline{X}}$. The second is given similarly by
tensoring on the left by $\cI$ over ${\cal O}_{\overline{X}}$ the left-hand side of
$\iota : c{\scriptstyle\otimes}d \mapsto
((c{\scriptstyle\otimes}1){\scriptstyle\otimes}(1{\scriptstyle\otimes}d)$; it uses the structure
of left ${\cal O}_{\overline{X}}$-module of $\cI$.

Also, as in \cite{BeO} 4.4, one recalls the notion of PD-differential operator.
Let $\cM$ and $\cN$ be two 
${\cal O}_{\overline{X}}$-modules.

 A PD-differential  operator
${\cM} \to {\cN}$ of order $\leq n$ (resp. HPD-differential operator) is a ${\cal
O}_{\overline{ X}}$-linear map
${\cal J}^{n}\otimes {\cM}\to {\cN}$ (resp.${\cal J}\otimes {\cM}\to {\cN}$). Every
PD-differential operator $\delta : {\cal J}^{n}\otimes {\cM}\to {\cN}$ induces a classical
differential operator $\delta^{b} : {\cM} \to {\cN}$ of order $n$ with "cologarithmic zeroes"
along $Y$.

The importance of these notions for us stems from the following theorem whose proof runs exactly as
in the ``{classical}'' case (\cite{BeO} Theorem 4.12). For that, we introduce the notion of a
quasi-nilpotent connection in the sense of \cite{BeO} 4.10 (but in our log setting, again):

\begin{de}
A connection $\nabla$ on $\cM$
is quasi-nilpotent if for any local section
$s$ of $\cM$ with local coordinates $x_1,\ldots,x_d$ on $ X$ such that $Y$ is defined by the
equation $x_1\ldots x_a=0$, there exists a positive integer
$k$ such that $\prod_{0\leq j\leq k-1}(\nabla(x_{i}\partial/\partial x_{i})-j)^{k}(s)=0$ for $1\leq i\leq a$ and
$(\nabla(\partial/\partial x_{i}))^{k}(s)=0$ for $a+1\leq i\leq d$).
\end{de}

\begin{thm} The data of an HPD stratification on
$\cM$ is equivalent to the data of a logarithmic integrable connection $\nabla$ on
$\cM$ wich is quasi-nilpotent.
\end{thm} 

  Then, Grothendieck's
 linearization functor $L$ is defined as follows.
Let $\cH$ be the category of ${\cal
O}_{\overline{X}}$-modules with morphisms given by HPD-differential operators and $\cC$ to the
category of crystals over $\overline{X}_0/S)_{crys}^{log}$. 
For any sheaf $\cM$ of ${\cal O}_{\overline{X}}$-modules,  
we endow the ${\cal O}_{\overline{X}}$-module ${\cal J}\otimes _{{\cal O}_{\overline{X}}}{\cM}$.
with the
HPD-stratification $\epsilon_{L(\cM)}$
$${\cal J}\otimes {\cal J}\otimes {\cM} \stackrel{\iota\otimes id_{\cM}}{\longrightarrow} {\cal
J}\otimes {\cal J}\otimes {\cM} 
\stackrel{id_{\cal J}\otimes f}{\longrightarrow} {\cal J}\otimes{\cM} \otimes  {\cal J}$$

where $f:{\cM} \otimes  {\cal J}\rightarrow {\cal J}\otimes{\cM}$ interchanges the factors. In other
words, the HPD-stratification is given by:
$$ (a{\scriptstyle\otimes}b){\scriptstyle\otimes} (c{\scriptstyle\otimes}d) {\scriptstyle\otimes} m
\mapsto
((ac{\scriptstyle\otimes}b){\scriptstyle\otimes}m{\scriptstyle\otimes}(1{\scriptstyle\otimes}d)$$ 
 
\begin{de} The covariant functor $L:\cH\rightarrow \cC$ is defined by:
\begin{itemize} 
\item For any sheaf $\cM$ of ${\cal O}_{\overline{X}}$-modules,  
$L({\cM})$ is the crystal corresponding to the ${\cal O}_{\overline{X}}$-module with
HPD-stratification $({\cal J}\otimes _{{\cal O}_{\overline{X}}}{\cM},\epsilon_{L(\cM)})$. 

\item For an HPD-differential operator
 $\varphi :{\cM}\to {\cN}$ ( that is, an ${\cal O}_{\overline{X}}$-linear map
$\varphi :  {\cal J}\otimes {\cM}\to {\cN}$), $L(\varphi) :L(\cM)\to L(\cN)$ is the morphism of
crystals corresponding to the ${\cal O}_{\overline{X}}$-linear morphism compatible with
HPD-stratifications, given by the composition:
$${\cal J}\otimes {\cM} \stackrel{\iota\otimes id_{\cM}}{\longrightarrow} {\cal J}\otimes {\cal
J}\otimes {\cM} 
\stackrel{id_{\cal J}\otimes \varphi}{\longrightarrow} {\cal J}\otimes {\cN}$$
\end{itemize}
\end{de}

We refer to \cite{BeO} Sect.2,Sect.6 for more details. Note that since
$\cal J$ is locally free, the functor $L$ is exact.

The correspondence between crystals on
$(\overline{X}_0/S)^{log}_{crys}$ and $\cO_{\overline{X}}$-module 
${\cal M}$ endowed with a quasi-nilpotent integrable connection with logarithmic singularities,
is then given by the following rule:
Let $pr_1,pr_2 : {\cal J}\to {\overline X}$ be the two canonical projections. 
If $\cV$ is a crystal on $(\overline{X}_0/S)^{log}_{crys}$, let ${\cal M}={\cV}_{\overline X}$ be the
evaluation of $\cV$ on $\overline X$. The defining condition of a crystal produces an isomorphism:
$$\epsilon : pr_2^*{\cal M}\simeq pr_1^*{\cal M}$$
This induces an integrable quasi-nilpotent logarithmic connection on $\cal M$ as explained above.
Conversely, by theorem 4, every logarithmic integrable connection on
$\cM$ wich is quasi-nilpotent induces an HPD stratification on
$\cM$. If $(U,T,M_T,i,\delta)$ is an object of the crystalline site, then by smoothness, etale locally on
$T$, the morphism $({\overline X}_0,D_0)\to ({\overline X},D)$ extend to a morphism $h : (T,M_T)\to
({\overline X},D)$. We define $\cV_T$ to be $h^*{\cal M}$. If we have two such $h_i : (T,M_T)\to
({\overline X},D)$ ($i=1,2$), then there exist $h^\prime : (T,M_T)\to ({\cal J},M_{\cal J})$ such that
$h_i=h^\prime pr_i$ and $\epsilon$ induces an isomorphism $h^*_1{\cal M}\simeq h^*_2{\cal M}$. Thus
$\cV$ is well defined.

It is not hard from the classical case (Theorem
6.12 of \cite{BeO}), to deduce the following crystalline Poincar\'e lemma.

\begin{lem} Let $\cV$ be a crystal on $(X_0/S)_{crys}^{log}$ and ${\cal M}$ the associated
${\cal O}_{\overline{X}}$-module with its integrable connection. Then the complex of crystals
$L(\cM\otimes\Omega^{\bullet}_{\overline{X}}(log\,{Y}))$ is a resolution of $\cV
$.
\end{lem}      

Finally, the same argument as in the classical theory (\cite{BeO} Sect.5.27) shows also the following
useful lemma:

\begin{lem} Let $\cM$ be a sheaf of ${\cal O}_{\overline{ X}}$-modules and ${\cal I}(Y)$ the ideal
of definition of
$ Y$ . Then: 
$$Ru_{*}L(\cM)=\cM \hbox{ and } Ru_{*,c}L(\cM)=\cM\otimes {\cal I}(Y).$$
\end{lem}

Combining Lemma 1 and 2 above, we deduce:
$$Ru_*{\cal V}\cong {\cal M}\otimes_{{\cal
O}_{\overline{X}}}\Omega^\bullet_{\overline{X}/S}(log\,D)$$ and
$$Ru_{*,c}{\cal V}\cong {\cal M}\otimes_{{\cal
O}_{\overline{X}}}\Omega^\bullet_{\overline{X}/S}(-log\,D).$$

\subsection{The Gauss-Manin connection}
As in section 4.1, $\overline{X}$ is a smooth $S$-scheme (not necessarily proper), $X$ an $S$-open scheme of $\overline{X}$
such that $D=\overline{X}-X$ is a divisor with normal crossings over $S$. Let $f: \overline{\cX}\rightarrow \overline{X}$ be a
proper morphism such that $\overline{\cX}$ is smooth over $S$, $f$ is smooth over $X$
and ${\cal D}=\overline{\cX}\times_{\overline{X}}D$ is a relative divisor with normal crossings (such $f$ is
called semi-stable, see
\cite{Ill}). We have a relative de Rham complex with logarithmic poles 
$$\Omega^\bullet_{\overline{\cX}/\overline{X}}(log{\cal D}/D)=\Omega^\bullet_{\overline{\cX}/S}(log{\cal D})/
f^*\Omega^\bullet_{\overline{X}/S}(log D)$$
As explained in \cite{Katz} (see also \cite{Ill}), we have a Gauss-Manin connection with logarithmic poles along $D$, 
on the locally free $\cO_{\overline{X}}$-module:
$${\cal E} =R^*f_*(\Omega^\bullet_{\overline{\cX}/\overline{X}}(log{\cal D}/D))$$
The restriction of $\cal E$ to $X$ is the usual Gauss-Manin sheaf $R^*f_{\vert{\cal X}*}\Omega^\bullet_{{\cX}/X}$ and $\cal E$
is the Deligne's canonical extension to $\overline X$. The Gauss-Manin connection on $\cal E$ is integrable and if $\cO_S$ is
killed by a power of $p$, then this connection is quasi-nilpotent (\cite{Katz}).

\vfill\eject

\section{BGG resolutions for crystals}

 Let
$B=T.N$ resp. $Q=M\cdot U$ be the Levi decomposition of the upper triangular subgroup of $G$, resp.
of the Siegel parabolic, viewed as group schemes over $\Z$. We keep the notations of the introduction for the weights of $G$.
 Let
${\bf V}=<e_g,\ldots,e_1,e_1^*,\ldots,e_g^*>$ be the standard
$\Z$-lattice on which
$G$ acts; given two vectors $v,w\in \bV$, we write $<v,w>={}^tvJw$ for their symplectic product. 
$Q$ is the stabilizer of the standard lagrangian lattice $\bW=<e_g,\ldots,e_1>$; we have
${\bf V}=\bW\oplus
\bW^*$; $M=L_I$ is the stabilizer of the decomposition $(\bW,\bW^*)$;  one has
$M\cong GL(g)\times
\bG_m$. Let
$B_M=B\cap M$ be the standard Borel of $M$. Let
$\Phi$, resp. $\Phi_M$  be the set of roots of
$(G,B)$, resp. $(M,B_M)$ and let $\Phi^M=\Phi-\Phi_M$. We denote by $\Phi^{\pm}$, resp.
$\Phi^{\pm}_M$, 
$\Phi^{M\pm}$, the set of positive/negative roots in $\Phi$, resp. $\Phi_M$, $\Phi^M$.

\vskip 5mm

\subsection{Weyl modules over $\Z_p$}

   From this section on, the notations
$\g$,
$\q$, (and
$\m$ but there should not be confusion with the maximal ideal of the Hecke algebra) stand for the Lie algebras
over
$\Z$ of the corresponding group schemes. The Kostant-Chevalley algebra 
$\U=\U(\g)$ (resp. $\U(\q)$, $\U(\m)$) is the subring 
 of the rational enveloping algebra
$U(\g_\Q)$ (resp. $U({\q}_\Q)$, resp. $\U(\m_\Q)$) generated over $\Z$ by ${X^{n}\over n!}$
 with $X\in {\bf
g}_\alpha$,
$\alpha\in\Phi$ (resp.
$\alpha\in\Phi-\Phi^{M-}$, $\Phi_M$) , $n\geq 0$ an integer. 
There is a natural ring epimorphism 
$\U(\q) \to \U(\m)$. A $\g$-stable lattice of a $G_\Q$-representation which is
$\cal U$-stable is called $\g$-admissible (see \cite{Bo}, Sect.VIII.12.7 and 8) same thing for a
$\m$-lattice which is $\U(\q)$-stable.

\subsubsection{Admissible lattices}

In this section, we explain how one can construct Weyl modules over $\Z_{(p)}$ by plethysms when the highest
weight is $p$-small: $\vert \lambda\vert<p$. This construction is used in Appendix II to give a
construction by plethysms of the crystals (resp. filtered vector bundles) over a toroidal compactification of
the Siegel variety over $\Z_p$, associated to irreducible representations whose highest weights are
$p$-small.

\vskip 5mm

If $\lambda$ is a fundamental weight, then the irreducible representation
$V_\lambda$ of $G$ has a canonical admissible lattice $V(\lambda)_\Z$ for the Chevalley order ${\g}$
\cite {Bo} p.206. For another dominant weight $\lambda\in X^+$, several admissible lattices exist
over $\Z$. However, given an prime $p$, we have shown in \cite{PT}, Sect.1.2, that for
$\lambda=(a_g,\ldots,a_1;c)$ such that $a_g+a_{g-1}+g+(g-1) <p$, these lattices all coincide
after tensoring by the localization 
$\Z_{(p)}$ of
$\Z$ at $p$.
Note that our condition $\vert\lambda+\rho\vert<p-1$ implies $a_g+a_{g-1}+g+(g-1) <p$.

For such a weight, let us recall the construction by plethysms of this unique admissible
$\Zp$-lattice
$V_{\lambda,
\Z_p}$. It will be used systematically in the sequel as it fits well in the construction of sheaves
over the Siegel modular variety.

 Let $s= \vert\lambda \vert$; hence $s<p$. For any $(i,j)$ with $1\leq i<j\leq n$, let 
 $\phi_{i,j}:\bV^{\otimes s}\rightarrow\bV^{\otimes (s-2)}$
 the contraction given by
$$v_1\otimes\ldots\otimes v_s\mapsto \langle v_i,v_j\rangle v_1\otimes\ldots\otimes
\hat{v}_i\otimes\ldots\otimes\hat{v}_j\otimes\ldots\otimes v_s ;$$
Let $\psi\in \bV^{\otimes 2}$ be the image of the symplectic form $<\,,\,>\in (\bV\otimes \bV)^*$ via
the identifications
$$ (\bV\otimes \bV)^*\,\cong\,  \bV^*\otimes \bV^*\,\cong\, \bV\otimes \bV$$
the last one being given by $\bV\cong \bV^*, v\mapsto <v,->$.

We consider for any $s\geq 2$ the maps $\psi_{i,j}:\bV^{\otimes s-2}\rightarrow
\bV^{\otimes s}$ obtained by inserting $\psi$ at $i$th and $j$th components. Observe that $\psi_{i,j}$ is injective.
Let $\theta_{i,j}=\psi_{i,j}\circ \phi_{i,j}\in End(\bV^{\otimes s})$.   Let $\bV^{<s>}$ be the submodule of
$\bV^{\otimes s}$ defined as intersection of the kernels of the $\theta_{i,j}$'s (note that $Ker\,
\theta_{i,j}=Ker\,\psi_{i,j}$).

As we shall see below, for $p> 2\cdot g$, $\bV^{<s>}_{\Z_{(p)}}$ is the image of $\bV^{\otimes s}$ by an idempotent in
the
$\Z_p$-algebra generated by the $\theta_{i,j}$'s inside $End_{\Z_{(p)}}( \bV^{\otimes s})$. 
Finally, by applying the Young symmetrizer
$c_\lambda=a_\lambda\cdot b_\lambda$ (see
\cite{FH} 15.3 and 17.3), whose coefficients are in ${\Z_{(p)}}$, to $\bV^{<s>}\otimes \Z_{(p)}$, one obtains the
sought-for  lattice $V_{\lambda,
\Z_{(p)}}$. 

\vskip 5mm

\begin{lem}  There exists an idempotent $e_s$ in the $\Z[{1\over g}]$-subalgebra of

$End_{\Z[{1\over g}]}(\bV^{\otimes s})$ generated by the $\theta_{i,j}$'s ($1\leq i<j\leq g$, such
that
$$\bV^{<s>}=e_s\cdot \bV^{\otimes s}.$$
\end{lem}

\noindent {\bf Proof:} Let 
$$\Phi=\bigoplus\phi_{i,j}:\bV^{\otimes s}\rightarrow \bigoplus_{1\leq i<j\leq s} \bV^{\otimes
(s-2)}$$Thus, 
$$\bV^{<s>}=Ker\,\Phi.$$ Similarly, put
$$\Psi:\sum_{i<j}\psi_{i,j}:\bigoplus_{1\leq i<j\leq s} \bV^{\otimes (s-2)}\rightarrow \bV^{\otimes
s}.$$ and
$$\Theta=\Psi\circ\Phi=\sum_{1\leq i<j\leq s}\theta_{i,j}.$$ Since 
$$\Phi\circ \Psi=(\times g),$$ we see that ${1\over g}\cdot \Theta$ is an idempotent. It belongs to
the  $\Z[{1\over g}]$-algebra generated by the $\theta_{i,j}$'s.

Thus,

$$\bV^{\otimes s}=\bV^{<s>}\oplus Im\,\Psi,\quad x=(x-{1\over g}\cdot \Theta(x))+{1\over g}\cdot
\Theta(x).$$

This decomposition of $\Z_{(p)}$-modules is $G$-stable. We put
$e_s=Id-{1\over g}\cdot \Theta$. This is the desired projector to
$\bV^{<s>}$.

To conclude:

\begin{cor}
For any prime $p$ which does not divide $2\cdot g$ and such that $p> s=\vert
\lambda\vert$, the module
$V_{\lambda,\Z_{(p)}}$ obtained by Construction 5.1 is  the image of
$\bV^{\otimes s}_{\Z_{(p)}}$ by an idempotent in the $\Z_{(p)}$-subalgebra of $End_{\Z_{(p)}}(\bV^{\otimes s})$
generated by permutations and the $\theta_{i,j}$'s. This algebra commutes to the $G$-action.
\end{cor}

 We apply a similar construction for a
$B_M$-dominant weight $\mu$ of $M$ with $\vert\mu\vert<p$. We denote by $W_{\mu,\Z_{(p)}}$ the canonical admissible
lattice of $W_{\mu}$ over $\Z_{(p)}$ given by the Young symmetrizer. 
It can be regarded as a $\U(\q)$-module via $\U(\q)\rightarrow \U(\m)$. 

\begin{lem} The subcategory of the category of $M$-representations, free and of finite rank over $\Z_p$, consisting of
representations of highest weight $<p$ is semisimple. 
\end{lem}
{\bf Proof} We have to show that there is no nontrivial  extensions in this subcategory. Let $\lambda$ and $\mu$ be
two $M$-dominant weights such that $\vert\lambda\vert<p$ and $\vert\mu\vert<p$.
$\lambda$ and $\mu$ are not in the same orbit for the action of the affine Weyl group (\cite{J},Part II,6.1).
 Let $W_\lambda$ and $W_\mu$ be the corresponding canonical admissible lattices over $\Z_p$, then
$Ext^{1}(W_\lambda,W_\mu)=0$ by the linkage principle (\cite{J},Part II,6.17, see also \cite{PT}, Sect.1.4,
Proposition). 
\bigskip 

\subsubsection{The BGG complex}

We are interested in a variant
of the ``{BGG complex}'' constructed in \cite{BGG} 
where one replaces the Borel subgroup by the parabolic $Q$. Over the field $\Q$, it is
defined in \cite{CF} Chapter VI, Prop.5.3 as the eigenspace for the infinitesimal character 
$\chi_{\lambda+\rho}$ inside the standard bar resolution of $V_{\lam,\Q}$: 
$$D(\lam)_\Q := \U_\Q \otimes_{\U(\q)_\Q} (\Lambda^\bullet(\g/\q)\otimes V_{\lam,\Q}).$$

Following \cite{BGG}, we show in \cite{PT} that this BGG complex admits a natural $\Z_{(p)}$-structure in
terms of integral Verma modules:

$$C(\lam)_{\Z_{(p)}}=\bigoplus_{w\in W^M}\U\otimes_{\U(\q)} W_{w(\lam+\rho)-\rho,\Z_{(p)}}$$

and we prove the following result. 
Let $D(\lam)_{\Z_{(p)}} := \U_{\Z_{(p)}} \otimes_{\U(\q)_{\Z_{(p)}}} (\Lambda^\bullet(\g/\q)\otimes V_{\lambda,\Z_{(p)}})$
be the standard $\Z_{(p)}$-complex, a natural $\Z_{(p)}$-version of the standard bar resolution over $\Q$ 
of $V_{\lambda,\Q}$.   

\begin{thm} 
Let $\lam\in X^+$ and let $p> \vert\lambda+\rho\vert$.  
Then there is a canonical morphism of complexes 
$j: C(\lam)_{\Z_{(p)}}\hookrightarrow D(\lam)_{\Z_{(p)}}$ such that 
\begin{itemize} 
\item it is injective with ${\Z_{(p)}}$-flat cokernel,
\item $Im(j_{\Q})$ is the BGG complex over $\Q$.
\end{itemize}
\end{thm}

{\bf Remarks:} 1) The BGG complex mentioned here is a variant for the parabolic $Q$ of the one defined in lemma 9.8 of
\cite{BGG} in the Borel case. For details concerning the differential maps, see Sect.2 of \cite{PT}.

2) The bound on $\lambda$ needed for proving this theorem is actually looser than $(\sum_{i=1}^g a_i)+d<p$: it is enough that
$a_g+a_{g-1}+g+(g-1)<p$.

3) We do not claim that these complexes are exact, as they are not. However, as we will see in
Sect.3.5 after applying the functor
$L$ to a sheaf construction (Sect.5.4 below), we will transform the dual of
$C(\lambda)_\bullet$ into a resolution of the sheafification of the dual of
$V_{\lambda,\Z_{(p)}}$. 

\vskip 5mm

\subsubsection{Kostant-Chevalley algebra and universal enveloping algebra}

We fix the same notations as in 5.1. In particular, $\U$ is the Kostant-Chevalley algebra of $\g$ over $\Z$. $\U$ can
be identified with the algebra $Dist(G)$ of distributions of $G$ (\cite{J}Part II,1.12). Recall that 
$$Dist(G) =\bigcup_{n\geq 0}(\Z[G]/{\cal M}^{n+1})^{*}$$
where $\cal M$ is the maximal ideal of regular functions vanishing at the unit element.
Let $\tilde\U$ be the universal enveloping algebra of $\g$. By the universal property of $\tilde\U$, we have a
natural homomorphism $\gamma : \tilde\U\to \U=Dist(G)$. $\gamma$ is injective and it is surjective over $\Z_p$ when
restricted to the $<p$-step of the filtrations of $\tilde\U$ resp. $\U=Dist(G)$:
$$\gamma: \tilde\U^{<p}\cong \U^{<p}
$$

 It will imply the following lemma:  

\begin{lem} Let $\U$ and $\tilde\U$ be the Kostant-Chevalley algebra and universal enveloping algebra
over $\Z_p$ respectively and $V_p$, $W_p$ be two $Q$-representations over $\Z_p$ whose semisimplifications 
have $p$-small highest weights
(a suficient condition on the highest weights is $\vert\lambda_i\vert<p$ ), then the canonical map 
$$Hom_\q(V_p,{\tilde{\cal
U}}\otimes _{{\tilde{\cal U}}({\q})} {W_p}) \to  Hom_\q(V_p,{\cal U}\otimes _{{\cal U}({\q})} {W_p}) $$ 
induced by $\gamma$, is an isomorphism.
\end{lem}

\noindent{\bf Proof:}

By PoincarŽ-Birkhoff-Witt over $\Z_p$, we have 
$${\tilde{\cal
U}}\otimes _{{\tilde{\cal U}}({\q})} {W_p}={\tilde{\cal
U}}\u^-\otimes _{\Z_p} {W_p}$$
where $\u^-$ is the unipotent radical of the  parabolic Lie algebra opposite of $\q$. 
It is enough to show 
$$Hom_\q(V_p,{\tilde{\cal U}}(\u^-)\otimes _{\Z_p} {W_p}) =
Hom_\q(V_p,{\tilde{\cal
U}}(\u^-){}^{<p}\otimes _{\Z_p} {W_p})
$$
Recall the decomposition of $W_p$
as a direct sum of
$\t$-eigenmodules $W_\sigma$ is valid over $\Z_p$ by diagonalizability of tori over any base.

For any $H\in \t$, $\underline{X}^{\underline{n}}\in\tilde{\cal U}(\u^-)$ ($\underline{n}=(n_\alpha)_{\alpha\in\Phi^{M-}}$) and $w\in W_\sigma$, we have
$$H\cdot (\underline{X}^{\underline{n}}\otimes w)=(\mu-\sum_{\alpha\in \Phi^{M+}}n_\alpha\alpha)\cdot
(\underline{X}^{\underline{n}}\otimes w)
$$
 For any $\q$-equivariant
$\phi:V_p\rightarrow {\tilde{\cal U}}(\u^-)\otimes _{\Z_p} {W_p}$, the image of a highest weight vector $v\in V_p$ is
of the form
$$\phi(v)=\sum_i\underline{X}_i^{\underline{n}_i}\otimes w_i\quad \hbox{\rm with}\, w_i\in W_{\mu_i}$$
Comparing the weights we have relations of the type
$$\lambda=\mu_i- \sum_{\alpha\in \Phi^{M+}}n^{(i)}_\alpha\alpha$$
by increasing the coordinates of $n^{(i)}$, we can assume that $\mu_i$ is the highest weight of $W_p$, hence is $p$-small.
Solving a linear system of inequations, we see that for any $\alpha\in \Phi^{M+}$,
$n^{(i)}_\alpha<p$ as desired.

\subsection{$p$-adic integral automorphic vector bundles} 

 Let $f:A\rightarrow X$ be the universal principally
polarized abelian variety over $X$ (with a $U$-level structure). Recall that
$R^1f_*\Omega_{A/X}^\bullet$ is endowed with the Gauss-Manin connection, which is integrable and quasi-nilpotent (see
Section 4.3). Let $\overline{X}$ be a toroidal compactification of $X$ over $\Z_p$.
 Let $\overline{X}_{n}=\overline{X}\otimes \Z/p^{n}\Z$; let $(\overline{X}\otimes \F_p/(\Z/p^{n}\Z))^{\rm log}_{\rm
cris}$ be the logarithmic crystalline site associated to the scheme
$\overline{X}\otimes \F_p$ and its divisor at infinity. Note that $\overline{X}\otimes \F_p$ is a toroidal
compactification of $X\otimes \F_p$.
 As recalled in Sect.4.1 above, there is an equivalence of category 
between crystals on this
site and locally free ${\cal O}_{\overline{X}_n}$-modules endowed with an integrable and
``{quasi-nilpotent}'' logarithmic connection.
Let ${\bf Rep}_{\Z_p}(G)$, resp. ${\bf Rep}_{\Z_p}(Q)$, be the category of algebraic representations of $G$, resp. $Q$, on
finitely generated free modules. Consider the respective full subcategories ${\bf Rep}_{\Z_p}^{\leq p-1}(G)$ and ${\bf
Rep}_{\Z_p}^{\leq p-1}(Q)$ consisting in objects whose highest weights are $p$-small (in fact, whose highest weights
$\mu$ satisfy $\vert \mu\vert\leq p-1$).

For each $n\geq 1$, let ${\bf \cV}_n^{\nabla}$, resp. $\overline{\bf \cV}_n^{\nabla}$ be the category of locally free ${\cal
O}_{X_n}$-modules, resp. ${\cal O}_{\overline{X}_n}$-modules, endowed with an integrable and ``{quasi-nilpotent}'',
resp. integrable, ``{quasi-nilpotent}''logarithmic connection, and 
${\bf \cF}_n$, resp. $\overline{\bf \cF}_n$ that of locally free ${\cal
O}_{X_n}$-modules, resp. ${\cal O}_{\overline{X}_n}$-modules endowed with a fitration
with locally free graded pieces. 

 The goal of this section is to define for each $n\geq 1$ two functors
$$\overline{V}_{\Z/p^n\Z}:{\bf Rep}_{\Z_p}^{\leq p-1}(G)\rightarrow \overline{\bf \cV}_n^{\nabla}$$ 
and another
$$\overline{F}_{\Z/p^n\Z}:{\bf Rep}_{\Z_p}^{\leq p-1}(Q)\rightarrow \overline{\bf \cF}_n$$ 

We first define functors on ${\bf Rep}_{\Z_p}(G)$, resp. ${\bf Rep}_{\Z_p}(Q)$ with values in vector bundles over $X_n$. Then 
we proceed to show that these vector bundles extend
 to
$\overline{X}_n$ provided they come from representations in ${\bf Rep}_{\Z_p}^{\leq p-1}(G)$ resp. ${\bf
Rep}_{\Z_p}^{\leq p-1}(Q)$.
\vskip 1cm

\subsubsection{``{Flat Vector bundles}'' on $X$}

Let us define
$${V}_{\Z/p^n\Z}: {\bf Rep}_{\Z_p}(G)\rightarrow {\bf \cV}_n^{\nabla}$$

Let $\cO_X^{2g}$ be the trivial vector bundle of rank $2g$ on $X$ endowed with the canonical
symplectic pairing (see section 5.1) and its natural action of
$G$ on the left. Let us put

$$\cT=
\underline{Isom}_{X}(\cO_{X}^{2g},(R^1 f_*\Omega^\bullet_{A/X})^\vee)$$
where the isomorphisms are symplectic similitudes.
It is an algebraic $G$-torsor over $X$ for the right action 
$$\cT\times G\rightarrow  \cT,\quad (\phi,g)\mapsto
\phi\circ g.$$
 For any
$V\in {\bf Rep}_{\Z_p}(G)$, we define $\cV$ as the contracted product 
$$\cV=\cT\stackrel{G}{\times} V$$
that is, the quotient of the cartesian product by the relation
$(\phi,g\cdot v)\sim (\phi\circ g, v)$. It is a vector bundle on $X$ hence over
$X_n$ for any
$n\geq 1$.
 
\vskip 5mm

\noindent{\bf Fact:} 1) $\cV$ is equipped with a connection of the desired type. 

2)  The image of the standard
representation is
$(R^1 f_*\Omega^\bullet_{A/X})^\vee$.

\noindent{\bf Proof:} 1) Let $\cA=R^1f_*\Omega_{A/X}^\vee$; we consider the (dual) Gauss-Manin connection:

$$\nabla: \cA\rightarrow  \cA\otimes_{\cO_X}\Omega_X$$
It is symplectic in the sense that for two sections $f,g$ of $\cA$, we have
$$<\nabla f,g>+<f,\nabla g>= d<f,g> $$ 
where the symplectic product is extended to 
$$\cA\otimes\cA\otimes \Omega_X\rightarrow \Omega_X$$
Therefore, given a point $\phi$ of $\cT$ over an $X$-scheme $Y$, we can transport $\nabla$ to an element $\nabla_\phi$ of
$\g\otimes
\Omega_X\subset End_{\cO_Y}(\cO_{Y}^{2g})\otimes_{\cO_X} \Omega_X$ defined by the diagram
$$\begin{array}{lcc}\cA_Y&\stackrel{\nabla}\rightarrow& \cA_Y\otimes\Omega_X\\\phi\downarrow
&&\downarrow\\
\cO_Y^{2g}&\stackrel{\nabla_\phi}\rightarrow &\cO_Y^{2g}\otimes \Omega^1_X\end{array}$$
Given $(V,\rho_V)\in {\bf Rep}_{\Z_p}(G)$, the representation $\rho_V$ viewed on the Lie algebra $\g$ enables us to define 
$$\nabla_{V,\phi}=(\rho_V\otimes Id_{\cO_Y}\otimes_{\cO_X}) Id_{\Omega_X} \circ \nabla \in End(V)\otimes \cO_Y\otimes_{\cO_X}
\Omega_X$$
It is a connection on
$ V\otimes
\cO_Y$. For $Y=\cT$, and $\phi$ the canonical point of $\cT$, we can descend this connection to the contracted product
because
$$\nabla_{\phi\circ h}=h^{-1}\circ \nabla_\phi \circ h$$
The resulting $\nabla_\cV$ is integrable and quasi-nilpotent because it is so for the Gauss-Manin connection.

2) Consider the morphism of $X$-schemes
$$\cT\times \cO_X^{2g}\rightarrow \cA,\quad (\phi,v)\rightarrow \phi(v)$$
It descends to the contracted product since $\phi\circ g (v)=\phi(g.v)$. It defines therefore a morphism of vector bundles over
$X$: $\cV_{st}\rightarrow \cA$. This morphism is an isomorphism over $\cT$ and $\cT\rightarrow X$ is faithfully flat, therefore
it is an isomorphism over $X$.

\subsubsection{Filtered Vector bundles on $X$}

The definition of the second functor
$${F}_{\Z/p^n\Z}:{\bf Rep}_{\Z_p}(Q)\rightarrow{\bf \cF}_n$$
is similar. We endow $\cO_{X}^{2g}$ with the standard symplectic pairing plus a Hodge filtration given by the Siegel parabolic
$Q$. 

$$\cT_H=
\underline{Isom}_{H,X}(\cO_{X}^{2g},(R^1 f_*\Omega^\bullet_{A/X})^\vee)$$
where the isomorphisms are symplectic similitudes respecting the Hodge filtrations.
$\cT_H$ is an algebraic $Q$-torsor over $X$. For any $W\in {\bf Rep}_{\Z_p}(Q)$, let 
$$\overline{\cW}=\cT_H\stackrel{Q}{\times} W$$
It is a vector bundle on $X$ hence over $X_n$ for any $n\geq 1$. It comes equipped with a filtration of the
desired type. The image of the standard representation is
$(R^1 f_*\Omega^\bullet_{A/X})^\vee$ with its standard filtration. The proof of these two assertions is very similar to the
one in the previous section.

\noindent{\bf Remark:} 1) In fact, by the same construction, one can define functors $V_{\Z[1/N]}$ and $F_{\Z[1/N]}$ such that
$V_{\Z/p^n\Z}=V_{\Z[1/N]}\otimes
\Z/p^n\Z$ and similarly for $F$.

2) Every $M$-representation gives rise to a $Q$-representation by letting the unipotent radical act trivially on $W$.

 \subsubsection{Comparison with the transcendental definitions}

Let $\widetilde{\cT}=G(\Q)\backslash  G(\A)\times G(\C)/UU_\infty$, the left action of
$G(\Q)$ on
$G(\A)\times G(\C)$ being diagonal, while the right one of $UU_\infty$ being only on the $G(\A)$-factor; the first
projection
$pr_1:G(\A)\times G(\C)\rightarrow G(\A)$ induces  a structure of principal
$G(\C)$-bundle over the analytic Siegel variety $S_U$  $\overline{pr_1}:\widetilde{\cT}\rightarrow S_U$. Moreover, let
$\check{\cZ}$  be the compact dual domain of the Siegel half-space $\cZ$. Let
$c\in GSp_{2g}(\C)$ be the standard Cayley matrix which defines the Cayley transform $\beta: \cZ\hookrightarrow \check{\cZ}$.
Consider the twisted multiplication
$$\mu:G(\A)\times G(\C)\rightarrow G(\C), (g,g')\mapsto g' c\cdot g_\infty\cdot c^{-1}$$ 
for
$g=(g_f,g_\infty)\in G(\A)$; it induces a morphism
$\overline{\mu}:\widetilde{\cT}\rightarrow 
\check{\cZ}$.

Recall the transcendental definition of
the automorphic vector bundle associated to $V\in {\bf Rep}_\C(Q)$: one forms the contracted product
$\check{\cV}=G(\C)\stackrel{Q(\C)}{\times} V$, which is a vector bundle over $\check{\cZ}$; then one
forms its pull-back
$\beta^*(\check{\cV})$ to
$\cZ$ by the Cayley transform
$\beta:\cZ\hookrightarrow
\check{\cZ}$. One takes the product $\beta^*(\check{\cV})\times G_f/U$, and one defines the 
holomorphic vector bundle 
$\widetilde{\cV}\rightarrow
S_U$ by
$$\widetilde{\cV}=G(\Q)\backslash(\beta^*(\check{\cV})\times G_f/U)\rightarrow G(\Q)\backslash
(\cZ\times G_f/U)=S_U.$$ We refer to $V\mapsto \widetilde{\cV}$ as the transcendental
construction. It is valid for $V\in {\bf Rep}_\C(Q)$ as well. 

The definition of $\widetilde{\cT}_H$ is slightly more subtle. We start from the $Q(\C)$-bundle
$\cQ:G(\C)\rightarrow \check{\cZ}$. We form its pull-back $\beta^*(\cQ)\rightarrow\cZ$ by $\beta$. It still carries an
equivariant action
of $G(\Q)$ on the left. We then form

$$\widetilde{\cT}_H=G(\Q)\backslash \beta^*(\cQ)\times G_f/U.$$
which is a $Q(\C)$-torsor over $S_U$.

\begin{lem} Over $\C$, the functor $V_{\C}$, resp. $F_{\C}$  is canonically isomorphic to the one
defined by the standard transcendental construction.
\end{lem}

\noindent{\bf Proof:} 1) Case of $V_\C$: we prove two statements

\begin{enumerate}
\item There is a canonical isomorphism of $G(\C)$-principal bundles $\widetilde{\cT}\rightarrow \cT$.
\item The transcendental construction can be described as
$$\widetilde{\cV}=\overline{pr_1}_*\circ \overline{\mu}^* \check{\cV}=\widetilde{\cT}\stackrel{G(\C)}\times V$$ 

\end{enumerate}

1. Let $V_{st}$ be the standard representation of $G$. We recall first that the pull-back by
$\cZ^{\prime}=\cZ\times G_f/U\rightarrow S_U$ of the vector bundle
$\cA$ endowed with the dual Gauss-Manin connection is isomorphic to the vector bundle of the local
system $\cZ^{\prime}\times V_{st}$ endowed with its obvious flat connection. Therefore, the pull-back
of $\widetilde{\cT}$ is isomorphic to 
$ \underline{Isom}_{\cZ^{\prime}}(\cZ^{\prime}\times V_{st},\cZ^{\prime}\times
V_{st})=\cZ^{\prime}\times G(\C)$, with action of
$G(\Q)$ diagonally on the left. Hence, by quotienting by $G(\Q)$, we obtain a canonical isomorphism
$\widetilde{\cT}\cong \cT$.

2. Let $V\in {\bf Rep}_\C(G)$. In this situation, only the $\cC^\infty$-structure of
$\widetilde{\cV}$ matters (indeed, only the structure of the underlying locally constant sheaf).  
On one hand, it is well-known that $\widetilde{\cV}$ is the vector bundle,  associated to the
$V$-covering
$G(\Q)\backslash (\cZ^{\prime}\times V)\rightarrow S_U$.
On the other hand, the pull-back by
$G(\C)\times\cZ^{\prime}\rightarrow
\widetilde{\cT}$ of
$\widetilde{\cT}\stackrel{G(\C)} \times V$ identifies to $\cZ^{\prime}\times V$; it is endowed 
with a free action of $G(\Q)$ (diagonally on the left), and of $U$ on the right. The resulting
quotient is again the vector bundle associated to the
$V$-covering
$G(\Q)\backslash (\cZ^{\prime}\times V)\rightarrow S_U$ as desired.

\vskip 5mm

2. Case of $F_\C$: From the definition of $\widetilde{\cT}_H$, it is clear that for any $V\in {\bf Rep}_\C(Q)$, 

$$\widetilde{\cV}=\widetilde{\cT}_H\stackrel{Q(\C)}\times V.$$

Moreover, there is a canonical isomorphism $\widetilde{\cT}_H\cong \cT_H$ of holomorphic $Q(\C)$-bundles. Indeed,
the pull-back by $\cZ^{\prime}\rightarrow S_U$ of $\cT_H$ 
$$\underline{Isom}_{\cZ^{\prime}}(\beta^*\cV_{st},\beta^*\cV_{st})= \beta^*\cQ\times G_f$$
hence, by quotienting, the desired isomorphism.

\subsubsection{$\Z_p$-Integral extension to $\overline{X}$ for $p$-small weights}

We have the diagram

$$(5.2.1) \quad \begin{array}{lcc}X_{\Q_p}&{\hookrightarrow} &X_{\Z_p}\\j\downarrow&\stackrel{k}\searrow&\downarrow\\

\overline{X}_{\Q_p}&
\stackrel{i}\hookrightarrow & \overline{X}_{\Z_p}
\end{array}$$

On one hand, for any $Q$-representation $W$, we have constructed a vector bundle $\cW$ over $X_{\Z_p}$; on the other hand, M.
Harris (\cite{Ha2}) has defined a functor from
$Q$-representations defined over $\Q$ to vector bundles over $\overline{X}_{\Q}$ coinciding with ours on $X_{\Q_p}$. 
We first glue the vector bundles $\overline{\cW}_{\Q_p}$ with $\cW_{\Z_p}$ into a vector bundle $\widetilde{\cW}_{\Z_p}$ over 
the cofibered product $\widetilde{X}_{\Z_p}=\overline{X}_{\Q_p}\cup_{X_{\Q_p}}X_{\Z_p}$.

Then, we observe that $\widetilde{X}_{\Z_p}=\overline{X}_{\Z_p}-D_{\F_p}$ is an open subset with complement of
codimension $2$ in $\overline{X}_{\Z_p}$. Therefore, by \cite{EGA4} Cor.5.11.4, the direct image of
$\widetilde{\cW}_{\Z_p}$  is a coherent sheaf on $\overline{X}_{\Z_p}$. Let us see it is locally free, at least
if $V$ has $p$-small highest weight. By dŽvissage, it is enough to consider irreducible $M$-representations
with such $p$-small highest weight. By Appendix II, it is enough to consider the standard representation. In
that case, the coherent sheaf on $\overline{X}_{\Z_p}$ is
$Lie(\cG/\overline{X})^\vee$, which is locally free. This concludes the proof.

\vskip 1cm

In particular, for any dominant weight
$\lambda$, we have attached to the representation $V_\lambda$ of $G$ of highest weight $\lambda$ 
 a vector bundle ${\cal
O}_{\overline{X}_n}$-module 
$\overline{\cV}_{\lambda,n}$ on
$\overline{X}_n$ together with a
 connection with logarithmic poles along $D_n$, hence a
logarithmic crystal 
 ${\overline\cV}_{\lambda,n}$ on $(\overline{X}/(\Z/p^{n}\Z))^{\rm log}_{\rm cris}$. Moreover, it carries a natural filtration
since $V_\lambda$ is also a $Q$-representation.

\subsubsection{Differential operators over $\Z_{(p)}$}
Let $V$ and $W$ be two rational representations
of $Q$, and $\cal V_{/\Q}$, $\cal W_{/\Q}$ the corresponding automorphic
vector bundles over
$X_\Q$ (see previous subsection) and $\overline{\cal V}_{/\Q}$,
$\overline{\cal W}_{/\Q}$ their canonical
extension to the toroidal compactification $\overline X$. According to
Proposition 5.1 of
\cite{CF} VI.5, we have a functorial homomorphism
$$\Psi :  Hom_{ U(\g_\Q)}( U(\g_\Q)\otimes_{\U(\q_\Q)} V,{ U(\g_\Q)}\otimes _{\U(\q_\Q)} {W}) \to
{\rm Diff. Operators} (\overline{\cal W}^{\vee}_{/\Q},\overline{\cal
V}^{\vee}_{/\Q}).$$

Actually, in Proposition 5.1 of Chapt.VI, the construction of $\Psi$ is
explained over $\bC$. The
$\Q$-rationality statement is explained in Remark 5.2 following the Proof
of Proposition 5.1 of Sect.
VI.5. We now prove a variant thereof over $\Z_{(p)}$.

We treat first the case of degree $0$ differential operators by refering to 5.2.2:

\begin{lem} Let  $V$, $W$ be two $Q$-representations of
$p$-small highest weights (in fact, $\vert\lambda_V\vert$ and $\vert\lambda_W\vert <p$ is enough), $V_p$
and $W_p$ their canonical $\cal U$-stable lattices  and $\overline{\cal
V}_n$, $\overline{\cal W}_n$
the corresponding automorphic vector bundles over $\overline{X}_{n}$,
$n>0$. There is a functorial injective homomorphism

$$Hom_{\q}( V_p, W_p) \to Hom_{{\cal O}_{{\overline X}_n}}(\overline{\cal
W}^{\vee}_n,\overline{\cal V}^{\vee}_n)$$

compatible with the transcendental construction. 
 \end{lem}

Then, the case of general differential operators can be treated as follows

\begin{lem} Let  $V$, $W$ be two irreducible $Q$-representations of
$p$-small highest weights, $V_p$
and $W_p$ their canonical $\cal U$-stable lattices  and $\overline{\cal
V}_n$, $\overline{\cal W}_n$
the corresponding automorphic vector bundles over $\overline{X}_{n}$,
$n>0$. Then  $\Psi$  induces for
each $n>0$, a homomorphism
$$Hom_{\cal U}({\cal U}\otimes_{{\cal U}({\q})} V_p,{\cal U}\otimes _{{\cal
U}({\q})} {W_p}) \to {\rm P.D.Diff. Operators} (\overline{\cal
W}^{\vee}_n,\overline{\cal V}^{\vee}_n)$$
 \end{lem}

{\bf Proof:} We start with operators of order one. 
Note that the de Rham differential $d :
{\cal O}_{\overline{X}_n}\to \Omega^{1}_{\overline{X}_n}$ is the image by
$\Psi$ of the obvious map
$\delta: \U\otimes_{\U(\q)}\g_{\bZ_p}/\q_{\bZ_p}\to  \U\otimes_{\U(\q)}\bZ_p$
(compare with \cite{CF} VI, remark 5.2). By Lemma 10, this implies that each homomorphism
$\phi: V_p \to {\cal U}\otimes _{{\cal U}({\q})} {W_p}$  of degree one is mapped by $\Psi$ to a $\Z_p$-integral differential operator 
of order one. Indeed any $\phi$ as above factors as $\phi=\delta\otimes Id_{W_p}\circ (Id_{\tilde{\U}}\otimes\psi)$ for a $\psi\in
Hom_\q(V_p,\g/\q\otimes W_p)$.

Recall that $\tilde{\cal U}$ denotes the universal enveloping algebra of $\g$.  We have seen in Lemma 8 that
by $p$-smallness of the highest weights, the natural algebra homomorphism $\gamma:{\tilde{\cal U}}\to{\cal U}$
 induces a bijection between
$Hom_\g(V_p,{\cal U}\otimes _{{\cal U}({\q})} {W_p})$ and $Hom_\g(V_p,{\tilde{\cal U}}\otimes _{{\tilde{\cal
U}}({\q})} {W_p})$. Now, as a corollary of PBW over $\Z_p$ for $\tilde{\U}$, we see that every element $\phi\in Hom_\g(V_p,{\tilde{\cal U}}\otimes _{{\tilde{\cal
U}}({\q})} {W_p})$ of degree $m>1$ factors as $\phi=\delta\otimes Id_{W_p}\circ \psi$ where $\psi$ has degree $m-1$:
fix a basis $(X_\alpha)_{\alpha\in \Phi^{M-}}$ of $\u^-$; for $v\in V_p$ and $\phi (v)=\sum_{i} \underline{X}^{\underline{n}^{(i)}}\otimes w_i$, put  
 $\psi (v)=\sum_{i}\sum_{\alpha\in\Phi^{M-}} \underline{X}^{\underline{n}^{(i)}}-1_\alpha\otimes X_\alpha\otimes w_i$
where $1_\alpha$ is the family $(\delta_{\alpha,\beta})_{\beta\in\Phi^{M-}}$.
The conclusion follows by induction on $m$.

\subsection{The Hodge filtration on automorphic sheaves}

\subsubsection{The geometric aspect}

This paragraph is a recollection of well-known facts about the Hodge filtration in the
automorphic setting (see \cite{Del} Sect.5).

 Let ${\underline S}=R_{\C/\R}{\bf G}_{m}$ and $h_{0} : {\underline S}(\R)\to
G(\R)$ the homomorphism defined by 
$$z=x+iy\in\C^{\times}\to
\left(
\begin{array}{cc}
xI_{g} & yI_{g}\\
-yI_{g} & xI_{g}
\end{array}
\right)
=xI_{2g}+yJ_{2g}\in G(\R)$$
The $G(\R)$-orbit $\cZ$ of $h_{0}$ is analytically isomorphic to a double copy of the Siegel upper
half-plane  of genus $g$. The pair $(G,\cZ)$ defines a family of Shimura varieties "\`a la Deligne",
isomorphic to our  Shimura varieties $S_{U}$ for various level structures $U$. If $V$ is a real
representation of
$G$ and $h\in X$, then the composition $h : {\underline S}(\R)\to G(\R)\to GL(V)$ defines a real
Hodge structure $h_{V}$ on $V$ (\cite{Del}). Let $F_{h}$ be the filtration on
$V_{\C}=V_{\R}\otimes\C$ deduced from $h_{V}$.  For $V={\g}$ the adjoint representation,
$F^{0}_{h}({\g}_{\C})$ is a Lie algebra of a parabolic subgroup $P(h)$ of $G_\C$. The mapping
$h\to P(h)$ identifie $\cZ$ as an open subset of its compact dual ${\breve \cZ}=G(\C)/Q(\C)$. Now, for
general $V$, the mapping $h\to F_{h}$ define a $G(\C)$-equivariant filtration (the Hodge filtration)
on the constant fibre bundle $\cZ\times V_{\C}$. Dividing by $G(\Q)$ and $U$, we get a
filtration on the coherent sheaf $\cal V$ over $S_U$, associated to the representation $V$.
Moreover, if $\overline {\cal V}$ is the canonical extension of $\cal V$ to some toroidal
compactification of $S_U$, then this filtration has a canonical extension to $\overline{\cal V}$ (cf. 5.2.4).
In the case where $V$ is the standard representation of $G$, then, as explained above, 
$\overline{\cV}^\vee=R^1\bar{f}_*\Omega^\bullet_{\bar{A}/\overline{X}}(log\,\infty_{\bar{A}/\overline{X}})$
and the Hodge filtration on the dual is the classical one given by

$$(5.3.1)\qquad F^{2}(\overline{\cV}^\vee)=0\subset  F^{1}(\overline{\cV}^\vee) = 
\bar{f}_*\Omega^{1}_{\bar{A}/\overline{X}}(log\,\infty_{\bar{A}/\overline{X}})\subset
F^{0}(\overline{\cV}^\vee)=\overline{\cV}^\vee$$

Then, for a represention $V_\lambda$ associated to a dominant weight $\lambda$ 
of $G$, we can use Weyl's invariant theory as in
Appendix II, to describe the Hodge filtration on $\overline{\cV}_\lambda^\vee$. 

More precisely, we use the fact that each $\overline{\cV}_\lambda^\vee$ is a direct summand of some
higher direct image of the logarithmic de Rham complex over a toroidal compactification of the
 $s$-fold product of the universal abelian variety \cite{CF}p.234. 

Recall that for a complex $K^{\bullet}$, the notation $K^{\bullet \geq  i}$ denotes the subcomplex of
$K^{\bullet}$ equal to
$K^{\bullet}$ in degre $\geq  i$ and zero elsewhere. 

If ${\overline
f}_s:\overline{Y}\to {\overline X}$ is such a toroidal compactification, then the coherent sheaf 
$${\cal F}=R^\bullet {\overline f}_{s\,*}\Omega^\bullet_{{\overline Y}/{\overline X}}$$ 
is endowed
with the Hodge filtration 
$$Fil^i{\cal F}=R^\bullet {\overline f}_{s\,*}\Omega^{\bullet,\geq i}_{{\overline Y}/{\overline
X}}$$

Then, for a dominant weight $\lambda$ such that $\vert \lambda\vert=s$ 
we endow the sheaf $\overline{\cV}_\lambda$ 
by the filtration : 
$$Fil^i \overline{\cV}^\vee_\lambda={\overline \cV}^\vee_\lambda\cap Fil^i{\cal F}$$

In particular, if the highest weights of the
representation $V$ are $p$-small and as in Sect.5.2, this filtration is algebraic and has a model over $\Z_p$.
 Finally, we note that this filtration is compatible with tensor product, duality, etc.\\   

Let $\overline{\cV}^\vee_{\lambda,n}$ be the ${\cal
O}_{\overline{X}_n}$-module obtained by reduction mod $p^n$ of the module 
$\overline{\cV}^\vee_{\lambda}$.

\begin{de} The Hodge filtration on the de Rham complex
$${\overline\cV}^\vee_{\lambda,n}\otimes_{{\cal
O}_{\overline{X}_n}}\Omega^{\bullet}_{\overline{X}_n/\bZ/p^n}(\log\infty)$$

 is defined by:

$$F^{i}({\overline\cV}^\vee_{\lambda,n}\otimes_{{\cal
O}_{\overline{X}_n}}\Omega^{\bullet}_{\overline{X}_n/\bZ/p^n}(\log\infty))=\sum_{j} 
F^{j}({\overline\cV}^\vee_{\lambda,n})\otimes_{{\cal
O}_{\overline{X}_n}}\Omega^{\bullet}_{\overline{X}_n/\bZ/p^n}(\log\infty)^{\geq i-j}$$
\end{de}

\subsubsection{The group-theoretic aspect}

Let $H=\hbox{diag}(0,\ldots,0,-1,\ldots,-1)\in\hbox{Lie}T\subset\g$  (with $g$
$0$'s and $g$ $-1$'s).  For any rational $Q$-representation $V$, for any $i\in \Z$, let $V^i$ be
 the sum of the generalized $H$-eigenspaces with
eigenvalues $\geq  i$. This defines a decreasing filtration
$\{V^{i}\}$ on $V$.
Two cases are of particular interest for us:
\begin{itemize}
\item  $V$ is an irreducible 
$M$-representation with highest weight $\mu$; the filtration  is given by 
$V^{\mu(H)+1}=0\subset V^{\mu(H)}=V$. For instance, the standard representation $V_0$ of $M$ is filtered by
$0=V_0^1\subset V_0^{0}=V_0$ while its twisted contragredient $V_1=V_0^\vee\otimes \nu$ is filtered by
$0=V_1^0\subset V_1^{-1}=V_1$. 

\item  $V=V_{\lambda}$ is an irreducible 
representation of $G$ associated to the dominant weight $\lambda$. Then the filtration given by $H$
can also be defined by
plethysms from the $2$-step filtration of the standard representation $V_{st}$:
$F^{-1}=V_{st}$, $F^0=V_0$ is its unique simple $Q$-submodule (in fact, an $M$-module), and $F^1=0$. 

\end{itemize}

\noindent {\bf Fact:} In the construction $V\mapsto \cV$ of the coherent sheaf attached to a $Q$-representation,
the $H$-filtration defined above gives rise to a decreasing filtration on $\cV$. When $V$ is a $G$-representation, 
it coincides with the Hodge filtration given by $F_{h_{0}}$.

\noindent{\bf Proof:} Consider the dual filtration
$$(5.3.2.1)\qquad Fil^i\,\cV^\vee=\{\varphi:\cV\rightarrow \cO_X; \varphi(Fil^j \cV)\subset Fil^{i+j}
\cO_X\}$$ where the unit object $\cO_X$ is endowed with the trivial filtration: $Fil^0\cO_X=\cO_X$
and $Fil^j\cO_X=0$ for any $j>0$; When $V$ is the complex standard
representation $V_{st}\otimes\C$ of $G_\C$,  the dual of the
$H$-filtration coincides with the Hodge filtration (given by $F_{h_{0}}$) on $\cV^\vee$,
 indeed, the dual of the $H$-filtration reads:

$$(5.3.2.2)\qquad Fil^0\cV^\vee=\{\varphi; \varphi(Fil^1\cV)=0\}=\cV^\vee,$$
$$
Fil^1\cV^\vee=\{\varphi; \varphi(Fil^0\cV)=0\}=\cV_1^\vee,\quad {\rm and}\,
Fil^2\cV^\vee=0.$$
This is the Hodge filtration (5.3.1).

\vskip 5mm  

We can still define the $H$-filtration as above for a $Q$-representation $V$ defined over $\Z_p$ instead of
$\C$; the $V^i$'s are $\Z_p$-summands in $V$.

In particular, we endow the
standard bar resolution of $V_{\lambda,\,\Z_p}$ (say, for $\vert\lambda+\rho\vert<p-1$)
$$D(\lam) := (\U_{\Z_p}
\otimes_{\U(\q)_\Zp} (\Lambda^\bullet(\g/\q)\otimes V(\lam)_\Zp))$$

with the $H$-filtration.

Let 
$$C(\lam)_{\Z_p}=\bigoplus_{w\in W^M}\U\otimes_{\U(\q)} W_{w(\lam+\rho)-\rho,\,\Z_p}$$
 be the BGG
complex introduced in Sect.5.1.2 attached to $V_{\lambda,\,\Z_p}$. the $H$-filtration is given by 

$$F^{i}C(\lambda)_{\Z_p}=\bigoplus_{w\in W^{M},w(\lambda+\rho)(H)-\rho(H)\geq i}
\U\otimes_{\U(\q)} W_{w(\lam+\rho)-\rho ,\Z_p} $$

Then the injection $j : C(\lam)_\Zp \hookrightarrow D(\lam)_\Zp$ is a filtered direct factor of
$D(\lam)_\Zp$ by \cite{PT}. 

Let us define the dual BGG complexes $\overline{\cK}_{\lambda,n}^\bullet$ and
$\overline{\cK}_{\lambda,n}^\bullet$. Their graded pieces are the coherent sheaves over
$\overline{X}_n$:

$$\overline{\cK}_{\lambda,n}^{i} = \bigoplus_{w\in W^{M},l(w)=i}\overline{\cal
W}_{w(\lambda+\rho)-\rho,n}^\vee $$

resp.

$$\overline{\cK}_{\lambda,n}^{i,sub} = \bigoplus_{w\in W^{M},l(w)=i} \overline{\cal
W}^{sub,\vee}_{w(\lambda+\rho)-\rho,n}$$

and the differentials are deduced by lemma 11 (sect. 5.2.5) from the BGG complex of Sect.5.1.2.
By dualizing the $H$-filtration, we obtain a natural  decreasing
filtration on $\overline{\cK}_{\lambda,n}^\bullet$, stable by the differentials, given by

$$F^{i}\overline{\cK}_{\lambda,n}^{\bullet}=\bigoplus_{w\in W^{M},w(\lambda+\rho)(H)+i\leq \rho(H)}
\overline{\cal W}^{\vee}_{w(\lambda+\rho)-\rho,n}$$

Recall that by the Theorem of \cite{PT}, the map $j$ has a retraction of filtered complexes, hence
the dual $j^\vee$ has a natural section; its
sheafification defines an injection of complexes of coherent
${\cal O}_{\overline{X}_n}$-modules:

$$\kappa : \overline{\cK}_{\lambda,n}^{\bullet} = \bigoplus_{w\in W^{M}} \overline{\cal
W}^{\vee}_{w(\lambda+\rho)-\rho,n}  \hookrightarrow {\overline\cV}^{\vee}_{\lambda,n}\otimes_{{\cal
O}_{\overline{X}_n}}\Omega^{\bullet}_{\overline{X}_n/\bZ/p^n}(\log\infty) $$

$$\kappa : \overline{\cK}_{\lambda,n}^{\bullet, sub} = \bigoplus_{w\in W^{M}} \overline{\cal
W}^{\vee,sub}_{w(\lambda+\rho)-\rho,n}  \hookrightarrow
{\overline\cV}^{\vee}_{\lambda,n}\otimes_{{\cal
O}_{\overline{X}_n}}\Omega^{\bullet}_{\overline{X}_n/\bZ/p^n}(-\log\infty) $$

We summarize the considerations of this section in the proposition

\begin{pro} The morphism $\kappa$ of complexes of vector bundles over $\overline{X}_n$ ($n\geq 1$) is
filtered.
\end{pro}

\subsection{BGG resolution}

 We denote by
${\cal D}_n$ the logarithmic divided power envelope of the diagonal immersion $\overline{X}_n\to
\overline{X}_n{\hat\times}_{\Z/p^n}\overline{X}_n$ where 
$\overline{X}_n{\hat\times}_{\Z/p^n}\overline{X}_n$ is the fiber product in the category of
logarithmic schemes. Let $p_1$ and $p_2$ be the two canonical projections ${\cal D}_n\to
\overline{X}_n$. Finally, for any $B_M$-dominant weight $\mu$ of $M$, such that
$\vert\mu\vert<p$, let
$L(\overline{\cal W}_{\mu,n})$ be the logarithmic crystal on
$(\overline{X}/\Z/p^{n})^{\rm log}_{\rm cris}$ corresponding to $p_{1}^{*}{\overline{\cal W}}_{\mu,n}$ (Sect.
4.2 for $L$ and 5.2.    for ${\overline{\cal W}}_{\mu,n}$). 
For simplicity, in the sequel, we drop the index $n$ in the notations of the sheaves, thus
we write ${\cal W}_\mu$ for ${\cal W}_{\mu,n}$. Note that we can not consider the situation 
over $\bZ_p$ because we need a
nilpotent base for our crystalline arguments.

\begin{pro} 
Let $\lambda$ be a $B$-dominant weight of $G$, such that $\vert\lambda+\rho\vert<p$;

(i) There is a resolution in the category of logarithmic crystals on\\
$(\overline{X}_0/(\Z/p^n\Z))_{crys}^{log}$:
$$ 0\to {\overline{\cV}^{\vee}_\lambda} \to L(\overline{\cK}_\lambda^{0}) \to 
L(\overline{\cK}_\lambda^{1})
\to  
\ldots
$$ where 
$$\overline{\cK}_\lambda^{i} = \bigoplus_{w\in W^{M},l(w)=i} \overline{\cal
W}^{\vee}_{w(\lambda+\rho)-\rho}.$$

(ii) There is a canonical filtered quasi-isomorphism of complexes of logarithmic crystals

$$L(\overline{\cK}_\lambda^{\bullet}) \to L(\overline{\cV}^{\vee}_{\lambda}\otimes_{{\cal
O}_{\overline{X}_n}}\Omega^{\bullet}_{\overline{X}_n/\Z/p^n}(\log\infty))$$
\end{pro} 

{\bf Proof:} We transpose the proof given in \cite{CF}, VI,Sect.5 for the complex case in a
$\Zp$-setting. 

By Lemma 11, each $\g_{\Z_{(p)}}$-morphism of order
$1$: 
$$\U\otimes_{\U(\q)}W_1\to \U\otimes_{{\cal U}(\q)}W_2$$
 induces a logarithmic differential
operator of order $1$, $\overline{\cal W}_2^\vee\to
\overline{\cal W}_1^\vee$ for the corresponding locally free
${\cal O}_{\overline{X}_n}$-module; therefore, it induces a morphism of crystals $L(\overline{\cal
W}_2^\vee)\to L(\overline{\cal W}_1^\vee)$. We deduce from theorem 5 (section 5.1.2), that there is a complex of
crystals  

$$0\to
{\overline{\cV}^{\vee}_\lambda} \to L(\overline{\cal K}_\lambda^{0})
\to L(\overline{\cal K}_\lambda^{1})
\to   \ldots. $$
 
On the other hand, we know that
$$0\to {{\overline\cV}^{\vee}_\lambda} \to 
L({\overline\cV}^{\vee}_{\lambda}\otimes_{{\cal
O}_{\overline{X}_n}}\Omega^{\bullet}_{\overline{X}_n/\Z/p^n}(\log\infty))$$

is a resolution of ${{\overline\cV}^{\vee}_\lambda}$. Indeed, the exactness of the complex is the
crystalline Poincar\'e's lemma (actually, its logarithmic version: bottom of p.221 of
\cite{Ka}, see our section 4.2, lemma 4 above).

By Theorem D of \cite{PT} (Theorem 5 of section 5.1.2 here),
$\overline{\cal K}_\lambda^{\bullet}$ is a subcomplex of
${\overline\cV}^{\vee}_{\lambda}\otimes_{{\cal
O}_{\overline{X}_n}}\Omega^{\bullet}_{\overline{X}_n/\Z/p^n}(\log\infty)$. Therefore, $\overline{\cal
K}_\lambda^{\bullet}$ is  a resolution of
${{\overline\cV}^{\vee}_\lambda}$.  
This proves statement (i) of the theorem. The second assertion follows from the fact that
$H$ commutes with $Z\g$. As explained in Section 5.1.2 above.

\begin{thm}  The natural morphisms 

$$\overline{\cK}_{\lambda}^{\bullet}\to \overline{\cV}_{\lambda}^\vee\otimes_{{\cal
O}_{\overline{X}_n}}\Omega^{\bullet}_{\overline{X}_n/\Z/p^n}(\log\infty)$$ 
and

$$\overline{\cK}_{\lambda}^{\bullet,sub}\to {\overline\cV}_{\lambda}^\vee\otimes_{{\cal
O}_{\overline{X}_n}}\Omega^{\bullet}_{\overline{X}_n/\Z/p^n}(-\log\infty)$$

are filtered quasi-isomorphisms of complexes of coherent sheaves on
$\overline{X}_n$.

\end{thm}

{\bf Proof:} One applies $ Ru_*$ resp. $Ru_{*,c}$ to both members of the quasi-isomorphism (ii) of
Prop. 4; then one makes use of the fact that $Ru_* L(\cV)\cong \cV$ for any
$\cO_{\overline{X}_n}$-module $\cV$ and the properties of the $L$-construction recalled in Section 4.2.

\vfill\eject

\section{Modulo $p$ crystalline representations}

\subsection{Etale sheaves associated to crystals}

Let $k$ be a perfect field of char. $p>0$, $W=W(k)$ the ring of Witt
vectors with coefficients in $k$
and $K$ the fraction field of $W$. $K^{ac}$ is a fixed algebraic closure of
$K$ and $G_{K}=Gal(
K^{ac}/K)$ is the associated Galois group. Let  $Rep_{\Z_p}(G_K)$ be the
category of
$G_K$-modules of finite type over $\Z_p$ and $MF_{W}^{[0,p-2]}$ that
of finitely generated W-modules $M$ endowed with a filtration
$(Fil^{r}M)_r$ such that
$Fil^{r}M$ is a direct factor, $Fil^{0}M=M$ and
$Fil^{p-1}=0$ together with semi-linear maps $\varphi_r : Fil^{r}M\to M$
such that the restriction
of
$\varphi_r$ to $Fil^{r+1}M$ is equal to $p\varphi_{r+1}$ and satisfying the
strong divisibility
condition : $M=\sum_{i\in\Z}\varphi_{r}(Fil^{r}M)$. Recall that by the
theory of Fontaine-Laffaille
\cite{FoL}, we have a fully faithful functor

$$V_{\rm cris} : MF_{W}^{[0,p-2]}\longrightarrow Rep_{\Z_p}(G_K)$$

A $p$-adic representation is called
crystalline if it is in the essential image of $V_{\rm cris}$.

In our setting, we are interested in the subcategory $MF_{k}^{[0,p-2]}$ of
filtered modules $M$ such
that $pM=0$. $MF_{k}^{[0,p-2]}$ is an abelian category and the objects are
in particular $k$-vector
spaces. The restriction of the functor $V_{\rm cris}$ to $MF_{k}^{[0,p-2]}$
can be describe as
follows: Let $S={\cal O}_{K^{ac}}/p{\cal O}_{K^{ac}}$, choose $\beta\in
K^{ac}$ such that
$\beta^{p}=-p$ and for
$i<p$, define a filtration $Fil^{i}S=\beta^{i}S$ and Frobenius
$\varphi^{i}(\beta^{i}x)=x^{p}$, then
as explained in \cite{Wa}, Prop.2.3.1.2', we have an isomorphism
$$V_{\rm cris}(M)\simeq Hom_{MF^{[0,p-2]}_{k}}(M,S)^*$$
Moreover, $V_{\rm cris}(M)$ is a finite dimension $\F_p$-vector space and
${\rm dim}_{\F_p}V_{\rm
cris}(M)={\rm dim}_{k}M$. Recall that
$V_{\rm cris}$ is the nice inverse of a not so nice contravariant Dieudonn\'e
functor
$\bD^*$: see \cite{Wa} p.219-223.

Let $\overline{X}$ be a smooth and proper scheme over $W$ of relative dimension
$d$ and $D$ a relative
divisor with normal crossings of
$\overline{X}$, we put $X=\overline{X}-D$. Faltings introduced in \cite{Fa}
relative
versions of the categories mentioned above: the category
${\cal R}ep_{\Z_{p}}(X\otimes K)$ of \'etales
$\Z_p$-sheaves over the generic fiber $X\otimes K$ and the category ${\cal
MF}^{\nabla}(\overline{X})$ of
filtered transversal logarithmic crystals over
$\overline{X}$. Moreover, we have a notion of "associated" between objects
of  ${\cal
R}ep_{\Z_{p}}(X\otimes K)$ and those of
${\cal MF}^{\nabla}(\overline{X})$.
To get a good theory over $\Z_p$, we need to consider only the full
subcategory
${\cal MF}^{\nabla, [0,p-2]}(\overline{X})$ of ${\cal
MF}^{\nabla}(\overline{X})$ of filtered crystals
$\cal F$ such that
$Fil^{0}{\cal F}={\cal F}$ and $Fil^{p-1}{\cal F}=0$ and we have to add
some other technical
hypothesis (cf Sect.4.2). Faltings \cite{Fa} (see also \cite{Tsu}) has defined a relative Fontaine
functor
$${\bf V}^* : {\cal MF}^{\nabla, [0,p-2]}(\overline{X}) \longrightarrow {\cal
R}ep_{\Z_{p}}(X\otimes
K)$$
In section 4.2 below, we will recall its definition in the special case where $\overline{X}$ is
actually of characteristic $p$. 

\begin{de}
For any ${\cal F}\in {\cal MF}^{\nabla, [0,p-2]}(\overline{X})$, we say that $\cal F$ and
${\bf V}^*({\cal F})$ are associated.
\end{de}

We have the following  theorem of Faltings
(\cite{Fa}Th.5.3) :
\begin{thm} Assume that there is a positive integer $a$ such that $a+d\leq p-2$
and $Fil^{a}{\cal F}=0$
where ${\cal F}\in {\cal MF}^{\nabla, [0,p-2]}(\overline{X})$ then, for $i=0,\ldots,2(p-2)$, we have
a natural and functorial
isomorphism of
$G_K$-modules:
$$H_{\rm et}^{i}(X\otimes K^{ac},{\bf V}^*({\cal F}))\cong V_{\rm
cris}^*(H^{i}_{\rm log-cris}({\bar
X},{\cal F}))$$
\end{thm}

\subsection{The mod $p$ case}
As we use only the mod $p$ version of the previous comparison theorem,
we only recall the notion of
associated sheaves and the comparison theorem in their mod.$p$ version, following \cite{Fa} and
\cite{Tsu}.

\subsubsection{Filtered modules}
Let $k$ be a perfect field of char. $p>0$, $W=W(k)$ the ring of Witt
vectors with coefficients in $k$
and $K$ the fraction field of $W$. $K^{ac}$ is a fixed algebraic closure of
$K$ and $G_{K}=Gal(
K^{ac}/K)$ is the associated Galois group.

Let $\overline{X}$ be a smooth and proper scheme over $W$ of relative
dimension $d$ and $D$ a relative
divisor with normal crossings of $\overline{X}$, we put $X=\overline{X}-D$. Let
$\overline{X}_{0}=\overline{X}\otimes_{W} k$ be the special fiber  of
$\overline{X}$ and $D_{0}$ the
induced divisor. If $F_{X_{0}} : {\cal O}_{\overline{X}_{0}} \to {\cal
O}_{\overline{X}_{0}}$ is the
absolute Frobenius, we denote by
$$\varphi_{\overline{X}_{0}} : F_{\overline{X}_{0}}^{-1}({\cal
O}_{\overline{X}_{0}})\to {\cal O}_{\overline{X}_{0}}$$
the ${\cal O}_{\overline{X}_{0}}$-linear homomorphism induced by
$F_{\overline{X}_{0}}$.

\begin{de}

We define the
category
${\cal MF}^{\nabla, [0,p-2]}_{k}({X_{0}})$ of strongly divisible filtered
logarithmic modules over
$\overline{X}_{0}$ with Hodge-Tate weights between
$0$ and $p-2$ as follows : an object is a quadruple
$({\cal F},{\cal F}^{i},\varphi^{i}_{\cal F},\nabla_{\cal F})$ where \\
$\bullet$  $\cal F$ is a quasi-coherent ${\cal O}_{\overline{X}_{0}}$-module.\\
$\bullet$ ${\cal F}^{i}$, $i=0,\ldots,p-1$, is a decreasing filtration of
$\cal F$ by quasi-coherent
${\cal O}_{\overline{X}_{0}}$-modules such that ${\cal F}^{0}={\cal F}$ and
${\cal F}^{p-1}=0$.\\
$\bullet$ $\varphi^{i}_{\cal F} : {\cal F}^{i} \to {\cal F}$ is a
$\varphi_{\overline{X}_{0}}$-linear
homomorphism such that the restriction of $\varphi^{i}_{\cal F}$ to ${\cal
F}^{i+1}$ is zero and
such that the induced map $$\oplus_i \varphi^{i}_{\cal F} : \oplus {\cal F}^{i}/{\cal F}^{i+1}
\to{\cal F}$$
is an isomorphism (condition of strong divisibility).\\
$\bullet$ $\nabla_{\cal F} : {\cal F}\to {\cal F}\otimes_{{\cal
O}_{\overline{X}_{0}}}\Omega^{1}_{\overline{X}_{0}}(log D_{0})$ is  a
quasi-nilpotent
integrable connection satisfying

 1) Griffiths transversality :
$\nabla_{\cal F}({\cal F}^{i})\subset {\cal F}^{i-1}\otimes_{{\cal
O}_{\overline{X}_{0}}}\Omega^{1}_{\overline{X}_{0}}(log D_{0})$ for
$i=0,\ldots,p-1$.

2) Compatibility with Frobenius :
$\displaystyle{\nabla_{\cal F}\circ \varphi^{i}_{\cal
F}=\varphi^{i-1}_{\cal F}\otimes
{d\varphi_{\overline{X}_{0}}\over p}\circ \nabla_{\cal F}\vert{\cal F}^{i}}$.

Where ${d\varphi_{\overline{X}_{0}}\over p}$ is by definition, the
reduction mod $p$ of
${d{\tilde\varphi}_{\overline{X}_{0}}\over p}$ for some lifting (wich exist
locally)
${\tilde\varphi}_{\overline{X}_{0}}$ over $W$ of $\varphi_{\overline{X}_{0}}$.

$\bullet$ $\cal F$ is uniform : \'etale locally, there is a log-immersion
$\overline{X}_{0}\to
\overline{Z}$ such that
$$({\cal F},{\cal F}^{i})\simeq \oplus_{\lambda\in\Lambda}({\cal O}_{{\overline
Z}^{DP}},J^{[i-e_{\lambda}]}_{{\overline Z}^{DP}}) \hbox{ with } e_{\lambda}\geq  0,
\vert\Lambda\vert<\infty$$
where ${\overline Z}^{DP}$ is the log-divided power envelope of the
immersion $\overline{X}_{0}\to
\overline{Z}$, $J_{{\overline Z}^{DP}}$ the coresponding PD-ideal.

\end{de}

A morphism of ${\cal MF}^{\nabla, [0,p-2]}_{k}({X_{0}})$ is an ${\cal
O}_{\overline{X}_{0}}$-linear
homomorphism compatible with filtrations and commuting with Frobenius and
connections.

\vskip 5mm

By \cite{Fa},Th.2.1, each ${\cal F}^{i}$ is locally free and locally (for
the Zariski topology) a
direct factor of $\cal F$. Moreover, any morphism of  ${\cal MF}^{\nabla,
[0,p-2]}_{k}({X_{0}})$ is strict
with respect the filtrations. We deduce from this that ${\cal MF}^{\nabla,
[0,p-2]}_{k}({X_{0}})$ is an
abelian category.

\subsubsection{The functor $\bf V^*$}
To a filtered module $\cal F$ as above, we associate an \'etale sheaf ${\bf
V}({\cal F})$ over $X\otimes
K$ as follows: \\
Let $\overline{U}=Spec(R)$ be an affine open irreducible subset of
$\overline {X}$,
$U=\overline{U}\times_{\overline{X}}X$,
$\overline{U}_{0}=\overline{U}\otimes_{W}k$ and assume that
$R$ is flat, of finite type over $W$ and that
$R/pR\neq 0$. Let
$\hat{R}$ be the
$p$-adic completion of $R$ and $\tilde{\hat{R}}$ be the union of all
normalizations of $\hat{R}$ in
finite sub-Galois extensions of an algebraic closure of the field of
fractions of $\hat{R}$ such that
the normalization of $\hat{R}[1/p]$ in such finite extension is unramified
outside $D$ (cf.
\cite{Fa},II,i)). On
$\tilde{\overline{U}_{0}}=Spec(\tilde{\hat{R}}/p\tilde{\hat{R}})$, we have
a canonical log-structure wich
is the inverse image of that of $\overline{U}$. Let $\tilde{Z}$ be the
PD-envelope of the
log-immersion
$\tilde{\overline{U}_{0}}\to \tilde{\overline{U}_{0}}\times
\overline{U}_{0}$, locally, for the etale
topology,
${\cal O}_{\tilde{Z}}$ is a PD-polynomial algebra over ${\cal
O}_{\tilde{\overline{U}_{0}}}$.

\vskip 5mm

Let $({\cal F},{\cal F}^{i},\varphi^{i}_{\cal F},\nabla_{\cal F})$ be an
object of $MF^{\nabla,
[0,p-2]}_{k}({X_{0}})$. As a crystal, we can evaluate $\cal F$ on
$\tilde{\overline{U}_{0}}$ and obtain
an ${\cal O}_{\tilde{\overline{U}_{0}}}$-module ${\cal
F}_{\tilde{\overline{U}_{0}}}$  endowed with a
descending filtration  ${\cal F}^{i}_{\tilde{\overline{U}_{0}}}$ .

Let $J_{\tilde{Z}}$ be the PD-ideal of the log-immersion $\overline{U}\to
\tilde{Z}$, we have
naturally an object $({\cal
O}_{\tilde{Z}},J^{[i]}_{\tilde{Z}},
\varphi^{i}_{\tilde{Z}},\nabla_{\tilde{Z}})$
of
$MF^{\nabla,[0,p-2]}_{k}({X_{0}})$ and where $\varphi^{i}_{\tilde{Z}}$ is
by definition, the reduction
mod $p$ of
${\tilde\varphi}_{\tilde{Z}}/ p^{i}$ for some lifting (wich exist locally)
${\tilde\varphi}_{\tilde{Z}}$ over $W$ of $\varphi_{\tilde{Z}}$. On the complex
$${{\cal F}_{\tilde{\overline{U}_{0}}} \otimes_{{\cal
O}_{\tilde{\overline{U}_{0}}}}
({\cal O}_{\tilde{Z}}
\otimes_{{\cal O}_{\overline{U}_{0}}}
\Omega^{\bullet}_{\overline{U}_0/k}}(log\infty))$$
deduced from the integrable connection $\nabla_{\cal F}$, we have a
filtration $Fil^{r}$ whose degree
$q$-part is given by
$$Fil^{r}({\cal
F}_{\tilde{\overline{U}_{0}}} \otimes_{{\cal O}_{\tilde{\overline{U}_{0}}}}
{\cal O}_{\tilde{Z}}
\otimes_{{\cal O}_{\overline{U}_{0}}}
\Omega^{q}_{\overline{U}_0/k}(log\infty))=({{\cal
F}_{\tilde{\overline{U}_{0}}} \otimes_{{\cal
O}_{\tilde{\overline{U}_{0}}}} {\cal O}_{\tilde{Z}})^{r-q}\otimes_{{\cal
O}_{\overline{U}_0}}}
\Omega^q_{\overline{U}_0/k}(log\infty)$$
and where
$$({\cal F}_{\tilde{\overline{U}_{0}}} \otimes_{{\cal
O}_{\tilde{\overline{U}_0}}} {\cal O}_{\tilde{Z}})^{r}=\sum_{i+j=r}{\cal
F}^{i}\otimes
J^{[j]}_{\tilde{Z}}$$

Moreover, a refined crystalline Poincar\'e lemma tells us that the natural
morphism
$${\cal F}^{i}_{\tilde{\overline{U}_0}}\to
Fil^{i}({\cal F}_{\tilde{\overline{U}_{0}}} \otimes_{{\cal
O}_{\tilde{\overline{U}_0}}}
{\cal O}_{\tilde{Z}}
\otimes_{{\cal O}_{\overline{U}_{0}}}
\Omega^{\bullet}_{\overline{U}_0/k}(log\infty))$$
is a ${\rm Gal}(\tilde{\hat{R}}/{\hat{R}})$-equivariant resolution of the ${\cal
O}_{\tilde{\overline{U}_{0}}}$-module
${\cal F}^{i}_{\tilde{\overline{U}_{0}}}$.

Put ${\cal F}_{\tilde{Z}}={{\cal F}_{\tilde{\overline{U}_{0}}} \otimes_{{\cal
O}_{\tilde{\overline{U}_{0}}}}} {\cal O}_{\tilde{Z}}$, it is a filtered
${\cal O}_{\tilde{Z}}$-module as
noted above. We define the Frobenius
$\varphi^{r}_{{\cal F}_{\tilde{Z}}} : {\cal F}^{r}_{\tilde{Z}} \to {\cal
F}_{\tilde{Z}}$, $(r<p)$, as
the unique additive map whose restriction to the image of  ${\cal
F}^{i}\otimes J^{[j]}_{\tilde{Z}}$
is $\varphi^{i}_{\cal F}\otimes \varphi^{j}_{\tilde{Z}}$ for all $i, j$
such that $i+j=r$.
Finally, we define the ${\rm Gal}(\tilde{\hat{R}}/{\hat{R}})$-equivariant
morphism of complexes:
$$\varphi^{r} :  Fil^{r}({\cal F}_{\tilde{Z}}
\otimes_{{\cal O}_{\overline{U}_{0}}}
\Omega^{\bullet}_{\overline{U}_0/k}(log\infty))\to
{\cal F}_{\tilde{Z}}
\otimes_{{\cal O}_{\overline{U}_{0}}}
\Omega^{\bullet}_{\overline{U}_0/k}(log\infty)$$

as induced from the morphisms:

$$\varphi^{i-q}_{{\cal F}_{\tilde{Z}}}\otimes\wedge^{q}d\varphi/p :
{\cal F}^{i-q}_{\tilde{Z}}\otimes \Omega^{q}_{\overline{U}_0/k}
\to {\cal F}_{\tilde{Z}}\otimes \Omega^{q}_{\overline{U}_0/k}$$

Now we are ready to define the Frobenius
$$\varphi^{i} : {\cal F}^{i}_{\tilde{\overline{U}_{0}}} \to {\cal
F}_{\tilde{\overline{U}_{0}}}$$
 by using the resolution of ${\cal F}^{i}_{\tilde{\overline{U}_{0}}}$ given
above.

For $i<p$, we define the ${\rm Gal}(\tilde{\hat{R}}/{\hat{R}})$-module
${\bf V}_{U}({\cal F},i)$ as
the kernel of
$$1-\varphi^{i} :H^{0}(\tilde{\overline{U}_{0}},{\cal
F}^{i}_{\tilde{\overline{U}_{0}}})\to
H^{0}(\tilde{\overline{U}_{0}},{\cal F}_{\tilde{\overline{U}_{0}}})$$

Let $E=\tilde{\hat{R}}/p\tilde{\hat{R}}$, choose $\beta\in K^{ac}$ such
that $\beta^{p}=-p$ and for
$i<p$, define a filtration $Fil^{i}E=\beta^{i}E$ and Frobenius
$\varphi^{i}(\beta^{i}x)=x^{p}$, then
as explained in \cite{Tsu} proof of prop. 4.3.4 or \cite{Fa},II,f), we have an
isomorphism
$${\bf V}_{U}({\cal F},i)^*\simeq Hom_{MF^{\nabla,
[0,p-2]}_{k}({X_{0}})}({\cal F}[i],E)$$

Where ${\cal F}[i]$ is the twisted module defined by ${\cal F}[i]^{j}={\cal
F}^{i+j}$ and
$\varphi_{{\cal F}[i]}^{j}=\varphi_{\cal F}^{i+j}$. Moreover, using this
description, we deduce
that ${\bf V}_{U}({\cal F},i)$ is finite of order $p^{h}$
(\cite{Fa},Th.2.4) where
$h=\vert\Lambda\vert$ and $\Lambda$ is the index set in the definition  of
a uniform filtered module.

By \cite{Fa}, II,g) or \cite{Tsu}(4.4), if we regard ${\bf V}_{U}({\cal F},i)$
as a finite
locally \'etale constant sheaf on $U\otimes_{W}K$, we can glue the local data
${\bf V}_{U}({\cal
F},i)$, for various "small" $U$ (cf. \cite{Tsu} 3.3.2). There is a unique
finite locally constant sheaf
${\bf V}_{X}({\cal F},i)$ on $X\otimes_{W}K$ such that the restriction to
"small" $U$ is ${\bf
V}_{U}({\cal F},i)$. Finally, we define the covariant comparison functor ${\bf V}$ by ${\bf V}({\cal
F})={\bf V}_{X}({\cal F},p-2)(2-p)$,  and its contravariant version  ${\bf V}^\star$ by  
${\bf V}^*({\cal F})={\bf V}({\cal F})^*$.

\subsection {Association modulo $p$ for Siegel varieties}

Let us come back to the case of Siegel varieties. 
Let $X_{/\Z[1/N]}$ be the moduli scheme classifying p.p.a.v. with level $U$-structure over
$\Z[1/N]$. Its toroidal compactification over $\Z[1/N]$ is denoted by $\overline{X}$ 
(for some choice of a smooth $GL(\Z^g)$-admissible polyhedral cone decomposition of the convex cone of all positive
semi-definite symetric bilinear forms on $\R^g$). We have $S_U=X\otimes_{\Z[1/N]} \C$. 
Recall that, to the representation
$V_{\lambda\,/\F_p}$ of
$G_{\F_p}=G\otimes \F_p$ of highest weight $\lambda$, one can associate an etale sheaf
$V_\lambda(\F_p)$ resp. $V_\lambda(k)$  over
$X\otimes\Q$  resp. its extension of scalars to $k$.
One possible construction of this etale sheaf is by the theory of the fundamental group: any
representation of the arithmetic fundamental group $\pi_1(X\otimes\Q,\overline{x})$ on a
finite abelian group
$V$ gives rise to an etale sheaf whose fiber at $\overline{x}$ is $V_{\lambda\,/\F_p}$. Let us
consider the structural map
$f:A\rightarrow X\otimes\Q$ given by the universal principally polarized abelian surface with
level structure of type $U$ (we assume here $U$ sufficiently deep). The sheaf
$R^1f_*\Z/p\Z$ is \'etale. It
corresponds to an antirepresentation of the fundamental group taking values in $G(\Z/p\Z)$. Then,
composing with the representation $G_{\F_p}\rightarrow GL(V_{\lambda\,/\F_p})$, we obtain an \'etale
sheaf denoted by
$V_\lambda(\F_p)$. Similarly for $V_\lambda(k)$, by considering the extension of scalars from $\F_p$ to $k$:
$G_k\rightarrow GL_k(V_\lambda(k))$.

For any dominant weight
$\lambda$  of $G$, we have thus obtained a  
$V_{\lambda}(\F_p)$ of
${\cal R}ep_{\F_{p}}(X\otimes K)$. On the other hand, if moreover 
$\vert\lambda+\rho\vert<p-1$, the crystal $\overline{\cal V}^\vee_\lambda$
constructed in Section 5.2 satisfies the conditions of Definition 7 which turn it into an
object of
${\cal MF}^{\nabla,[0,p-2]}(\overline{X})$. To verify this, one starts with the
standard representation. On 
$$\overline{\cV}^\vee_1=R^1\bar{f}_*\Omega^\bullet_{\bar{A}/\overline{X}}(log\,
\infty_{\bar{A}/\overline{X}}),$$
 the Gauss-Manin connection
satisfies Griffiths transversality for the Hodge filtration, compatibility to Frobenius and
uniformity. The tricky point is to verify the strong divisibility condition (section 6.2, definition 7). 
 It follows from the
degeneracy of the Hodge spectral sequence which is proven in
\cite{Fa},Th.6.2. For general $\lambda$, we use that
that
$\overline{\cal V}^\vee_\lambda$ is a
sub-object (and quotient) of a first direct image for some Kuga-Sato variety and the fact
that
${\cal MF}^{\nabla,[0,p-2]}(\overline{X})$ is a an abelian category. Note that the objects $\overline{\cal V}_\lambda$
(without dualizing) do not belong to this category, as their weights don't fit the bound.

\begin{thm} ( \cite{CF}Th.6.2(iii))
 $$\bV^*(\overline{\cal
V}^\vee_\lambda)=V_{\lambda}^\vee(\F_p)$$
that is, $V_{\lambda}(\F_p)$ and $\overline{\cal
V}_\lambda^\vee$ are associated.
\end{thm}

The proof is given in  \cite{CF}Th.6.2(iii). In fact, there, the
result is proven only in the $\Q_p$-coefficients case, but for
$\vert\lambda+\rho\vert<p-1$ the
proof is valid word for word in the integral context. The key argument is the existence of the
minimal compactification whose boundary has relative codimension $\geq  2$. The next section
gives more details about this.

\subsection{The Comparison Theorem} 

We will explain the relative comparison theorem Th.6.2 of Faltings \cite{Fa} in our particular setting.
 In fact we merely extend the
arguments sketched in \cite{CF}, p.241. Before going into our situation, we recall the method of \cite{Fa} (we
hope that more details will be given by the experts in the future).
 
\subsubsection{General setting}
For any $p$-adic ring $ R$ (smooth integral over $\bZ_p$), we denote by $\hat R$ its $p$-adic
completion, by $R_0=R\otimes_{\bZ_p}\bZ/p\bZ$ its reduction mod.$p$; let $F$ be the field of fractions of $\hat R$;
choose an algebraic closure $\overline{F}$ of $F$ and form 
$\overline{\hat R}$, union of of all the normalizations of $\hat{R}$ in
finite sub-Galois extensions of $\overline F$.
Put
$S=\overline{\hat R}/p\overline{\hat R}$.

Let $f : Y\rightarrow Spec(R)$ be a smooth and proper morphism of schemes of relative dimension $d<p-1$,
$Y_0=Y\otimes_{\bZ_p}\bZ/p\bZ$ the special fiber,  ${\overline Y}=Y\otimes_{R}\overline{\hat R}$, 
${\overline Y}_\eta=Y\otimes_{R}\overline{F}$, 
${\overline Y}_0=Y_0\otimes_{R_0}S$ and
$f_0 : Y_0 \rightarrow Spec(R_0)$, $\overline{f} : \overline{Y} \rightarrow Spec(\overline{\hat R})$,
 $\overline{f}_\eta : \overline{Y}_\eta \rightarrow Spec(\overline{ F})$, 
$\overline{f}_0 : \overline{Y}_0 \rightarrow Spec(S)$ the
corresponding morphisms. We have the following standard diagram :

$$\begin{array}{ccccc}
 \overline{Y}_\eta&\buildrel{\overline j}\over\hookrightarrow & \overline{Y}&\buildrel{\overline i}\over\hookleftarrow &{\overline Y}_0\cr 
\downarrow \overline{f}_\eta& &\downarrow \overline{f}&&\downarrow \overline{f}_0\cr
Spec({\overline F})&\hookrightarrow &Spec(\overline{\hat R})&\hookleftarrow &Spec(S)\cr
\end{array}$$

Let $R\Psi(S(1))={\overline i}^*R{\overline j}_*(S(1))$ be 
the "relative complex of $p$-adic vanishing cycles" for the constant sheaf $S(1)=\Z/p\Z(1)\otimes S$. 
This object is not explicitely introduced in \cite{Fa}, but as explained in \cite{Ill2}, we can rewrite the complex computing
\'etale cohomology as a complex of vanishing cycles. Then we have a "Kummer" map : 
$$R\Psi(S(1))\rightarrow \Omega^\bullet_{\overline{Y}_0/Spec(S)}.$$
 
Taking direct images, we obtain natural maps :

$$R^*f_{0*}(\Omega^\bullet_{Y_0/Spec(R_0)})\otimes_R S \rightarrow
R^*\overline{f}_{0*}(\Omega^\bullet_{\overline{Y}_0/Spec(S)})
\leftarrow R^*\overline{f}_{0*,et}(R\Psi) \simeq R^*\overline{f}_{\eta *,et}(S)$$

$$R^*\overline{f}_{0*,et}(R\Psi) \simeq R^*\overline{f}_{\eta *,et}(S)\leftarrow
R^*\overline{f}_{\eta *,et}(\bZ/p\bZ(1))\otimes_R S$$ 

Faltings (\cite{Fa}, page 72, see also recent corrections of the corresponding proof in informal notes by the author) shows that 
the second arrow is an "almost-isomorphism"; wich implies that
the sheaves  $R^*f_{0*}(\Omega^\bullet_{Y_0/Spec(R_0)})$ and $R^*\overline{f}_{\eta *,et}(S(1))$
are associated.

\subsubsection{Setting for Siegel varieties} The notations are those of section 6.3. Let
$U=Spec(R)\subset X$ be an affine open subset and $f :Y_U\to U$ be the restriction of $f_s
:Y=A\times_X\ldots\times_XA \rightarrow X$, where $A$ is the universal abelian variety,
 we assume $s<p-1$. Let $\hat{X}$  be
the formal completion of $X$ along the special fiber. Let 
${\hat f}  :{\hat Y}_U\to {\hat U}$ be
the base change of $f$ to the affine  formal scheme 
${\hat U}=Spf({\hat R})$. Over
$Spec(\hat{R}\otimes \bQ_p)$, we have two Žtales sheaves
$R^s{\hat f}_*\bZ/p\bZ(1)$ and
$\bV^*(R^s{\hat f}_*(\Omega^\bullet_{Y_U\otimes\F_p/U\otimes\F_p}))$. As explained in the
general setting subsection, there is a functorial isomorphism of Žtales sheaves:
$$R^s\hat{f}_*\bZ/p\bZ(1) \simeq  \bV^*(R^s \hat{f}_*(\Omega^\bullet_{Y_U\otimes\F_p/U\otimes\F_p}))$$ 
over $\hat{U}$.
By
functoriality, these local isomorphisms glue to a global one over $\hat X$.

Let  $X^*$ be the minimal compactification of $X$ over $\Z_p$. It is defined in \cite{CF} Th.2.5
Chapter V. It is projective, normal of finite type; its boundary admits a natural stratification whose
strata have codimension at least $2$ (since we assume $g\geq 2$). We apply Grothendieck's GAGA theorem
to deduce  that the isomorphism over $\hat{X}$ between the sheaves  $R^s{\hat f}_*\bZ/p\bZ(1)$ and
$\bV^*(R^s{\hat f}_*(\Omega^\bullet_{Y_U\otimes\F_p/U\otimes\F_p}))$ is algebraic.
 More
precisely, every Žtale covering of the formal scheme $\hat{X}$ is defined by an Žtale finite ${\cal
O}_{\hat X}$-algebra $\cal A$. Since the minimal compactification is normal and has boundary of
codimension
$\geq 2$, this algebra extends to $\hat {X^*}$ (\cite{EGA4}, Cor 5.11.4) and so define an algebraic 
Žtale covering of
$X$ whose base change to $\hat X$ is $\cal A$, we deduce an equivalence of topoi $X_{et}\simeq {\hat
X}_{et}$.  As the morphism
$f$ is  proper and smooth, the sheaf 
$R^s{\hat f}_*\bZ/p\bZ(1)$ on $\hat X$ is locally constant and so descend to $X$ and gives the sheaf
$R^s f_*\bZ/p\bZ(1)$. By construction, the sheaf 
$\bV^*(R^s{\hat f}_*(\Omega^\bullet_{Y_U\otimes\F_p/U\otimes\F_p}))$ is also locally constant and
also descend to $X$ and gives the sheaf
$\bV^*(R^sf_*(\Omega^\bullet_{Y_U\otimes\F_p/U\otimes\F_p}))$. 

Moreover, as $X_{et}\simeq {\hat X}_{et}$,  
every formal morphism between $R^s{\hat f}_*\bZ/p\bZ(1)$ and  
$\bV^*(R^s{\hat f}_*(\Omega^\bullet_{Y_U\otimes\F_p/U\otimes\F_p}))$ is algebraic. This shows that $R^s
f_*\bZ/p\bZ(1)$ is associated to $R^sf_*(\Omega^\bullet_{Y\otimes\F_p/X\otimes\F_p})$ for the
asociation without divisor at infinity and $R^s
f_*\bZ/p\bZ(1)$ is associated to $R^s{\overline f}_*(\Omega^\bullet_{{\overline
Y}\otimes\F_p/{\overline X}\otimes\F_p}(log\infty))$ for the association with divisor at infinity.

\vfill\eject

\section{Proof of Theorem 1}

\subsection{A lemma on modular representations}

Our reference for results used in this Section are \cite{Bo} VIII.13.2 and
\cite{J}, II.3. Let
$\hat{T}$ be the standard maximal torus in
$\hat{G}$. One has
$$\hat{T}=\{(t_1,\ldots,t_g,u;x); u^2=t_1\ldots t_g\}$$ The degree $2$
covering $\hat{G}\rightarrow
GO_{2g+1}$ induces on $\hat{T}$ the projection
$$(t_1,\ldots,t_g,u;x)\mapsto
diag(t_1,\ldots,t_g,xt_g^{-1},\ldots,xt_1^{-1},x)$$ We view the Weyl
group $W_{\hat{G}}$ as a subgroup of $\hat{G}_{/\Z}$ by using permutation
matrices in a standard way.
Let $W'$ be the subgroup of
$W_{\hat{G}}$ consisting in the permutations $w_B$ ($B\subset [1,g]$) acting by
$t^{w_B}=t'$ where $t=(t_1,\ldots,t_g,u;x)$ and
$t'=(t'_1,\ldots,t'_g,u';x)$ with $t'_i=t_i^{-1}$ if
$i\in B$, $t'_i=t_i$ if $i\notin B$, and $u'=u\cdot t_B^{-1}$ where
$t_B=\prod_{i\in B} t_i$.

Let $\hat{B}=\hat{T}.\hat{N}$ be the Levi decomposition of the standard
Borel subgroup $\hat{B}$.
Recall we assumed
${\bf GO(\omega)}$ for
$\overline{\rho}_\pi$. We can assume that
$\overline{\rho}_\pi(D_p)\subset \hat{B}(k)$.
 Throughout this section, we assume that

{\bf (RLI)} : there exists a split (non necessarily connected) reductive
Chevalley subgroup
$H$ of
$\hat{G}_{/\Z}$  with
$W'\propto \hat{T}\subset H$, and a subfield $k'\subset k$, of order say
$\vert k'\vert=q'=p^{f'}$
($f'\geq  1$), so that $H(k')_\nu\subset Im\,\overline{\rho}_\pi$  and
$\overline{\rho}_\pi(I_p)\subset H^0(k')$.
Where $H(k')_\nu$ is the subgroup of $H(k')$ consisting in elements whose
$\nu$ belongs to
$Im\,\nu\circ \bar{\rho}_\pi$.

{\bf Comment:} It has been pointed to us by R. Pink that if $H$ is
connected and $W'\propto
\hat{T}\subset H$, then $H$ should contain the derived group of $\hat G$;
then, {\bf (RLI)} becomes
in some sense an assumption of genericity for
$\pi$ and
$p$, but it cannot be verified in a single example for $g\geq  2$, hence
our insistance on
the possible disconnectedness of
$H$: it allows us to show the existence of concrete examples for the theorem.

Let $H^0$ be the neutral component of $H$ over $\Z$. Its semisimple rank is
$g$. Recall that in the
condition of  Galois ordinarity {\bf (GO)}, we introduced an element
$\hat{g}\in \hat{G}$ so that
 $$\rho_\pi(D_p)\subset \hat{g}\cdot
\hat{B}(\cO)\cdot \hat{g}^{-1}$$

 Recall the convention (in vigor since Sect.3.3) that we omit the
conjugation by $\hat{g}$, thus writing
$\hat{B}$, $\hat{N}$, $\hat{T}$ instead of $\hat{g}\cdot \hat{B}\cdot
\hat{g}^{-1}$ and so on.

The subdata
$(H^0,\hat{T},\hat{B}\cap H^0)$ in $(\hat{G},\hat{T},\hat{B})$ induce an
inclusion of the set of roots
of $H^0$ into that of $\hat{G}$: $\Phi_{H^0}^\pm\subset \Phi^\pm$. Let
$\Phi'=\Phi\cap
Vect_\Q(\Phi_{H^0})$ and $\Delta'$ a system of basis made of positive
simple roots for $\Phi'$. By
\cite{Bo} VI, $n^o\,1.7$, Prop.24, it can be completed into a basis
$\Delta$ of $\Phi$ contained in
$\Phi^+$. Note that $\Phi_{H^0}$ is a subsystem of maximal rank in $\Phi'$. Let
$\Delta_{H^0}$ be the basis of $\Phi_{H^0}$ contained in $\Phi_{H^0}^+$. A
priori, it could be
different from $\Delta'$ (not in the examples we have in view though). Let
$\hat{\varpi}_i$ be the fundamental weights in $X$ of $\hat{G}$. Let

$$\Phi_{H^0}^\perp=\{\lambda\in X; \langle
\lambda,\beta^\vee\rangle=0\, {\rm for}\,\beta\in \Phi_{H^0}\}$$

where $\alpha^\vee$ denotes the coroot corresponding to a root $\beta$. We
write $\hat{\varpi}=\hat{\varpi}_g$ for
the minuscule weight of $\hat{G}$; it is the highest weight of the spin
representation $V_{/\F_p}$
of $\hat{G}$.

Observe that $\Phi_{H^0}^\perp=\Z\cdot \nu$. Let $X'$
be the
$\Z$-module generated by
$\Delta'$. It is equal to
$\Z\hat{\varpi}_1\oplus\ldots\oplus\hat{\varpi}_{g}\subset X$, and one has
$$X=X'\oplus \Phi_{H^0}^\perp$$

The irreducible representations of $H^0$ over $k'$ (or over any perfect
extension of
$\F_p$) are classified by $X'{}^+\times \Phi_{H^0}^\perp$. We shall
consider certain
(absolutely) irreducible representations over $k'$ of the abstract group
$H^0(k')$.

Let $e={q'-1\over(\w,q'-1)}$. Recall that
by the formula  $\nu\circ \rho_\pi=\chi^-\w\cdot \omega_\pi$, the kernel
$N$ of the homomorphism
$$X\rightarrow Hom(\hat{T}(k')_\nu,k'{}^\times)$$
contains $$(q'-1)\cdot X'\oplus e\cdot \Z\cdot \nu$$
A fortiori, we can assume equality holds:
$$Im\,\bar{\rho}_\pi=k'^{\times\, \w},\,{\rm hence,}\, N=(q'-1)\cdot
X'\oplus e\cdot \Z\cdot \nu.$$

 It results easily from Steinberg's theorem (see
Chapter II, Prop.3.15 and Coroll.3.17 of
\cite{J}) that the irreducible representations of the abstract group
$H^0(k')_\nu$ are classified by

$$X_{H,q'}=\{(v,a)\in X'{}^+\times [0,e-1]; 0\leq \langle
v,\beta^\vee\rangle\leq q'-1\,\mbox{\rm for all}\,\beta\in \Delta_{H^0}\}$$

For brevity, we call such weights $q'$-reduced, although the terminology is
not conformal to that of
Jantzen's book Chapter II, Section 3.

For $\mu\in X_{H,q'}$, we write $W(\mu)$ for the corresponding
$H^0$-representation and
$\Pi_{H^0}(\mu)\subset X$
 for its set of weights, resp. $\overline{\Pi}_{H^0}(\mu)\subset
Hom(\hat{T}(k'),k'{}^\times)$ the set
of their restrictions to $\hat{T}(k')_\nu$. Let
$\Pi_{\hat{G}}(\hat{\varpi})$ resp.
$\overline{\Pi}_{\hat{G}}(\hat{\varpi})$ the set of weights (resp. of the
functions on
$\hat{T}(k')$ that they induce) associated to the spin representation
$\bV_{/k'}$ of $\hat{G}$.

Recall that
${\Pi}_{\hat{G}}(\hat{\varpi})=\{\hat{\varpi}^{w'}; w'\in W'\}$ and that we
assumed
$W'\propto
\hat{T}\subset H$.

\begin{lem} For $p>5$, if $W(\mu)$ is a simple $H^0_{k'}$-module with
highest weight $\mu\in X_{H,q'}$
with
$\overline{\hat{\varpi}}=\overline{\mu}$ and
$\overline{\Pi}_{H^0}(\mu)\subset \overline{\Pi}_{\hat{G}}(\hat{\varpi})$,
then $\mu=\hat{\varpi}$.
\end{lem}

{\bf Remark:} For $p=5$, $\hat{G}=Spin(5)$ and $H\subset \hat{G}$,
isomorphic to $ SL(2)\times SL(2)$
via $\hat{G}\cong Sp(4)$,
$\mu=3\hat{\varpi}_2$, the lemma is false, hence the necessity of the
assumption $p>5$.

{\bf Proof:} Since $\overline{\mu}=\overline{\hat{\varpi}}$, one has
$\mu-\hat{\varpi}\in (q'-1)X$.

\bigskip

1) Let us first check that $\mu-\hat{\varpi}\in N\cap
\Phi_{H^0}^\perp=e\cdot \Phi_{H^0}^\perp$.

Let $\alpha\in\Delta_{H^0}$. We want $\langle
\mu-\hat{\varpi},\alpha^\vee\rangle=0$.
We start with a preliminary observation:

For any simple root $\alpha\in\Delta_{H^0}$, $\langle
\hat{\varpi},\alpha^\vee\rangle\in\{-1,0,1\}$. Indeed, this is true for any
fundamental weight $\hat{\varpi}$ of
$\hat{G}$. In particular for our minuscule weight
$\hat{\varpi}$.

Then, we distinguish three cases

\begin{itemize}

\item If $\langle
\hat{\varpi},\alpha^\vee\rangle=1$, we have
$\langle\mu,\alpha^\vee\rangle=1$ because $\mu$ is
$q'$-reduced.

\item If $\langle
\hat{\varpi},\alpha^\vee\rangle=0$; let us exclude the possibility $\langle
\mu,\alpha^\vee\rangle=q'-1$. Since $q'-1\geq  1$ we would have
$\mu-\alpha\in\Phi_{H^0}(\mu)$ as
the $\alpha$-string of $\mu$ has length $q'-1$. Hence by the assumption, we
could write
$\mu-\alpha= \hat{\varpi}^y+(q'-1)\lambda$ for some $y\in W'$ and
$\lambda\in X$.

But $\langle
\hat{\varpi}^y,\alpha^\vee\rangle\in\{-1,0,1\}$, and $\langle
\mu-\alpha,\alpha^\vee\rangle=q'-3$ hence $q'-1$ should divide $1,2$ or $3$
impossible
since
$q'-1>3$.

\item
Similarly, if $\langle
\hat{\varpi},\alpha^\vee\rangle=-1$, we must exclude $\langle
\mu,\alpha^\vee\rangle=q'-2$. Again $\mu-\alpha\in \Pi_{H^0}(\mu)$, hence
$\mu-\alpha\equiv
\hat{\varpi}^y\, mod. (q'-1)X$. But $\langle
\hat{\varpi}^y,\alpha^\vee\rangle\in\{-1,0,1\}$ and $\langle
\mu-\alpha,\alpha^\vee\rangle\equiv -3\,{\rm mod}\,(q'-1)$, hence $(q'-1)$
should divide $2,3$ or $4$;
impossible since
$q'-1>4$.

\end{itemize}

2) Thus, $\mu-\hat{\varpi}\in \Phi_{H^0}^\perp=e\cdot\Z\cdot\nu$ as
desired (actually it shows that
$\langle
\hat{\varpi},\alpha^\vee\rangle\geq  0$ for any $\alpha\in \Delta_{H^0}$).
 Since the components of
$\hat{\varpi}$ and
$\mu$ along $\nu$ are between $0$ and $e-1$, and
$\mu-\hat{\varpi} \, (mod.\,e)$ we conclude $\mu=\hat{\varpi}$. The lemma is
proven.

It is the main ingredient in the proof of the following result.

\begin{lem} Let $\sigma:\Gamma\rightarrow GL_k(W)$ be a continuous Galois
representation such that
 for
any
$g\in
\Gamma$, the characteristic polynomial of $\overline{\rho}_\pi(g)$
annihilates $\sigma(g)$. Assume
 that
$p-1>max(4,j_A)$, that
$\overline{\rho}_\pi$ satisfies {\bf GO($\omega$)} and {\bf (RLI)},

then, either $W=0$, or the two characters $1$ and $\omega^{-\w}$ restricted
to $I_p$
occur as subquotients of
$W$ viewed as an
$I_p$-module.
\end{lem}

{\bf Comment:} One could naturally ask whether the simpler assumptions that
$\overline{\rho}_\pi$ is
absolutely irreducible and for any $g\in \Gamma$ the characteristic
polynomial of
$\overline{\rho}_\pi(g)$ annihilates
$\sigma(g)$ are sufficient to conclude that all constitutents of $\sigma$ are
copies of $\overline{\rho}_\pi$. This statement is true for $g=1$, but, it
is false for $g=2$. A
counterexample has been found by J.-P. Serre. He lets $\Gamma$ act on
$\F_p^4$ through the so-called
cuspidal representation of the non-split central extension ${}_2A_5$ of the
icosaedral group $A_5$. It
is four-dimensional, symplectic and absolutely irreducible. Then,
$(W,\sigma)$ is one of the two irreducible $2$-dimensional of this group.
This is why we introduced
{\bf (RLI)}. This assumption is not satisfied in the example there. Also,
thanks to the ordinarity
assumption {\bf (GO)}, we focused our attention on the highest weight of
$\overline{\rho}_\pi$ (which is a local information at $p$) rather than the
global representation
$\overline{\rho}_\pi$ itself.

\bigskip

{\bf Proof:} Assume $W\neq 0$; let $\Gamma'$ be the inverse image by
$\overline{\rho}_\pi$ of $H(k')$
in
$\Gamma$ and
$\Gamma''$ the kernel of $\overline{\rho}_\pi$ restricted to $\Gamma'$.
Then $\sigma(\Gamma'')$ is a
nilpotent $p$-group. Thus, replacing $W$ by its submodule fixed by
$\sigma(\Gamma'')$, still
denoted by $W$, one can assume that
$W$ is a non-zero module on which
$\Gamma'$ acts through $H(k')_\nu$:

$$\begin{array}{lcc}
\Gamma'&\rightarrow&GL_k(W)\\\downarrow\overline{\rho}_\pi&\nearrow&\\H(k')&
&\end{array}$$

We first treat the case of $\omega^{-\w}$. Let $H^0$ be the neutral
component of $H$. Let
$\tilde{W}=Ind_{H^0(k')_\nu}^{H^0(k')} W$. It is an $H^0(k')$-module, and
for any
$t\in\hat{T}(k')_{\nu}$, the action of $t$ on $\tilde{W}$  is annihilated
by $\prod_{w\in
W'}(X-\varpi^{w}(t))$.
By
 Steinberg theorem (\cite{J} Sect
II.3.15), the space
$W$ viewed as $H^0(k')$-module has a subquotient $W(\mu)$ which comes from
an algebraic
 simple
$H^0_{k'}$-module corresponding to a
$q'$-reduced highest weight $\mu$. We associate to this representation the sets
${\Pi}_\mu$ resp. $\overline{\Pi}_\mu$ as above. By the assumption
$W'\subset H$, one can assume that
$\overline{\Pi}_{H^0}(\mu) \subset\overline{\Pi}_{\hat{G}}({\hat{\varpi}})$ and
$\overline{\hat{\varpi}}= \overline{\mu}$ (if $
\overline{\mu}=\overline{\hat{\varpi}}^{w'}$ for some $w'\in W'$,
simply replace $W(\mu)$ by $W(\mu^{w'{}^{-1}})$ which also occurs as
$H^0_{k'}$-subquotient of $W$).
By the previous lemma, for
$p>5$, we have ${\hat{\varpi}}=\mu$. Let $x$ be a highest weight vector in
$W(\mu)$ for
$H^0_{\F_p}$. It is fixed by
$H\cap \hat{N}(k)$. Since $I_p\subset \overline{\rho}_\pi^{-1}(H^0(k))$,
the action of $I_p$ on $x$
is through its image by $\hat{\varpi}_g\circ(\overline{\rho}_\pi\,
mod.\hat{N})$. By the assumption {\bf (GO)}, and
Lemma 3, this character is equal to
$\omega^{-\w}$ on $I_p$ which therefore occurs as a subquotient of
$W\vert_{I_p}$. To treat the case of the trivial character, we consider
instead of the highest weight $\mu$ by the lowest
weight
$\mu'$ of $W(\mu)$; we can assume that
$\overline{\mu}'=\overline{\hat{\varpi}}^{w_0}$ where $w_0$ is the longest
element of
$W_{\hat{G}}$. Let $N_{H^0}$ be the unipotent radical of a Borel of $H^0$
adapted to $(GO)$. On the lowest weight quotient
$W(\mu)/N_{H^0}\cdot W(\mu)$, $\overline{\rho}_\pi$ acts by
$\hat{\varpi}^{w_0}\circ(\overline{\rho}_\pi\, mod.\hat{N})$,
which is trivial by (3.3.2).QED

\subsection{Deducing Theorem 1 from Theorem 6}

Recall we have fixed $\lambda=(a_g,\ldots,a_1;c)$ with $c=a_g+\ldots +a_1$ and $\vert\lambda+\rho\vert<p-1$. We have the following
reduction steps:
 
\noindent 1) By Poincar\'e duality, and self-duality of the Hecke operators for $\ell$ prime to $N$, Statement (i) of
Theorem 1 is equivalent to the vanishing of
$$H_*^j(S_U,V_\lambda(k))_\m=0\quad for \,q<d$$
where $\star=c,\emptyset$. These modules are artinian over $\cH_\m$, so by Nakayama's lemma, it is enough to show that their
$\m$-torsion vanishes:

$$(7.2.1)\quad H_*^j(S_U,V_\lambda(k))[\m]=0\quad for\,*=\emptyset\, or \, c\, {\rm and}\, q<d$$ 
which we will prove below.

\noindent 2) Then, Statements $(ii)$ and $(iii)$ are easy consequences of $(i)$ as can be seen by induction on $q<d$ using the
long exact sequences
$$0\rightarrow V_\lambda(\cO)\rightarrow V_\lambda(\cO)\rightarrow
V_\lambda(\cO/{\varpi}\cO)\rightarrow 0$$
and
$$0\rightarrow V_\lambda({\varpi}^{-1}\cO/\cO)\rightarrow V_\lambda(K/\cO)\rightarrow
V_\lambda(K/\cO)\rightarrow 0$$

For instance, from the latter, one obtains, with obvious notations: if
$H^{q-1}_*(K/\cO)_\m=0$, then $H^q_*({\varpi}^{-1}\cO/\cO)_\m\rightarrow H^q(K/\cO)_\m[{\varpi}]$ is an
isomorphism; hence by Nayamama's lemma, assertion one implies that
$H^q_*(K/\cO)_\m$ vanishes for $q<d$. 

Note that since
$p>j_A>a_g\ldots\geq  a_1\geq  0$, one knows that $V_{\lambda\,\F_p}$ is absolutely irreducible (see
for instance Proposition II.3.15, p.222, of \cite{J}).

3) As in section 6.3,  $X_{/\Z[1/N]}$ is the moduli scheme classifying p.p.a.v. with level $N$ structure over
$\Z[1/N]$. Its toroidal compactification over $\Z[1/N]$ is denoted by $\overline{X}$. Let
$V_\lambda(\F_p)$ resp. $V_\lambda(k)$  be the \'etale sheaf over
$X\otimes\Q$  in $\F_p$- resp. $k$-vector space corresponding to $V_{\lambda\,\F_p}$. Using the
etale-Betti comparison
isomorphism (and its equivariance for algebraic correspondences), Theorem 1
will be proven if we show the vanishing of the etale cohomology groups corresponding to $(7.2.1)$. 

\vskip 5mm

This interpretation as Žtale cohomology allows us to view $H^j_*(S_U,V_\lambda(\F_p))$ as a
$\F_p[Gal(\overline{\Q}/\Q)]\times {\cal H}_K$-module:
$$H^j_*(X,V_\lambda(\F_p))\cong H^j_{et,*}(X\otimes
\overline{\Q},V_\lambda(\F_p)).$$ 

{\bf Remark:} The $\F_p$-coefficients are useful to apply Fontaine-Laffaille and Faltings theory,
while the $k$-coefficients will come in when we localize at the maximal ideal $\m$ of ${\cal
H}_K(\cO)$.

Let $\overline{\cV}^\vee_\lambda$ be the object of ${\cal MF}^{\nabla,[0,p-2]}(\overline{X})$ associated to
$\lambda$ as in Section 5.2. Recall that $\overline{\cV}^\vee_\lambda$ has a filtration of length $\vert\lambda\vert$;
since
 $d+\vert\lambda\vert <p-1$
and since $\overline{\cV}^\vee_\lambda$ and $V_\lambda(\F_p)^\vee$ are associated (Theorem 8 above, section 6.3), 
we can apply Th.5.3 of
\cite{Fa} (see theorem 7, section 6.1), so that for any $j\geq 0$ : $H^j_{et,*}(X\otimes\overline{\Q}_p,V_\lambda(\F_p))^\vee$
is the image by the Fontaine functor $\bV^*$ of
$H^j_{cris,*}(S_U\otimes k,\cV_\lambda^\vee)$. 

Note that since we work mod.$p$ instead of
mod.$p^n$, we have
$$H^j_{log-crys,*}(X\otimes k,\cV_\lambda^\vee)=H^j_{log-dR,*}(X\otimes \F_p,\cV_\lambda^\vee)$$

We have constructed in Section 5.3.2 a filtered complex of coherent sheaves
$\overline{\cal K}_\lambda^\bullet$ on $\overline{X}\otimes \F_p$ by functoriality from the BGG
resolution of the
$G_{\F_p}$-module $V_{\lambda_{\F_p}}$. It follows from Theorem 6 that there are isomorphisms of
filtered
$\F_p$-vector spaces:
$$ H^j_{log-dR}(X\otimes{\F_p},\cV_\lambda^\vee)\cong H^j(\overline{X}\otimes
\F_p,\overline{\cal K}_\lambda^\bullet)$$ 
and
$$ H^j_{log-dR,c}(X\otimes \F_p,\cV_\lambda)\cong H^j(\overline{X}\otimes
\F_p,\overline{\cal K}_\lambda^{\bullet}{}^{sub})$$
  
where $\overline{\cal
K}^\bullet_\lambda$ resp. $\overline{\cal K}^\bullet_\lambda{}^{sub}$ denotes the canonical, resp.
subcanonical Mumford extension of the filtered complex of sheaves
${\cal K}_\lambda^\bullet$. The resulting filtration on the right-hand side is called
the
$F$-filtration; it corresponds via these isomorphisms to the Hodge filtration on the left-hand side.
The weights of this filtration can be computed as in 
\cite{RT2} (who treats the case $g=2$): Let us consider the map
$$W_G\rightarrow \Z,\, w\mapsto p(w)=-(w(\lambda+\rho)(H)-\rho(H))$$
where $H=diag(0,\ldots,0,-1,\ldots,-1)$.  Let $W_M$ be the Weyl group of the Levi
subgroup
$M$ of the Siegel parabolic. Observe that this map factors through
the quotient $W_M\backslash W_G$; this
quotient is in bijection with the set
$W^M$ (cf.p.229 of
\cite{CF}). By Theorem 6, sect. 5.4, we have
$$gr^pH^j_{log-dR,*}=\bigoplus_{w\in W^M,p(w)=p,\ell(w)\leq j-p} H^{j-\ell(w)}(\overline{X}\otimes
\F_p,\overline{\cal W}^\vee_{w(\lambda+\rho)-\rho})$$
Note that, unfortunatly, $p$ is not a good notation for the degre of our Hodge filtration.
 The image $p(W_G)$ of $p$ is therefore the set of possible weights occuring in
$H^j_{crys,*}$ for $j\leq d$. Moreover, $p$ is injective on
$W_M\backslash W_G$, and its values are exactly the
$j_B$ ($B\subset A$).  The set of possible lengthes $\ell(w),\,w\in
W^M$ is $[0,d]$. For each
$j<d$, let us consider the set $W^M(j)=\{w\in W^M; \ell(w)\leq j\}$; the key observation is that for
$j<d$, $W^M(j)$  
 does not contain the unique element $w\in W^M$ such
that
$\ell(w)=d$, namely the one acting by $(a_g,\ldots,a_1;c)\mapsto (-a_g,\ldots,-a_1;c)$. But this
element is the unique one for which $p(w)$ takes on its maximal value:
$j_A$.  Hence, this maximal weight
does not occur in
$H^j_{log-dR,*}(X\otimes
\F_p,\cV_\lambda^\vee)$ for $j<d$.

On the other hand, under assumptions {\bf (Gal)} and {\bf (GO)},
$\overline{\rho}_\pi$ is ordinary with weights given by
$j_B$ for all subsets $B\subset A$; in particular
$j_A$ and
$0$ indeed occur with multiplicity one; actually, even if we replaced {\bf (GO)} 
by geometric ordinarity, it would result from 
 lemma 3, Sect.3.3, that $0$ and $j_A$ do occur in $\rho_\pi$). Now, consider the global Galois
representation
$\sigma^j$ on $W_j= H^j_*(X\otimes\overline{\Q},V_\lambda(k))[\m]$, the kernel of $\m$ in the
module
$H^j_*(X\otimes\overline{\Q},V_\lambda(k))$. The Eichler-Shimura relations
imply for any $g\in Gal(\overline{\Q}/\Q)$, the characteristic polynomial of
$\overline{\rho}_\pi(g)$ annihilates
$\sigma^j(g)$. Our lemma 13 sect. 7.1, shows, assuming {\bf (RLI)}, that this implies that $W_j$ admits
$\overline{\rho}_\pi$ as subquotient. This is a contradiction since the maximal weight
$j_A$ occurs in
$\overline{\rho}_\pi$ but not in $W_j$.

\subsection{Examples}

 Let $F$ be a real quadratic field with Galois group $\{1,\sigma\}$. Let $\Gamma_F={\rm Gal}(\overline{\Q}/F)$. Let $f$ be a
holomorphic Hilbert cusp form for $GL(2)_{/F}$ of weight
$(k_1,k_\sigma)$,
$k_1,k_\sigma\geq  2$, $k_1= k_\sigma+2m$ for an integer $m\geq 1$. Assume it is a new form of conductor $\n$ which is eigen
for Hecke operators $T_v$ ($v$ prime to $\n$); denote by $a_v$ the corresponding eigenvalues. Since the weight of
$f$ is not parallel,
$f$ does not come from $\Q$. Let
$f_\sigma$ be the inner conjugate of $f$ by $\sigma$. Let $\epsilon$ be the finite order part of its central character. We
assume that $\epsilon$ factors through the norm map. Starting from
\cite{Y}, a series of works have established that
$f$ admits a holomorphic theta lift
$\pi$ to
$G(\bA)$ where
$G=GSp(4)$ (see \cite{R1} and \cite{R2}). Since $f$ does not come from $\Q$,
$\pi$ is cuspidal; moreover, by Th.6.2 of \cite{Vi}, $f$ admits also a
globally generic lift, hence $\pi$ is stable at every place. It occurs in the
$H^3$ of a Siegel variety of some level, say $N$, with coefficient system of highest weight
$\lambda=(a,b;c)$ where $a=k_\sigma+m-3$,
$b=m-1$, and $c=a+b$. At the moment, the level $N$ of $\pi$ can only be said to be multiple of $N(\n)D_F$  where
$D_F$ is the discriminant of
$F$; this should be the
exact level of
$\pi$, but this can not yet be established in general.

Let
$\Q(f)=\Q[a_v]_v$ be the number field generated by the eigenvalues of $f$; one can take
$\Q(f)$ as field of definition of
$\pi$ (although this may not be the smallest possible one, as pointed out by Prof. Yoshida). 
For any prime $\p_f$ of $\Q(f)$ prime to $\n$, the $\p_f$-adic Galois representation associated to $\pi$ exists; it is given by
$$(7.3.1)\quad \rho_\pi=Ind^F_\Q
\rho_f$$
it is absolutely irreducible. The
conductor of
$\rho_\pi$ is 
$Norm(\n)\cdot D_F$; this results from the fact that $\n$ is also the (prime-to-$p$ part of the) conductor of
$\rho_f$ by Carayol's theorem.

Indeed, $\pi$ is motivic: by Theorem 2.5.1 of \cite{BR}, for any imaginary quadratic field $F'$, there exists a motive
$M_{f,F'}$ defined over
$F\cdot F'$, of rank
$2$ over some extension 
$\Q(f,F')$ of $F'\cdot \Q(f)$; the motives $M_{f,F'}$ are ``{associated to
$f$}'': they give rise to a compatible system of $\lambda$-adic representations of $\Gamma_F$, which is associated
to
$f$. Its Hodge-Tate weights are $0$ and $k_1-1$ above $Id_{F'}$, and $m$ and $m+k_\sigma-1$ above $\sigma\otimes Id_{F'}$.

\noindent{\bf Remark:}
In
fact there should exist
$M_f$ defined over
$\Q$, of rank
$2$ over
$\Q(f)$, associated to
$f$ in the above sense. 

Then we consider for each imaginary quadratic $F'$
$$(7.3.2)\quad M_{\pi,F'}=Res^{F\cdot F'}_{F'} M_{f,F'}$$
$M_{\pi,F'}$ is defined over $F'$, of rank $4$ over $\Q(f,F')$; it is pure of weight $\w=k_1-1$
and the four Hodge-Tate weights $0<m<m+k_\sigma-1<k_1-1$ do occur. These motives define a compatible system of degree $4$
$\lambda$-adic representations of $\Gamma$, associated to
$\pi$.

\noindent{\bf Remark:} Similarly, there should exist $M_\pi$ defined over $\Q$, of rank $4$ over $\Q(f)$ with those Hodge-Tate
weights, associated to
$\pi$.

\vskip 5mm
In the CM case, we restrict our attention to the situation where $f$ is a theta series coming from a biquadratic 
extension
$M=EF/F$, 
$E$ imaginary quadratic. Let ${\rm Gal}(E/\Q)=\{1,\tau\}$,
${\rm Gal}(F/\Q)=\{1,\sigma\}$ and ${\rm Gal}(M/\Q)=\{1,\sigma,\tau,\sigma\tau\}$. We write
$f=\theta(\phi)$ where $\phi$ is a Hecke character of infinity type
$n_1+n_\sigma\sigma+n_{\sigma\tau}\sigma\tau+n_\tau\tau\in \N[Gal(EF/\Q)]$, such that

$$(*)\quad n_1+n_\tau=n_\sigma+n_{\sigma\tau}\quad{\rm and}\quad n_1>n_\sigma>n_{\sigma\tau}>n_\tau$$
and of conductor $\f$ prime to $p$
in
$M$. In that case, one has
$a=n_\sigma-n_\tau-2$, $b=n_1-n_\sigma-1$ and $c=n_1+n_\tau-3$; indeed, since $n_\tau=(c-a-b)/2$, we
see that the condition $n_\tau=0$ is equivalent to $c=a+b$, in which case one has
 $n_1=\w$,
$n_\sigma=k_\sigma-1+m$, $n_{\sigma\tau}=m$ (and $n_\tau=0$). We assume in fact in the sequel a condition slightly stronger
than
$(*)$, namely:

$$(**)\quad \phi^{(1+\tau)\cdot (1-\sigma)}=1\quad{\rm and}\quad n_1>n_\sigma>n_{\sigma\tau}>n_\tau$$ 
Under these assumptions, we say that $f$ is of $(2,2)$-CM type.

\noindent {\bf Remark:} If $(*)$ is satisfied for a character $\phi$, then $(**)$ is satisfied for $\phi^4$.

\vskip 5mm

Let $I=I_f$ be the ring generated by the normalized eigenvalues $a^0_v=\{v\}^{-m\cdot\sigma}\cdot a_v$ ($v$ prime to
$\n$) of $f$ in $\Q(f)$. The $a^0_v$'s are eigenvalues for the divided Hecke operators $T_0(v)=\{v\}^{-m\cdot
\sigma}\cdot T_v$ as introduced by Hida in the beginning of Sect.3 of \cite{Hi}. By Th.4.11 of \cite{Hi}, these eigenvalues
are still integral. Let
$\iota_p:I\hookrightarrow K\subset\overline{\Q}_p$, a
$p$-adic embedding, and $\p_f$ the prime of
$I$ associated to $\iota_p$. We assume hereafter that $p$ splits in $F$: $p\cdot \cO_F=\q\cdot \q^\sigma$, $\q\neq
\q^\sigma$.

 Recall that by a Theorem of Wiles (Th.2.2.2 of \cite{W0}) and a Proposition of Hida (Prop.2.3 of \cite{H}), if
$$ord_p(\iota_p(a^0_\q))=0\quad {\rm resp.}\, ord_p(\iota_p(a^0_{\q^\sigma}))=0$$
(that is, $ord_p(\iota_p(a_\q))=0$ resp. $ord_p(\iota_p(a_{\q^\sigma}))=m$),
then, the decomposition group $D_\q\subset \Gamma_F$ at $\q$ preserving $\iota_p$ is sent by $\rho_{f,\p_f}$, resp.
$\rho_{f_\sigma,\p_f}$ to a Borel subgroup of $GL(2)$; moreover,   $\rho_{f,\p_f}$ resp.  $\rho_{f_\sigma,\p_f}$ restricted to
the inertia subgroup
$I_\q$ has a $1$-dimensional unramified quotient. 

We put $k'=I/\p_f$. Let $J$ be the subring generated by the 
$(a_v,a_{v^\sigma})$ in $\Q(f)\times \Q(f)$.  For $p$ prime
to the index of $I$ in its normalization, and of $J$ in its normalisation, we can view $\rho_{\pi,\p_f}|_F=(\rho_f,
\rho_{f_\sigma})$ as taking values in $GL(2,I_{\p_f})\times GL(2,I_{\p_f})$. Let $X\subset k'^\times$ be the
subgroup generated by the reduction of
$Nv^{k_1-1}\cdot\epsilon(v)$ for all finite places $v$ prime to $\n p$. Let  
$$\overline{\cH}^0=\{(g,g')\in GL_2(k')\times GL_2(k'); det\,g=det\,g'\in X\}$$ 
the two factors being exchanged by $\sigma$, and
$$\overline{\cH}=\{1,\sigma\}\propto\overline{\cH}^0.$$
Similarly, let $\overline{\cH}_{CM}$ be the image by the spin representation of
  $\{g\in \hat{T}(k')\propto W'; \nu(g)\in X\}$

\begin{pro} For $f$ as above and $k_1>k_\sigma>2$, with Nebentypus of order at most $2$, there exists a (non-effective) finite
set
$S$ of finite places in $\Q(f)$ such that, for any
$p\notin S$, splitting in $F$, for which a $\p_f|p$ is ordinary for $f$ and $f_\sigma$, the image of
$\overline{\rho}_{\pi,\p_f}:\Gamma\rightarrow GL_{k'}(\overline{V})$ is equal to:
\begin{itemize}
\item $\overline\cH $, if $f$ is not CM,
\item contains a subgroup of $\cH_{CM}$ of index at most $gcd(p-1,n_1\cdot n_\sigma)$ if $f$ is of $(2,2)$-CM type.
\end{itemize}
\end{pro}

\noindent{\bf Comment:} In other words, let $H$ the subgroup of
$\hat{G}$ whose image by the spin representation is
${}^L (Res^F_\Q GL(2))$ (in the non-CM case) resp. ${}^L(Res^M_\Q M^\times)$ in the $(2,2)$-CM case.
Then, in both cases, the assumption {\bf (RLI)} of sect 7.1, is satisfied for $H$.

\noindent{\bf Proof:} Assume first that $f$ has no CM. We follow the method of proof of Ribet's thesis \cite{R}. More
precisely, we apply Th.3.1 of \cite{R}. We change its statement by replacing 
$\F_p^{k-1}$ there by our subgroup
$X$; since $X\subset \F_p^\times$, the proof of Th.3.1 runs identically.
Let $\overline{G}={\rm Im}\,\overline{\rho}_{\pi,\p_f}|_F$. In order to apply Th.3.1 as in Th.5.1 and 6.1 of
\cite{R}, we have to check
\begin{itemize}
\item (a) For almost all $p$ splitting in $F$ and ordinary as above, $\overline{\rho}_{f,\p_f}$ and
$\overline{\rho}_{f_\sigma,\p_f}$  act irreducibly on
$k'^2$ and their images have order divisible by
$p$,
\item (b) For almost all $p$ as above, there exists $\gamma\in \overline{G}$ such that $(Tr\, \gamma)^2$ generates $k'\times
k'$.
\end{itemize}

(a) If $\overline{\rho}_f$ is reducible, we have
$$\overline{\rho}_f\equiv\left( \begin{array}{cc}\overline{\chi}_1&*\\0&
\overline{\chi}_2\end{array}\right)\,{\rm mod.}\,\p.$$
Let us define a global character $\psi$ of conductor dividing $\n\cdot p$ by
$$\psi_{gal}\cdot\omega^{1-k_1}= \overline{\chi}_1/\overline{\chi}_2 .$$
Let $\psi_{\q}$, resp. $\psi_{\q^\sigma}$ be the restriction of $\psi$ to $I_{\q}$ resp.
to $I_{\q^\sigma}$. By the
ordinarity of
$\rho_f$ at $\q$ and $q^\sigma$, we see that
$\psi_{\q}=1$ or 
$\omega^{2(k_1-1}$ and 
$\psi_{\q^\sigma}=\omega^{2m}$ 
or
$\omega^{2(k_1-1)-2m}$. Let $\epsilon$ be a fundamental unit of $F$. Consider the numbers
$$\epsilon^{2m\cdot\sigma}-1,\epsilon^{[2(k_1-1)-2m]\cdot
\sigma}-1,\epsilon^{2(k_1-1)+2m\cdot\sigma}-1,\epsilon^{2(k_1-1)+[2(k_1-1)-2m]\cdot\sigma}-1;$$ 
If $\q$ is prime to these numbers, we see by global class-field theory that the global character
$\psi$ cannot exist.

{\bf Remark:} This reflects the fact that no congruence between $f$ and an Eisenstein series can
occur, as there are no non-zero Eisenstein series with non-parallel weight.

To assure that $p$ divides the order of ${\rm Im}\,\overline{\rho}_f$, one
proceeds as in Lemma 5.3 of
\cite{R} to exclude all entries of the list of prime-to-$p$ order subgroups in $GL_2(k')$. We have to modify the proof in case
(ii) as follows. Since $\overline{\rho}_f$ is totally odd, we would obtain a totally imaginary quadratic extension $M/F$, of
relative Galois group say,
$\{1,\tau\}$, and a ray-class group character
$\overline{\lambda}:Cl_{M,\f\cdot p}\rightarrow
\overline{\F}_p^\times$ (for some ideal $\f$ of $M$) such that $\overline{\rho}_f=Ind^M_F\,\overline{\lambda}^{gal}$, with
$Norm_{M/F}(\f)D_{M/F}|\n\cdot p$. One can lift $\overline{\lambda}$ into a Hecke character $\lambda$ of $M$ of type adapted to
$k$, so that the theta series $\theta(\lambda)$ belongs to $M_k(\Gamma_0(\n\cdot.p,\epsilon)$ and 

$$(C)\quad f\equiv \theta(\lambda) {\rm mod.}\p$$
 here again, we use the ordinarity of $f$ at $p$:
\begin{itemize}
\item first, if 
$D_{M/F}$ is divisible by $\q$ or $\q^\sigma$, $\theta(\lambda)$ cannot be ordinary at $\q$ (because $k_1$ and $k_\sigma$ are
greater than $1$); therefore the field
$M$ can only ramify above $\n$: this leaves a finite set of possibilities for $M$. 

\item Moreover, by Hida's $p$-stabilization lemma (Lemma 7.1 of Bull.SMF 1995), since $k_1$ and $k_\sigma$ are
greater than $2$ (the cohomological weight $(k_1-2,k_\sigma-2)$ is regular), congruence
$(C)$ can only occur if
$\lambda$ has conductor prime to
$p$. 
\end{itemize} 
In conclusion, if $p$ is prime to all congruence numbers $C_{\theta(\lambda')}$ for
the set of Hecke characters $\lambda'$ of the fields $M$, of the right infinity type, and conductor $\f$ with
$Norm(\f)D_F$ dividing
$\n$ is finite set. If $p$ is prime to the product of the congruence numbers $C(\theta(\lambda'))$ associated to these Theta
series, case (ii) does not occur.

\noindent{\bf Remark:} Note that these congruence numbers
should be given as the algebraic part of the  special value of the Hecke $L$-function $L_M(\lambda'\lambda'^{[\tau]},k)$. This
is the hypothetical converse of a general divisibility result of Hida-T. (Ann.Sci.ENS 1993). It is known at the moment only
for
$F=\Q$ (Hida Inv.64, 1981), but it is conjectured for any totally real field
$F$. 

\noindent   To treat case (iii), we follow closely the argument on p.264 of \cite{R}: if there were infinitely many $\p$
satisfying case (iii), then by using Cebotarev density theorem, one would find a set of positive density of $v$'s  
satisfying  
$a_v^2= 4\cdot Nv^{k_1-1}$. Since
$k_1$ is odd, this condition implies that
$v$ ramifies in $\Q(f)$ or is degree $2$ over $\Q$. This set has density zero in $F$. This is a contradiction. Thus, the set
of
$p$'s in case (iii) must be finite.

(b) As in \cite{R}, we proceed in two steps:
\begin{enumerate}
\item We establish the equality $\overline{G}=\cH^0$ for some prime $\p_f$,

\item We deduce from 1) the existence of $\gamma\in \overline{G}$ as desired for almost all ordinary $p$'s splitting in $F$.
\end{enumerate}

Let $p$ a rational prime, $\p|p$ in $\overline{\Q}$ dividing $\p_f$ and $\q$. We assume that it satisfies (a), that it splits
completely in
$\Q(f)$ and  that
$f$ and
$f_\sigma$ are ordinary at
$\q$. We assume furthermore that for any quadratic Dirichlet character $\chi$ mod.$\n$, there exists
$v$ prime to
$Norm(\n)$
such that $a_v\not\equiv \chi(v)\cdot
a_{v^\sigma}\, {\rm mod}.\p_f$.

These conditions are satisfied if
$\p_f$ is prime to all congruence numbers for all pairs  $f,f_\sigma\otimes\chi$ (for
the Hecke algebra of level $Norm(\n)^2$, generated by Hecke operators outside $Norm(\n)$); indeed the eigensystems of
$f$ and the
$f_\sigma\otimes \chi$, for any $\chi$ mod.$\n$
are mutually distinct. Indeed, if $a_v=a_{v^\sigma}\chi(v)$, for almost all $v$s, then $\chi$ descends to $\Q$. It defines a
quadratic extension $F'/\Q$. Let $E=F\cdot F'$. Let $f_E$ be the base change of $f$ to $E$. If $\tau$ generates $Gal(F'/\Q)$,
the weight of $f_E$ is $k_1(1+\tau)+k_\sigma(\sigma+\sigma\tau)$. The assumption implies that
$f_E=(f_{\sigma})_E=(f_E)_\sigma$; hence $f_E$ should descend to $F'$. This is absurd since its weight is not invariant by
$Gal(E/F)=\{1,\sigma\}$. So these congruence numbers are not zero, and thus can be avoided.

\noindent {\bf Claim:} for such $p$,  $\overline{G}=\overline{\cH}^0$.

\noindent {\bf Proof:} If, not, Th.3.8 of \cite{R} (or rather, its proof) implies that there exists a quadratic
character
$\chi$ of conductor dividing
$\n\cdot p$ such that

$$\overline{\rho}_f\sim \overline{\rho}_{f_\sigma}\otimes\chi.$$

This implies first $a_v\equiv \chi(v)\cdot
a_{v^\sigma}\, {\rm mod}.\p$ for all
$v$'s prime to $Norm(\n)p$. Moreover, by ordinarity of the Galois representations at $p$ (existence of an unramified line),
it also implies that $\chi$ is unramified at
$p$. Since $\chi$ is unramified at $p$, this is a contradiction by the choice of $\p$. 

\vskip 1cm

In fact, for $p$ as above and splitting totally in $\Q(f)$, we even have as in Lemma 5.4 of \cite{R},
a stronger result:

Let $$\overline{\bf \cH}^0=\{(g,g')\in GL_2(I/pI)\times GL_2(I/pI);\,det(g)=det(g')\in X\}$$
and
$$\overline{\bf G}= Im(Gal(\overline{F}/F)\rightarrow \overline{\bf \cH}^0$$
Then,
$$(*)\quad \overline{\bf G}=\overline{\bf \cH}^0.$$

2. Let $p_0$ be a prime
satisfying the conditions of 1 and splitting totally in $\Q(f)$, so that $(*)$ holds.  There exists
$x\in
\overline{\bf \cH}^0$ such that
$Tr(x)^2$ generates
$I/p_0I\times I/p_0I$ over
$\F_{p_0}$. Therefore, by Cebotarev density theorem, there are infinitely many finite places $v$ 
such that the image of
$(a_v^2,a_{v^\sigma}^2)\in I\times I$ in $I/p_0I\times I/p_0I$ generates this ring. For any 
such $v$, by Nakayama's
lemma,$(a_v^2,a_{v^\sigma}^2)$ generates the ring $I_{(p_0)}\times I_{(p_0)}$ over $\Z_{(p_0)}$ , 
hence $\Q(f)\times \Q(f)$
over $\Q$. Fix such a $v$;  let $J=I[(a_v^2,a_{v^\sigma}^2)]$; it is of finite index in $I\times I$.
for any prime $\p$ not dividing the index of $J$ in $I\times I$, we put $\gamma=\overline{\rho}_
{\pi,\p}(Fr_v)$; it belongs to
$\overline{G}$ and $tr(\gamma)^2$ generates $k'\times k'$ over
$\F_p$ (for $k'=I/\p$).  For those
$\p$'s, we conclude that $\overline{G}=\cH^0$.QED.

 In the $(2,2)$-CM case, let $f=\theta(\phi)$. For any $p$ and any $p$-adic field $K$ (with valuation ring $\cO$
and residue field $k$)  containing the field
$\Q(\phi)$ of values of
$\phi$, we still denote by $\phi=\phi^{gal}:Gal(\overline{M}/M)\rightarrow K^\times$ the $p$-adic Galois avatar of the Hecke
character
$\phi$. Thus, we have
$$\rho_\pi=Ind^{M}_\Q(\phi).$$
Let $T\subset G=GSp_4\subset GL(4)$ 
be the standard torus of $G$; the homomorphism $\psi: Gal(\overline{M}/M)\rightarrow GL_4(\cO)$ given by
$\psi=diag(\phi,\phi^\sigma,\phi^{\sigma\tau},
\phi^\tau)$ takes values in $T(\cO)$ by $(**)$. We have $\rho_\pi|_M\cong\psi$. Let $I_\phi$ be the ring of integers of
$\Q(\phi)$; denote by
$k'$ the subfield of
$k=\cO/(\varpi)$ image of $I_\phi$ by the reduction map $\cO\rightarrow k$. 

We claim that for almost all
$p$'s which split totally in
$M$, the image $\Psi$ of
$\psi$ contains a subgroup of index $\leq n_1\cdot n_\sigma$ of $A=\{t\in T(k'); \nu(t)\in X\}$. 

Observe that $\Psi\subset A$ and
$\nu(\Psi)=\nu(A)$. Moreover, since the conductor
$\f$ of
$\phi$ is prime to $p$, we see by class-field theory that the restriction of $\psi$ to the compositum of inertia subgroups
above
$p$ contains
$diag(a^{n_1},b^{n_\sigma},a^{n_1}\cdot b^{-n_\sigma},1)$. Since $k'^\times$ is
cyclic, we conclude. QED

\noindent{\bf Remark:}
Note that in the $(2,2)$-CM case, $p$ is ordinary for
$f$ and $f_\sigma$ at $\p$ if and only if $p$ splits in $M=E\cdot F$.

\vskip 5mm

\begin{cor} if $p\notin S$, splits in $F$, is ordinary for $f$ and $f_\sigma$ (at some $\p_f|p$), and
is greater than $max(5,\w+1)$, $(\pi,p)$ satisfies all the assumptions of Theorems 1 and 2.
\end{cor}  

Calculations communicated to us by H. Yoshida \cite{Y2} establish that the unique level one Hilbert cusp form over
$F=\Q(\sqrt{5})$ of weight $(14,2)$ (hence $m=6$) admits a non-zero cuspidal theta lift $\pi$ which is a classical holomorphic
Siegel cusp form of level $5$ and weight $8$ (that is, $a=b=5$, $c=10$). The motive associated to $\pi$ is rank four
with Hodge weights $0,6,7,13$. 
\begin{itemize}
\item The field $\Q(f)$ is equal to $F$ and the order $I_f$
is maximal. 
\item The prime
$31$ is greater than the motivic weight $\w=13$; 
\item it splits in
$F$:
$$ (31)=\p\p^\sigma,\quad \p=({13+3\sqrt{5}\over 2}),$$ 
\item $\p$ is ordinary for $f$ and $f_\sigma$,
\item  the image $\Psi$ of $\overline{\rho}_\pi$ is equal to
$$ \{1,\sigma\}\propto\{(g,g')\in
GL_2(\F_{31})\times GL_2(\F_{31})| det\,g=det\,g'\in (\F_{31}^\times)^{13}\}.$$
\end{itemize}

The verification of this last point uses Th.3.1 of \cite{R}; the main points are
\begin{itemize}
\item to show, for $\F_{31}=I_f/\p$ that: 

$$\Psi_f={\rm Im}\,\overline{\rho}_f =\{g\in GL(2,\F_{31}); det\,g\in (\F_{31}^\times)^{13}\}.$$

Indeed, $\Psi_f$  contains a unipotent element: consider the degree $2$ prime
$\lambda=(3)$ in $F$; the number
$a_\lambda^2-4N(\lambda)^2$ has order one at $\p$. By \cite{Sh} Lemma 1, this ensures the existence of a
unipotent element. 
$\Psi_f$ is not contained in a Borel: there is a prime $\q$ above $11$ such that $\overline{\rho}_f(Fr_\q)$ is elliptic. 

\item To find a $\gamma\in \Psi$ such that $Tr(\gamma)^2$ generates $I_f/(31)$ over $\bF_{31}$. Take for that the prime
$\q$ above $11$ as above and
$$\gamma=(\overline{\rho}_f(Fr_\q),\overline{\rho}_f(Fr_{\q^\sigma})\in GL_2(\F_{31})\times GL_2(\F_{31}).$$
One has $Tr(\gamma)^2=(28,1)\in \F_{31}\times \F_{31}$, which generates $\F_{31}\times \F_{31}$ over $\F_{31}$ 
.
\end{itemize}

This provides therefore an explicit example of a couple $(\pi,p)$ satisfying all our assumptions.  Other potential examples for
 the same $F$ and $f$ are
$p=19,29$; indeed, they satisfy all the conditions above, except that non-trivial unipotent elements have not been found in the
limit of the calculations of
$a_\lambda^2-4N(\lambda)^2$ (namely, $\lambda$ dividing at most $31$).

Yoshida \cite{Y2} also found that for $F=\Q(\sqrt{13})$, the unique level one Hilbert cusp form
of weight $(10,2)$ lifts to a nonzero holomorphic scalar-valued Siegel cuspform of level $13$, weight $(6,6)$
($a=b=3$) with $\Q(f)=F$, and $I_f$  maximal. The rank $4$ motive associated to $\pi$ has Hodge weights $0,4,5,9$. 
The primes
$p=17$ and $29$ are greater than $\w=9$, split in
$F$;  they are ordinary for $f$ and $f_\sigma$. The image of Galois contains $\{(x,y)\in \F_{p^2}\times\F_{p^2};
N(x)=N(y)\in\F_p^{9}\}\propto \{1,\sigma\}$. However, in the limit of the calculations ($\lambda$ dividing at most $61$) no
unipotent has been found in the image for those primes.
It would be interesting to find examples of cusp forms $f$ of the minimal possible weight, namely $(4,2)$. The theta lift $\pi$
would then occur in middle degree cohomology with constant coefficients: $a=b=0$, and the Hodge-Tate weights of
$\rho_\pi$ would be $0,1,2,3$.

\vfill\eject

\section{Proof of Theorem 2}
The main tool in the proof of Th.2 is the minimal compactification $j:X\hookrightarrow X^*$ (see 8.1 below). This
compactification is far from being smooth (for $g>1$), but it has some advantages over toroidal compactifications; namely, the
strata have a very simple combinatoric and, as a consequence, the Hecke correspondences extend canonically to the boundary.
Let us consider the long exact sequence of the boundary:

$$\begin{array}{ccccccc}\ldots & H^d_c(S_U,V_\lambda(\cO))&\rightarrow & H^d(S_U,V_\lambda(\cO))&\rightarrow
H^d_{\partial}(S_U,V_\lambda(\cO))&\ldots\\& \parallel&&\parallel&\parallel&\\\ldots &
H^d_{et}(X^*_{\overline{\Q}},j_!V_\lambda(\cO))&\rightarrow &
H^d_{et}(X^*_{\overline{\Q}},R j_* V_\lambda(\cO))&\rightarrow H^d_{et}(\partial X^*_{\overline{\Q}},R j_*
V_\lambda(\cO))&\ldots\end{array}
$$

In this section, we shall repeatedly use the standard spectral sequence for an Žtale sheaf $\cal F$ on $X^*$, and a diagram
$j:X\hookrightarrow X^*\hookleftarrow Y:i$
 
$$H^\bullet(Y,i^*R^\bullet j_*{\cal F})\Rightarrow H^\bullet (Y,i*Rj_*{\cal F}).$$
It will allow us to study (localization at $\m$ of)
$H^\bullet (Y,i^*R^\bullet j_*{\cal F})$, rather than the hypercohomology of the complex $i^*R j_*{\cal F})$.

We will thus be left with the study of the Galois action on the boundary cohomology group
$$
H^\bullet_{et}(\partial X^*_{\overline{\Q}},R^\bullet j_*
V_\lambda(\cO))
$$
 in order to show that its localization at $\m$ vanishes.
First, let us recall the description of $ X^*_{{\Q}}$ and the form of the spectral sequence attached to
its stratification. 

\subsection{The minimal compactification} 

The arithmetical
minimal compactification
$X^*=X_g^*$ of $X=X_g$ is defined in non-adelic terms in Th.2.3 of Chapter V of
\cite{CF}.
 It is a normal projective scheme over
$\Z[1/N]$. We are only
interested in the generic fiber $X^*_\Q=X^*\otimes \Q$. In this setting, an adelic definition can be found in
\cite{PinkD} or \cite{Pink} Sect.3 for a general reductive group $G$; let us describe the strata adelically for
$G=GSp(2g)$. We need some notations. For $r=1,\ldots,g$, let $P_r=M_r\cdot
U_{P_r}$ be the standard maximal parabolic of $G$ associated to the simple root $\alpha_{g-r+1}$ (see Sect.3.2.2). Its Levi
group
$M_r$ is isomorphic to
$GL(r)\times GSp(2g-2r)$ (recall that $GSp(0)=\G_m$ by convention). We decompose it accordingly into a product of group
schemes over $\Z$:
$M_r=M_{r,\ell}\times M_{r,h}$, where the index
$\ell$, resp. $h$, denotes the linear, resp. hermitian part of $M_r$. Thus, $M_{r,h}\cong GSp(2g-2r)$ admits a Shimura
variety, which is a Siegel variety of genus $g-r$, while $M_{r,\ell}$ does not. Let
$\kappa_r:P_r\rightarrow M_r=P_r/U_{P_r}$ and let
$P_{r,h}$ be the inverse image of
$M_{r,h}$ by $\kappa_r$. Let $K_{r,h}$ be the standard maximal compact times center in $M_{r,h}(\R)$, and 
$\cZ_{g-r}=M_{r,h}(\R)/K_{r,h}$ be the Siegel space of genus
$g-r$ (it has two connected components $\cZ_{g-r}^{\pm}$ 
); then the compactified symmetric space $\cZ_g^*$ can be described set-theoretically as:
$$\cZ_g^*=\bigsqcup_{r=0}^g G(\Q)\times^{P(\Q)}  \cZ_{g-r}$$

therefore,
 
$$S_U^*=G(\Q)\backslash \cZ_g^*\times G(\A_f)/U.$$

For any subgroup $V_r\subset P_r(\A)$, let us  denote by $V_{r,h}$ its projection
to $M_{r,h}(\A)=P_r(\A)/M_{r,\ell}(\A)\cdot U_{P_r}(\A)$. Then, by simple manipulations we obtain

$$(8.1.1)\quad S_U^*=\bigsqcup_{r=0}^g\bigsqcup_{\dot{x}}
S_{g-r,{}^x U_{r,h}}$$ 

where
\begin{itemize}
\item $\dot{x}$ runs over the finite set $P(\Q)P_{r,h}(\A_f)\backslash G(\A_f)/U$, and $x$ denotes an arbitrary
representative of
$\dot{x}$ in $G(\A_f)$; for later use, we may and do choose $x$ so that its $p$-component $x_p$ is trivial;
\item we have put
${}^x U_r=x\cdot U\cdot x^{-1}\cap P_r(\A)$,
\item we have
$$ S_{g-r,U_r^x}=M_{r,h}(\Q)\backslash M_{r,h}(\A)/{}^x U_{r,h}=M_{r,h}(\Q)\backslash \cZ_{g-r}\times M_{r,h}(\A_f)/{}^x U_{r,h}.$$
\end{itemize}
Note that the disjoint union is set-theoretic, not topological;
see below though.

For each $\dot{x}$, a standard application of ths Strong Approximation Theorem shows that the connected components of
$S_{g-r,{}^x U_{r,h}}$  are indexed by a system
$\{m_{f,h}\}$ of representatives in $M_{r,h}(\A_f)$ of the (finite) set of double cosets $M_{r,h}(\Q)\backslash
M_{r,h}(\A)/{}^x U_{r,h}\cdot M_{r,h}(\R)^+$, where $ M_{r,h}(\R)^+$ denotes the subgroup of $ M_{r,h}(\R)$ of elements with
positive similitude factor. Recall that we have assumed that $U$ is good; the condition $\nu(U)=\hat{\Z}^\times$ implies that
for any $r\geq 1$, the set $M_{r,h}(\Q)\backslash
M_{r,h}(\A)/{}^x U_{r,h}\cdot M_{r,h}(\R)^+$ has only one element. That is, $S_{g-r,{}^x U_{r,h}}$ is connected.

Let 
$$\Gamma_{M_{r,h}}(x)=M_{r,h}(\Q)\cap {}^x U_{r,h}\times M_{r,h}(\R)^+),$$
then, we have a canonical identification

$$S_{g-r,{}^x U_{r,h}}=\Gamma_{M_{r,h}}(x) \backslash \cZ_{g-r}^+$$
this is a Siegel variety of genus $g-r$.

By \cite{PinkD} Sect.12.3, the decomposition $(8.1.1)$ of $S_U^*$ into locally closed subsets canonically descends to $\Q$ into
a stratification of
$X^*_{\Q}$. We have
$$\partial X^*_\Q=X_1\sqcup\ldots \sqcup X_g$$
where the stratum $X_r$ is defined over $\Q$. Actually, 
$$(8.1.2)\quad X_r=\bigsqcup_{\dot{x}} X_{r,x}$$ 
with $\dot{x}\in P(\Q)P_{r,h}(\A_f)\backslash G(\A_f)/U$
and where $X_{r,x}$ is the canonical descent to $\Q$ of $S_{g-r,{}^x U_{r,h}}$.
$(8.1.2)$ is a disjoint union in the Zariski topology.

Recall
For the Zariski topology of $X^*$, one has $\overline{X}_i\supset X_j$ for $i<j$ and 
$$\overline{X}_i-\overline{X}_{i+1}= X_i.$$

\subsection{Spectral sequence associated to the stratification}
To the stratification $\partial X^*_\Q=\overline{X}_1\supset\ldots \overline{X}_g\supset \barX_{g+1}=\emptyset$ is associated a
spectral sequence in Betti or Žtale cohomology
$$(8.2.1)\quad E_1^{p-1,q}=H_c^{p-1+q}(\barX_{p}-\barX_{p+1},k_p^*R j_*V_\lambda(k))\Rightarrow H^{p-1+q}(\partial
X^*_{\overline{\Q}},R j_*V_\lambda(k))$$
where $k_r:X_r\hookrightarrow \partial X^*$ denotes the locally closed embedding of $X_r=\overline{X}_{r}-\overline{X}_{r+1}$. 
It is compatible with algebraic correspondences preserving the stratification. It is mentioned as a remark in Milne, Etale Coh.
Chap.III, Remark 1.30. We don't know a complete reference for it, hence we sketch the proof:
Given a stratification on a scheme $Y$, by closed subsets $Y=Y_0\supset Y_1\supset\ldots\supset Y_{n+1}=\emptyset$,
given a complex of etale sheaves $\bV$ on $Y$ with constructible cohomology, we consider for $p<q$ the
open immersion
$j_{pq}:Y_p-Y_q\hookrightarrow Y_p$ and the closed immersion $i_{pq}:Y_q\hookrightarrow Y_p$. Let $\bV_p=i_{0p}^*\bV$; we have
$\bV_q=i_{pq}^*\bV_p$ for any
$p<q$. We have short exact sequences 
$$0\rightarrow j_{pq,!}\bV_p|_{Y_p-Y_q}\rightarrow \bV_p\rightarrow i_{pq,*}i_{pq}^*\bV_p\rightarrow 0$$
This yields a stratification on the complex $\bV$:
$$0\subset j_{01\,!}(\bV|_{Y-Y_1})\subset j_{02\,!}(\bV|_{Y-Y_2}) \subset\ldots\subset j_{0p\,!}(\bV|_{Y-Y_p})\subset\ldots
\bV$$

Note that for any $p\geq 1$: 
$$j_{0p\,!}(\bV|_{Y-Y_p})/j_{0,p-1\,!}(\bV|_{Y-Y_{p-1}})\cong i_{0,p-1\,*}j_{p-1,p\,!}\bV_{p-1}|_{Y_{p-1}-Y_p},$$
hence,
$$E_1^{p,q}=H^{p+q}_c(Y_{p-1}-Y_p,(i_{0,p-1}^*\bV)|_{Y_{p-1}-Y_p})$$
as desired.

Let us apply this sequence to our stratification. We have for any $r\geq 1$:
$$\barX_r-\barX_{r+1}=\bigsqcup_{\dot{x}} X_{r,x}$$
So,
$$(8.2.2)\quad E_1^{r-1,s}=\bigoplus_{\dot{x}}H^{r-1+s}_c(X_{r,x}, R j_*
V_\lambda(k)|_{X_{r,x}})$$

By the standard spectral sequence
$$H^\bullet_c(X_{r,x}, R^\bullet j_*
V_\lambda(k)|_{X_{r,x}})\Rightarrow H^\bullet_c(X_{r,x}, R j_*
V_\lambda(k)|_{X_{r,x}})$$
We are left with the study of $R^\bullet j_*V_\lambda(k)$.

\subsection{The restriction of the higher direct image sheaf to the strata}

It is easy to determine the restriction mentioned above on the analytic site (in Betti cohomology). The details are in
\cite{Har2} Sect.2.2.5. One finds that the complex restricted to the stratum is quasi-isomorphic to the complex (with trivial
differential) of locally constant sheaves on
$S_{g-r,{}^x U_{r,h}}$ associated to the $\Gamma_{M_{r,h}}( x)$-module:
$$H^\bullet (\Gamma_{M_{r,\ell}}( x),H^\bullet (\Gamma_{U_{P_r}}( x),V_\lambda(k)))$$
where 
$$\Gamma_{M_{r,\ell}}(x)=M_{r,\ell}(\Q)\cap ( {}^{x} U_{r,\ell}\times M_{r,\ell}(\R)),\quad for\,
{}^x U_{r,\ell}=\kappa_r({}^x U)\cap M_{r,\ell}(\A_f)$$
and
$$ \Gamma_{U_{P_r}}(x)=U_{P_r}(\Q)\cap ({}^{ x} U\cap U_{P_r}(\A_f)\times U_{P_r}(\R)).$$
The main result of \cite{Pink} is that, replacing the Betti site by the
Žtale site, this result remains true. More precisely, by Th.(5.3.1) of \cite{Pink}, the complex of Žtale sheaves $R^\bullet
j_*V_\lambda(\F_p)$ over
$X^*_{/\Q}$ restricted to
$X_{r,x\,/\Q}$ is quasi-isomorphic to the complex of etale sheaves with trivial differential
obtained by canonical construction from the representation of
$M_{r,h}\otimes \F_p$ on
$$H^\bullet (\Gamma_{M_{r,\ell}}(x),H^\bullet (\Gamma_{U_{P_r}}(x),V_\lambda(\F_p))).$$   
(and similarly for $k$ instead of $\F_p$). 
We then mention a mod.$p$ version of Kostant decomposition theorem. Recall we have chosen the representatives $x\in G(\A_f)$
so that $x_p=1$. This implies in particular that $\Gamma_{U_{P_r}}(x)$ is dense in $U_{P_r}(\Z_p)$. For any reductive
subgroup $M\subset G$, and any $(M,B\cap M)$-dominant weight $\mu$ of $T\cap M$, let $V_{M,\mu}$ be the Weyl $\Z_p$-module of
highest weight $\mu$ for $M$. 

\begin{lem} Assuming $p-1>|\lambda+\rho|$, then, for any $r\geq 1$, the semisimplification of the $\F_p\Gamma_{M_r}(x)$-module

$$H^q(\Gamma_{U_{P_r}}(x),V_\lambda(\F_p))$$

 is an
$M_r(\F_p)$-module whose decomposition into irreducible
$M_r$-modules is given by:
$$H^q(\Gamma_{U_{P_r}}(x),V_\lambda(\F_p))^{ss}=\bigoplus_{w"\in W^{P_r},\ell(w")=q} V_{M_r, w"(\lambda+\rho)-\rho}$$
\end{lem}

\noindent{\bf Proof:} Over $\Q_p$, the module itself is semisimple and the decomposition is given by Kostant's theorem. By 
Theorem C of \cite{PT}, for
$p$ as stated,
$H^\bullet (\Gamma_{U_{P_r}}(x),V_\lambda(\Z_p))$ is torsion-free.
Therefore $H^\bullet (\Gamma_{U_{P_r}}(x),V_\lambda(\Z_p))$ is a stable lattice in $H^\bullet
(\Gamma_{U_{P_r}}(x),V_\lambda(\Q_p))$. Then, the determination of its composition factors as
 $\Z_p[M_r(\F_p)]$-module, for $p$ as stated, is the content of Cor.3.8 of \cite{PT}.

\vskip 3mm

As already noted, we have $M_r=M_{r,\ell}\times M_{r,h}$. Let $T_\ell=T\cap M_{r,\ell}$ and $T_h=T\cap M_{r,h}$; note that
$T_\ell$ consists in the $t\in T$ of the form 
$$diag (t_g,\ldots,t_{g-r+1},1\ldots,1,t_{g-r+1}^{-1},\ldots ,t_1^{-1}),$$
 while
the maximal torus
$T_h$ of
$M_{r,h}$ consists in the elements
$$t=diag(t_g,\ldots,t_1,\nu\cdot t_1^{-1},\cdots,\nu\cdot t_g^{-1}))\in T$$ 
such that
$t_g=\ldots=t_{g-r+1}=1$. For
$\mu_{w"}=w"(\lambda+\rho)-\rho\in X^*(T)$, we denote the restrictions to $T_\ell$ resp. $T_h$ by
$\mu_{w",\ell}=\mu_{w"}\vert_{T_\ell}$, and $\mu_{w",h}=\mu_{w"}\vert_{T_h}$; since $\mu_{w"}$ is dominant for $(M,B\cap M)$,
$\mu_{w",\ell}$, resp.
$\mu_{w",h}$, is dominant for
$(M_\ell,B\cap M_\ell)$, resp. $(M_h,B\cap M_h)$. By Theorem 1 of
\cite{PT}, it follows from
$p-1>\vert\lambda+\rho\vert$, that the irreducible $M_{r/\Z_p}$-module $V_{M_r, \mu_{w"}}$ can be decomposed
as a tensor product of irreducible $\Z_p$-modules over $M_{r,\ell}$ resp. $M_{r,h}$:
$$V_{M_r,\mu_{w"}}=V_{M_{r,h},\mu_{w",h}}\otimes V_{M_{r,\ell},\mu_{w",\ell}}.$$

Therefore, as $M_{r,h}$-module, we have
$$(8.3.1)\quad  H^\bullet (\Gamma_{M_{r,\ell}}(x),H^\bullet (\Gamma_{U_{P_r}}(x),V_\lambda(\F_p))=$$
$$\bigoplus_{w"\in W^{P_r}}
H^\bullet (\Gamma_{M_{r,\ell}}(x),V_{M_{r,\ell},\mu_{w",\ell}})\otimes V_{M_{r,h},\mu_{w",h}}$$
Thus, the Žtale sheaf on $X_{r,x\,/\Q}$ associated to this representation of $M{r,h}$ is
$$(8.3.2)\qquad \bigoplus_{w"\in W^{P_r}}
H^\bullet (\Gamma_{M_{r,\ell}}(x),V_{M_{r,\ell}},\mu_{w",\ell})
\otimes V_{M_{r,h},\mu_{w",h}}(\F_p)$$
In particular, the Galois action on the Žtale cohomology over $X_{r,x}\otimes \overline{\Q}$
of this sheaf arises only from the second factors of each summand.

\subsection{``{Hodge-Tate weights}'' of the $E_1$-terms}
Recall that $x_p=1$, hence ${}^xU_{r,h}$ is of level prime to $p$, so that $X_{r,x}$ has good reduction at $p$.
For each $r\geq
1$, and each
$w"\in W^{P_r}$, let us determine the Hodge filtration of the crystalline representations
$$H^\bullet_c(X_{r,x}\otimes\overline{\Q}_p,V_{M_{r,h},\mu_{w",h}}(\F_p)).$$
We have ${\rm dim}\,X_{r,x}=d_r={(g-r)(g-r+1)\over 2}$. Since $d_r+\vert \mu_w\vert<p-1$, Faltings' comparison
Th.5.3 of \cite{Fa} applies. Again, as in Sect.7.2, one determines the weights using the modulo $p$ BGG complex
(quasi-isomorphic to de Rham by Cor.1 to Th.6). Let $Q(G_{g-r})$ be the Siegel parabolic of $G_{g-r}=M_{r,h}$ and $M(G_{g-r})$
its standard Levi subgroup. The weights are given by
$$-(w'(\mu_{w"}+\rho_h)-\rho_h)(H_h)=-w'(w"(\lambda+\rho)-\rho+\rho_h)-\rho_h)(H_h)$$
where $w'\in W_{G_{g-r}}^{M(G_{g-r})}$.
By the description of $T_h$ given above, we see that $H_h=H$ and
$w'(-\rho+\rho_h)=-\rho+\rho_h$, hence, the weights are
$$(8.3.1)\quad p(w)=-(w(\lambda+\rho)-\rho)(H)\quad {\rm for}\, w=w'\circ w"$$

\vskip 3mm

\noindent{\bf Claim:}  For $r\geq 1$ and $w"\in W^{P_r}$; let 
$$W_G(w")=\{w\in W_G; w=w'\circ w",\,for\, w'\in
W_{G_{g-r}}^{M(G_{g-r})}\}.$$
Then, the function $W_G(w")\rightarrow \N, w\mapsto p(w)$ cannot take both values $0$ and $\w$.

\vskip 3mm

\noindent{\bf Proof:} As already observed, the function $w\mapsto p(w)$ factors through $W_G/W_M$. 
We see that $p(w)=0$ if and only if $w\in W_M$ and $p(w)=\w$ if and only if $w\in w_0W_M$ where $w_0$ is the longest length
element of $W_G$. The latitude to modify the given $w"$ by $w'$ is too limited: 
If $w"W_{G_{g-r}}^{M(G_{g-r})}\cap W_M\neq \emptyset$, then $w"W_{G_{g-r}}^{M(G_{g-r})}\cap w_0W_M= \emptyset$ and conversely.

\subsection{Hecke algebras for strata}
Let $S$ be a finite set of primes
containing the level of all strata but not containing $p$. Let $\cH(G_g)^S=\bigotimes_{\ell\notin S} \cH(G_g)_\ell$, resp.
$\cH(M(G_g))^S=\bigotimes_{\ell\notin S} \cH(M(G_g))_\ell$ be the abstract Hecke algebras generated over $\Z$ by double classes
at all primes $\ell\notin S$, for $G_g=G$ resp. the Levi
$M(G_g)$ of the Siegel parabolic
$Q(G_g)$. Similarly, for each
$r\geq 1$, we introduce
$\cH(M_r)=\cH(M_{r,\ell})\otimes \cH(M_{r,h})$; by identifying $M_{r,h}$ to $G_{g-r}$, we also introduce
$\cH(M(G_{g-r}))$. For each prime $q\notin S$, by Satake isomorphism, we see that the fraction fields of the $q$-local Hecke
algebras over $\R$ fit in a diagram of finite field extensions:
$$\begin{array}{ccc} Fr(\cH(M(G_g))_q)_\R&\rightarrow &Fr(\cH(M_{r,\ell})_q\otimes \cH(M(G_{g-r}))_q)_\R\\\uparrow&&\uparrow\\
Fr(\cH(G_g)_q)_\R&\rightarrow&Fr(\cH(M_{r\ell})_q\otimes \cH(G_{g-r})_q)_\R\end{array}$$
It corresponds by Galois to the diagram of subgroups of $\S_g\propto
\{\pm1\}^g$:
$$\begin{array}{ccccc}&&\S_r\times
\S_{g-r}&&\\&\swarrow&&\searrow&\\\S_g&&&&\S_r\times(\S_{g-r}\propto\{\pm 1\}^{g-r})\\&\nwarrow&&\nearrow&\\
&&\S_g\propto\{\pm1\}^g&&\end{array}$$
The diagram of fields can be descended from $\R$ to $\Q$ by using twisted action of the Weyl groups as in Sect.VII.1
p.246 of
\cite{CF}. In particular, $\cH(M(G_g))_q$ and $\cH(M_{r\ell})_q\otimes \cH(G_{g-r})_q$ are linearly disjoint over $\cH(G_g)_q$:
$$(8.5.1)\quad Fr(\cH(M_{r,\ell})_q\otimes \cH(M(G_{g-r}))_q)=$$
$$Fr(\cH(M_{r\ell})_q\otimes \cH(G_{g-r})_q)\cdot
Fr(\cH(M(G_g))_q)$$ On the other hand, as a consequence of Satake isomorphism, the Hecke-Frobenius element $U_{q,G}$, resp.
$U_{q,G_{g-r}}$ generates
$Fr(\cH(M(G_g))_q)$ over $Fr(\cH(G_g)_q)$, resp. $Fr(\cH(M(G_{g-r}))_q)$ over $Fr(\cH(G_{g-r})_q)$ (see  Sect.VII.1 
of \cite{CF}). Moreover,
$U_{q,G}=1_{\cH(M_{r,\ell})}\otimes U_{q,G_{g-r}}$.
From $(8.5.1)$, we see that the minimal polynomial
$Irr(X,U_{q,G},\cH(G_g))$ is divisible by $Irr(X,1_{\cH(M_{r,\ell}}\otimes U_{q,G_{g-r}},\cH(M_{r\ell})_q\otimes
\cH(G_{g-r})_q)$.

The Hecke algebra $\cH(G_g)^S$ acts on each stratum $X_r=\bigsqcup_{\dot{x}}X_{r,x}$ by $\Q$-rational algebraic
correspondences.
Indeed, there is a surjective homomorphism of $\Z$-algebras
$$\phi_{g-r}:\cH(M(G_g)^S\rightarrow \cH(M(G_{g-r})^S,$$
$$ 
[G_g(\Z_q)\cdot diag(a_r,b_{2g-2r},c_r)\cdot G_g(\Z_q)]
\longmapsto $$
$$\left\{\begin{array}{cc}[G_{g-r}(\Z_q)\cdot diag(b_{2g-2r})\cdot G_{g-r}(\Z_q)]&{\rm
if}\quad
a_r\in T_{M_{r,\ell}}(\Z_q)\\0&\mbox{\rm if not.}\end{array}\right.
$$
See \cite{Fre}, Sect.IV.3.

On $S_{g-r},{}^xU_{r,h}$, we let the double class $[U\alpha U]$ act by the algebraic correspondence
associated to $\phi_{g-r}([U\alpha U])$. By the theory of canonical models, since $\nu
(U)=\hat{\Z}^\times$, these correspondences are defined over $\Q$.

 Let $\m$ be the maximal ideal of
$\cH(G_g)$ associated to
$\overline{\theta}_\pi$. Let
$$W^{r,s}=E^{r-1,s}_{1\,\m}=(\bigoplus_{\dot{x}}H_c^{r-1+s}(X_{r,x}\otimes\overline{\Q}, R^\bullet j_*
V_\lambda(k)|_{X_{r,x}})[\m]$$

\begin{lem} For any $q\notin S$, the characteristic polynomial of $\overline{\rho}_\pi$ annihilates the action of the geometric
Frobenius $Fr_q$ on $W^{r-1,s}$.
\end{lem}

\noindent{\bf Proof:} By Theorem 4.2, Chapt.VIII of \cite{CF}, we know that 
$$Irr(X, U_{q,G_{g-r}},\cH(G_{g-r})_q)$$ 
annihilates $Fr_q$ on $W^{r-1,s}$. By the divisibility
relation obtained above, we also have
$Irr(X,U_{q,G},\cH(G_g))|_{X=Fr_q}=0$ on  $W^{r-1,s}$. By definition of $\overline{\rho}_\pi$, we have
$char(\overline{\rho}_\pi(Fr_q))=Irr(X,U_{q,G},\cH(G_g))$, as desired.
 
\subsection{End of the proof}
By the previous lemma, we can apply Lemma 13 to $W^{r-1,s}$ (for $r\geq 1$): if $W^{r-1,s}\neq 0$,
both characters $1$ and
$\omega^{-\w}$ occur in $W^{r-1,s}|_{I_p}$. This contradicts the Claim in Sect.8.4. Thus, we have for any $s\geq 0$,
$E_{1\,\m}^{r-1,s}=0$. By $(8.2.1)$ and $(8.2.2)$, we conclude that for any $r\geq 1$ and any $s\geq 0$, $H^{r-1+s}(\partial
X^*, R^\bullet j_*V_\lambda(k))_\m=0$ as desired. By the long exact sequence of cohomology of 
the boundary, we obtain
$H^d_c(X,V_\lambda(\cO))_\m=H^d(X,V_\lambda(\cO))_\m$. We deduce the corollary:
\begin{cor} For $(\pi,p)$ as in Th.1, the natural maps induce an isomorphism
$$H^d_c(X,V_\lambda(\cO))_\m=H^d(X,V_\lambda(\cO))_\m$$ 

\end{cor}
This is the first part of theorem 2. 

\subsection{Intersection cohomology}
For the minimal compactification $j:X\hookrightarrow X^*$ and an etale sheaf $\cal F$ over $X$,
 we consider the intermediate extension $j_{!,*}{\cal F}$. By \cite{BBD}, prop. 2.1.11, 
we have the following description of this complex :

$$j_{!,*}{\cal F}=\tau_{<c_{g}}Rj_{g,*}\tau_{<c_{g-1}}Rj_{g-1,*}\ldots \tau_{<c_{1}}
Rj_{1,*}{\cal F}$$

where for $U_r=\coprod_{0\leq i\leq r} X_i$, we put $j_r:U_{r-1}\hookrightarrow U_r$ , $i=1,\ldots,g$,
  and $c_r$ is the codimension of the stratum $X_r$ in $\overline{X}_{r-1}$.

We have 
$$\begin{array}{ccccccc}\ldots & H^d_c(S_U,V_\lambda(\cO))&\rightarrow &
IH^d(S_U,V_\lambda(\cO))&\rightarrow IH^d_{\partial}(S_U,V_\lambda(\cO))&\ldots\\&
\parallel&&\parallel&\parallel&\\\ldots & H^d_{et}(X^*_{\overline{\Q}},j_!V_\lambda(\cO))&\rightarrow &
H^d_{et}(X^*_{\overline{\Q}},j_{!,*} V_\lambda(\cO))&\rightarrow H^d_{et}(\partial
X^*_{\overline{\Q}}, j_{!,*} V_\lambda(\cO))&\ldots\end{array}
$$

\begin{pro} $IH^\bullet_{\partial}(S_U,V_\lambda(\cO))_\m=0$. 
\end{pro}

The proof will be similar to the usual cohomology case: it relies on Pink's theorem, lemma 13 and
a variant of Claim 8.4. Some more induction is needed though, due to the successive truncations
involved in defining
$j_{!,*}V_\lambda$.

By the spectral sequence (sect.8.2) associated to our stratification, we are reduced to show
$$H_{c,et}^\bullet (X_{r,x},j_{!,*}V_\lambda(k))_{\m_r}=0$$

\begin{lem}
$H_{c,et}^d(X_{r,x},j_{!,*}V_\lambda(\F_p))$ is a subquotient of 
$$H^\bullet(X_{r,x},R^\bullet j_{g,*}\circ
R^\bullet j_{g-1,*}\ldots\circ R^\bullet j_{1,*}V_\lambda(\F_p))$$
\end{lem}

\noindent {\bf Proof:} We write the argument for $g=1$ and $2$.
For $g=1$, it follows directly from the second spectral sequence associated to the complex
$\tau_{<c_1}Rj_{1,*}V_\lambda(\F_p)$:
$$H_c^\bullet(X_r,\tau_{<c_1}R^\bullet j_{1,*}V_\lambda(\F_p))\Rightarrow
H_c^\bullet(X_r,\tau_{<c_1}R j_{1,*}V_\lambda(\F_p)).$$ 

For $g=2$, applying this ``{second spectral sequence}'' to
$$\tau_{<c_2}R j_{2,*}(\tau_{<c_1}R j_{1,*}V_\lambda(\F_p))),$$
 The group $H_{c,et}^d(X_{r,x},j_{!,*}V_\lambda(\F_p))$ admits a d\'evissage by subquotients of
$$H_c^\bullet(X_r,\tau_{<c_2}R^\bullet j_{2,*}\tau_{<c_1}R j_{1,*}V_\lambda(\F_p)).$$
The complex inside the cohomology is filtered, hence the cohomology itself is
filtered and its graded pieces are subquotients of
$$H_c^\bullet(X_r,\tau_{<c_2}R^\bullet j_{2,*}\tau_{<c_1}R^\bullet j_{1,*}V_\lambda(\F_p))$$
by the formalism of spectral sequences.

For the next lemma, some more notations are needed.
Let $j_{pq}:U_p\hookrightarrow U_q$ for
$p<q$; thus,
$j=j_{0,g}=\ldots =j_r\circ j_{0,r}$. Let $i_{p,q}:X_p\hookrightarrow  U_q$ denotes the
locally closed immersion of
$X_p$ in $U_q$ (composition of the closed immersion $i_p:X_p\hookrightarrow U_p$ followed by $j_{p,q}$).
Note that $j_{0,r}=i_{0,r}$.

Let 
$$W(r)=\prod_{s=0}^r W_{G_{g-s}}^{P_s}$$
(so, $W(0)=\{1\}$).
For $w(r)=(w_r,\ldots,w_1)\in W(r)$, the symbol 
$w(r)\cdot (\lambda+\rho(r))$ is defined by induction by
$$w(r+1)\cdot (\lambda+\rho(r+1))=w_{r+1}\cdot (w(r)\cdot((\lambda+\rho(r)))+\rho_{r+1}).$$
(recall that $\rho_r$ denotes the half-sum of positive roots of $G_{g-r}$ for the order deduced from
$(G_g,B_g,T_g)$) and $w_r\cdot(\lambda+\rho_r)=w_r(\lambda+\rho_r)-\rho_r$.
One sees by induction on $r$ that $\vert \lambda+\rho\vert_r<p-1$
implies $\vert w(r)\cdot (\lambda+\rho(r))\vert_r<p-1$ for any $r\geq 0$.

For $\underline{\alpha}=(\alpha_g,\ldots,\alpha_1)\in \Z^g$, $\alpha_i\geq 0$, let
$$\delta_{\underline{\alpha},r}=\left\{\begin{array}{cc} 0& {\rm if}\,
\alpha_{r+1}+\ldots +\alpha_g>0\\1&\mbox{\rm if not} \end{array}\right.$$

\begin{lem} The sheaf $R^{\alpha_g}j_{g-1,g,*}\circ\ldots\circ R^{\alpha_1}j_{0,1,*}V_\lambda(\F_p)$ is constructible finite Žtale; 
for $r=0,\ldots,g$, its restriction to the stratum $X_r$ comes by the canonical construction from a $\Gamma_{M_r}(x)$-module
whose semisimplification is

$$\delta_{\underline{\alpha},r}\cdot \bigoplus_{\underline{w(r)}\in J}
 V_{w(r)\cdot
(\lambda+\rho(r))}\otimes T_{\underline{\alpha},r,\lambda}$$ 
for some subset $J\subset W(r)$ and where $T_{\underline{\alpha},r,\lambda}$ is the Žtale sheaf associated to a finite-dimensional $\F_p$-vector
 space without
action of $G_{g-r}$; it depends only of $(\alpha_r,\ldots,\alpha_1)$; it does not
contribute to Galois action on the Žtale cohomology of $X_r$.
\end{lem}

If $F=can(V)$, we shall write $F^{ss}$ for the sheaf $can(V^{ss})$

\noindent {\bf Proof:} For each $r$, we consider the abelian category $\cC_r$ of constructible Žtale sheaves 
in $\F_p$-vector spaces
over $U_r$; let $\cA_r$ be the abelian category generated by the 
$j_{s,r,!}i_{s,s,*}F_s$ ($0\leq s\leq r$)
 where $F_s$ is locally constant on $X_s$, of the form stated in the lemma. Since these sheaves
are supported by the strata $X_s$ and since there are no non-zero morphisms between sheaves with
disjoint support, $\cA_r$ consists exactly in the objects mentioned.

Let $\cB_r$ be the  abelian subcategory of $\cC_r$ stable by extension generated by $\cA_r$. We first note
that the sheaves of the form $G=j_{r-i,i,!}i_{s,r-i,*}F_s$, $0\leq s\leq  r-i$
 are objects of $\cB_r$. Indeed, we have the short exact sequence:
$$0\rightarrow j_{r-i-1,r,!}i_{s,r-i-1,*}F_s\rightarrow G\rightarrow
j_{r-i,r,!}i_{r-i,r-i,*} i_{r-i,r-i}^*G\rightarrow 0 $$  

By Pink's theorem, the sheaf on the right is in $\cB_r$; by decreasing induction on $i$, the sheaf on
the left is in $\cB_r$ since for $i=r-s$, we have $j_{s,r,!}i_{s,s,*}F_s\in \cB_r$.  In particular,
the sheaves $i_{s,r,*}F_s$ are objects of $\cB_r$.

\vskip 5mm

\noindent{\bf Remark:}  If we could prove that any finite $\F_p\Gamma_{M_r}(x)$-module with 
$p$-small highest weight is algebraic, it
would follow from \cite{PT} that it is semisimple. Then, $\cA_r$ and $\cB_r\vert_{X_r}$ 
would be semisimple.

\vskip 5mm

To prove the lemma, it is enough to show by induction on
$r\geq 0$ the following statement

$$
(P_r)\qquad R^{\alpha_r}j_{r-1,r,*}\circ
\ldots\circ R^{\alpha_1}j_{0,1,*}V_\lambda(\F_p)
\in\cB_r$$

Indeed, if we assume $(P_r)$ and if we apply $R^{\alpha_g}j_{g,*}\circ
R^{\alpha_{g-1}}j_{g-1,*}\ldots\circ R^{\alpha_{r+1}}j_{r,r+1,*}$, we obtain the extra factor
$\delta_{\underline{\alpha},r}$ as desired.  

$(P_r)$ is obvious for $r=0$ (take $T_{\underline{\alpha},0,\lambda}=\F_p$). 

\vskip 5mm

For $r=1$, let
$\R_1=R^{\alpha_1}j_{0,1,*}V_\lambda$; we know that
$$\R_1^{ss}\vert_{X_1}=\bigoplus_{w_1\in W_{G_{1}^{P_1}}}V_{w_1\cdot
(\lambda+\rho_1)}\otimes T_{\underline{\alpha},1,\lambda}$$
by Pink's theorem.
Therefore, we have an exact sequence on $U_1$:
$$0\rightarrow j_{0,1,!}V_\lambda\otimes T_0 \rightarrow\R_1^{ss}\rightarrow i_{1,1,*}\bigoplus_{w_1\in
W_{G_{g-1}^{P_1}}}V_{w_1\cdot (\lambda+\rho_1)}\otimes T_1\rightarrow 0$$
for some $T_0$ and $T_1$ as desired (in fact $T_0=\F_p$ if $\alpha_1=0$ and $0$ otherwise).

\noindent {\bf Induction step:} Assume  that $(P_{r-1})$ holds. Note that $R^{\bullet}j_{r-1,r,*}$
 preserves $\cC_r$.
Let 
$$\R_{r-1}=R^{\alpha_{r-1}}j_{r-2,r-1,*}\ldots\circ R^{\alpha_1}j_{0,1,*}V_\lambda(\F_p).$$
By assumption there is a filtration $F^\bullet\R_{r-1}$ whose graded pieces are in $\cA_{r-1}$.

Hence $R^\bullet j_{r-1,r,*}\R_{r-1}$ will be in $\cB_r$ if for each $s$ between $0$ and $r-1$:  

$$(8.6.1)\qquad R^{\bullet}j_{r-1,r,*}j_{s,r-1,!}i_{s,s,*}F_s\qquad
\mbox{\rm is in}\,\cB_r$$

We prove this statement by induction on the length of the stratification of $X^*$. If true for $g-1$, 
we consider $X^*$ of length $g$. 
We can assume $s=0$, and we have to prove that 
$R^{\bullet}j_{r-1,r,*}j_{0,r-1,!}F_0 \in \cB_r$.
We prove in the Appendix that such a sheaf is constructible with respect to the natural stratification 
of $X^*$. Therefore, it remains only to show that for each $s\leq r$, the locally constant sheaf
$$R^{\bullet}j_{r-1,r,*}j_{0,r-1,!}F_0\vert_{X_s}$$
involves the desired highest weight representations.

We have the following
exact sequences
$$0\rightarrow j_{t,r-1,!}j_{0,t,*}F_0 \rightarrow j_{t+1,r-1,!}j_{0,t+1,*}F_0
\rightarrow j_{t+1,r-1,!}i_{t+1,t+1,*}F_{t+1}\rightarrow 0$$
where $t=0,\ldots,r-2$ and $F_t=i_t^*(j_{0,t,*}F_0)$. By induction hypothesis, we have
$R^{\bullet}j_{r-1,r,*} j_{t+1,r-1,!}i_{t+1,t+1,*}F_{t+1} \in \cB_r$. So,
$R^{\bullet}j_{r-1,r,*}j_{0,r-1,*}F_0$
is constructible for the natural stratification, and we need to see it has the correct highest weight
representations (in brief, is of type $\cB_r$).
  
This sheaf is the $E_2^{\bullet,0}$-term in the spectral sequence of composition of two functors
abutting at
$$R^\bullet j_{0,r,*}F_0$$

By  Sublemma 1 below, this abutment is of type $\cB_r$.
Let us check that for $q>0$,
$$E_2^{p,q}= R^p j_{r-1,r,*}R^q i_{0,r-1,*}F_0$$
belongs to $\cB_r$. We notice that $R^q i_{s,r-1,*}V_{w(s)\cdot(\lambda+\rho(s))}$ is supported on
$X_{1}\coprod\ldots
\coprod X_{r-1}$,
 hence we can apply the induction assumption to $X_1^*$ which has a stratification of length $g-1$; 
we obtain 
$$If\quad q>0,\qquad E_2^{p,q}\in \cB_r.$$
 
The conclusion follows then from sublemma 2.

\begin{sublem} Let $X^*$ be a space with a stratification $\Sigma$ of length $g$. For each $r=0,\ldots,g$,
let
$\cA_r$ be an abelian subcategory of locally constant sheaves on
$X_r$; assume that for any $s\leq r\leq g$, $i_r^* R^\bullet i_{s,*}$ sends $\cA_s$ to $\cA_r$. Let
$\cB$ be the smallest abelian category of $\Sigma$-constructible Žtale sheaves on $X^*$ which is stable by
extensions (that is, which is thick) and contains $j_{s,!}i_{s,s,*}F_s$ (for $s=0,\ldots, g$).Then 
$R^\bullet j_*$ sends $\cA_0$ to $\cB$.
\end{sublem}

\noindent {\bf Proof:} Let $V_0\in \cA_0$ and $F=R^\bullet j_* V_0$.

Consider the filtration 
$$F_g=j_!F\vert_{U_0}\subset\ldots\subset
F_r=j_{r,!}F\vert_{U_{g-r}}\subset\ldots F_0=F$$ 

The successive quotients are given by 

$$F_{i-1}/F_i \cong j_{g-i+1,!}i_{g-i+1,*}i^*_{g-i+1}F_{i-1}.$$ 
Note that $i^*_{g-i+1}F_{i-1}=i^*_{g-i+1}F$ belongs to $\cB$ by assumption.

\vskip 5mm

We conclude by the following trivial lemma.

\begin{sublem}
 Let $\cB$ be a full thick abelian subcategory of an abelian category $\cC$ which is stable by subobjects
and quotients. Let
$E_2^{p,q}\Rightarrow H^{p+q}$ in
$\cC$
be a spectral sequence
 concentrated in $p,q\geq 0$. Assume that
$E_2^{p,q}\in
\cB$ for any $E_2^{p,q}$, $q\neq q_0$, and $E_\infty^{p,q}\in \cB$ for any $p,q$, then $E_2^{p,q_0}\in \cB$.  
\end{sublem}

\noindent {\bf Proof:} 
By decreasing induction on the $r$ of the spectral sequence $E_r^{p,q}$.

From these two lemmata,  th.2.(ii) will follow if we show

\begin{lem} For any $s=1,\ldots, g$, we have
$$\bH=H^\bullet_c(X_s,V_{w(s)\cdot(\lambda+\rho(s))})_{\m_s}=0.$$
\end{lem}

\noindent{\bf Proof:} We apply Claim 8.4;  a formula similar to the formula there shows that the Hodge-Tate weights 
occuring in $\bH$ are
$$-w'_s\cdot w"_s\cdot \ldots w'_1\cdot w"_1\cdot (\lambda+\rho(s))(H)$$
that is,
$$ p(w)=-(w(\lambda+\rho)-\rho)(H)\quad {\rm for}\, w=w'_s\circ w"_s\circ \ldots w'_1\circ w"_1$$

As in 8.4, since $s\geq 1$, $0$ and $\w$ cannot
occur simultaneously  as weights for this cohomology group. On the other hand, by the Galois-theoretic argument 8.6 they
should, if
$\bH\neq 0$ by Lemma 13. We conclude
$\bH=0$.

\vskip 5mm

It is maybe useful to state in a single result an outcome of our proof of Theorems 1 and 2:

\begin{cor} Under the assumptions for $\pi,p, \m$  as before, we have:

$$H^\bullet_c(S_U,V_\lambda(\cO))_\m=IH^\bullet(S_U,V_\lambda(\cO))_\m=
H^\bullet(S_U,V_\lambda(\cO))_\m=H^d(S_U,V_\lambda(\cO))_\m.$$
\end{cor}

\noindent{\bf Comment:} This corollary requires {\bf (RLI)}, but does not require the regularity of
$\lambda$. When $\lambda$ is regular, we have already mentioned that
$$H^\bullet_{cusp}(S_U,V_\lambda(\C))=IH^\bullet(S_U,V_\lambda(\C))=H^\bullet_!(S_U,V_\lambda(\C))=
H^d_!(S_U,V_\lambda(\C)).$$
moreover, it seems plausible that for such a $\lambda$, for any $q<d$, $H^q(S_U,V_\lambda(\C))=0$. 
It might
result from Franke spectral sequence. It does indeed for $g=2$ (see Appendix A of \cite{TU}). If it were
true, harmonic analysis would provide a complex version of our theorem, without localization :

For , $q<d$,
$$ H^q_{cusp}(S_U,V_\lambda(\C))=IH^q(S_U,V_\lambda(\C))=H^q_!(S_U,V_\lambda(\C))=
H^q(S_U,V_\lambda(\C))=0$$
and
$$ H^d_{cusp}(S_U,V_\lambda(\C))=IH^d(S_U,V_\lambda(\C))=H^d_!(S_U,V_\lambda(\C)).$$
But of course
$$H^d_!(S_U,V_\lambda(\C))\neq H^d(S_U,V_\lambda(\C)).$$

\vfill\eject

\section{Application to a control theorem}  In this section, we want to apply Theorem 1 for improving
upon Theorem 6.2 of \cite{TU}. More precisely, we want to replace the non effective assumption on the
prime $p$ there, (namely, $p$ prime to the order of the torsion subgroups of $H^q(S_U,V_\lambda(\Z))$ for
$q=1,2,3$) by an ``{effective}'' assumption $p-1>max(a_2+a_1+3,4)$ which in particular is independent of the level (however,
we need to assume the mod. $p$ non-Eisensteiness condition {\bf (RLI)} which is far from being effective, but depends only on
$\overline{\rho}_\pi$). Note however that we need to localize the Hecke algebra at the maximal ideal given by
$\theta_\pi$ modulo ${\varpi}$. This is innocuous for questions of congruences between $\theta_\pi$ and characters
coming from other representations occuring in
$H^3$.

\bigskip

We
prefer to treat axiomatically the general case $G=GSp(2g)_\Q$ of an arbitrary genus $g$, assuming
conjectures (which are proven for $g=2$). Most notations in this section follow those of Section 7 of
\cite{TU}. Let
$\lambda=(a_g,\ldots,a_1;c)$ be a  dominant regular weight ({\it i.e.} $a_g>\ldots>a_1>0$) and
$\pi$ a cuspidal representation of level
$U$ occuring in $H^d(S_U,V_\lambda(\C))$. Recall that $B$
denotes the standard Borel subgroup $B$ of $G$ and $B^+$ its unipotent radical. Let $p$ be a prime not dividing $N$. for any
$n\geq  1$, let
$$U_0(p^n)=\{g\in U; g\,mod.p^n\in B(\Z/p^n\Z)\}$$
resp.
$$U_1(p^n)=\{g\in U;
g\,mod.p^n\in B^+(\Z/p^n\Z)\}$$
The $p$-component of $U_0(p^n)$ is the Iwahori subgroup (resp. strict Iwahori subgroup) of level $p^n$ and is denoted by
$I_n\subset G(\Z_p)$, resp. $J_n\subset G(\Z_p)$.
Let 
$S_1(p^n)$ resp. $S_0(p^n)$ be the Siegel variety associated to $U_1(p^n)$ resp. to 
$U_0(p^n)$.
For each $n\geq 
1$, let
$$\cW^q_{\lambda,n}=H^q(S_1(p^n),V'_\lambda(K/\cO))$$
where $V'_\lambda$
denotes the Iwahoric induction of
$\lambda$ that is the lattice in $V_\lambda (K)$  consisting in $\lambda^{-1}$-equivariant rational functions $f$ on
$G/B^+$ taking
integral values on the Iwahori subgroup $I_1$ of $G(\Z_p)$. Thus $V'_\lambda$ is $I_1$-stable (hence $J_n$-stable for any
$n\geq 1$). Note that it contains the
$G(\Z_p)$-stable lattice $V_\lambda$ defined similarly, but with the stronger condition $f(G(\Z_p))\subset \cO$. Let
$\cW^q_\lambda$  be the inductive limit over $n\geq  1$ of the $\cW^q_{\lambda,n}$. 

Let $\cW^\bullet_{\lambda,n}=\bigoplus\cW^q_{\lambda,n}$, resp.
$\cW^\bullet_\lambda=\bigoplus\cW^q_\lambda$. We introduce several abstract Hecke algebras:
Let 
$$D_p=\{d\in T(\Q_p)\cap M_{2g}(\Z_p)^{prim}| ord_p(\alpha(d))\leq 0\,\mbox{\rm for any positive root}\,\alpha\}$$
where $M_{2g}(\Z_p)^{prim}$ denotes the set of integral matrices with relatively prime entries. $D_p$ is a semigroup. Let
$\cH^{N}$, resp.$\cH^{N,I_n}$, resp.
$\cH^{N,J_n}$ be the abstract Hecke $\cO$-algebra outside $N$ and integral at $p$, resp. integral at $p$ of type $I_n$, resp.
integral at $p$ of type
$J_n$:
$$\cH^N=\bigotimes_{\ell\,prime\,to\,Np}\cO[G(\Q_\ell)//G(\Z_\ell)]\otimes
\cO[U_pD_pU_p//U_p],$$
$$\cH^{N,I_n}=\bigotimes_{\ell\,prime\,to\,Np}\cO[G(\Q_\ell)//G(\Z_\ell)]\otimes
\cO[I_nD_pI_n//I_n],$$
$$\cH^{N,J_n}=\bigotimes_{\ell\,prime\,to\,Np}\cO[G(\Q_\ell)//G(\Z_\ell)]\otimes
\cO[J_nD_pJ_n//J_n].$$
 
For any $n\geq 1$, there is a natural surjective homomorphism $\cH^{N,J_n}\rightarrow \cH^{N,I_n}(\cO)$, but
that there is no homomorphism $\cH^{N,I_1}(\cO)\rightarrow \cH^N$. Assume that $\pi$ satisfies the condition {\bf
(AO)} of automorphic ordinarity at $p$ (see introduction). Let us recall how one can transfer the character
$\theta_\pi:\cH^N\rightarrow \cO$ to a character
$\theta'_\pi:\cH^{N,I_1}\rightarrow
\cO$. The inclusion of lattices $V_\lambda\subset V'_\lambda$, together with the finite morphis $S_0(p)\rightarrow S_U$ give
rise to a morphism of sheaves $(S_U,V_\lambda(\cO))\rightarrow (S_0(p),V'_\lambda)$, hence a morphism on cohomology
$$\iota:H_*^\bullet(S_U,V_\lambda(\cO))\rightarrow H_*^\bullet(S_0(p),V'_\lambda(\cO)).$$
Moreover, the Hecke operators $T_{p,i}$, $i=1,\ldots,g$, defining the condition {\bf (AO)} act on these cohomology groups.
Observe however that for each $i$, $T_{p,i}$ act differently in prime-to-$p$ level (e.g. on $S_U$), and in level $p$ (e.g. on
$S_0(p)$). They define idempotents on these cohomology groups; let $e_0=\lim_{n\rightarrow\infty} (\prod_{i=1}^g
T_{p,i})^{n!}$ be the idempotent defined on $H_*^\bullet(S_U,V_\lambda(\cO))$, and $e=\lim_{n\rightarrow\infty} (\prod_{i=1}^g
T_{p,i})^{n!}$ defined on $H_*^\bullet(S_0(p),V'_\lambda(\cO))$ by the same formula (with a different meaning though).

\begin{lem} (Hida's stabilization lemma) If $\lambda$ is regular, the homomorphism 
$$H_*^\bullet(S_U,V_\lambda(\cO))\rightarrow H_*^\bullet(S_0(p),V'_\lambda(\cO)),\quad x\mapsto e.\iota(x)$$
induced by the diagram
$$\begin{array}{ccc}H_*^\bullet(S_U,V_\lambda(\cO))&\rightarrow& H_*^\bullet(S_0(p),V'_\lambda(\cO))\\\bigcup&&e \downarrow\\
e_0\cdot H_*^\bullet(S_U,V_\lambda(\cO))&& e\cdot H_*^\bullet(S_0(p),V'_\lambda(\cO))\end{array}$$
is an isomorphism sending an eigenclass for $\cH^N$ to an eigenclass for  $\cH^{N,I_1}$.
\end{lem}
  
\noindent{\bf Proof:} See Prop.3.2 of \cite{TU} (proven there for $GSp(4)$ over a totally real field: it generalizes directly
to arbitrary
$g$).

Denote by $\bh_\lambda(U;\cO)$, resp.$\bh_\lambda(U_1(p^n);\cO)$, resp. $\bh_\lambda(U_0(p^n);\cO)$,  the image of
$\cH^N$ in ${\rm End}_\cO(H^\bullet(S_U,V_\lambda(\cO)))$, resp. of $\cH^{N,J_n}$ in  ${\rm
End}_\cO(\cW^\bullet_n)$, resp.
$\cH^{N,J_n}$ in ${\rm
End}_\cO(H^\bullet(S_0(p^n),V'_\lambda(\cO)))$.
 By the lemma above for $*=\emptyset$, the character
$\theta_\pi:\bh_\lambda(U;\cO)\rightarrow \cO$ induces a character $\theta'_{\pi}:\bh_\lambda(U_0(p);\cO)\rightarrow \cO$;
hence (compatible) characters of
$\bh_\lambda(U_1(p^n);\cO)$ for any $n\geq 1$. Let
$$\bh_\lambda =projlim_n \bh_\lambda(U_1(p^n);\cO).$$
Note that $\bh_\lambda$ acts faithfully on $\cW^\bullet$.
Let    
$\m'=Ker\,\overline{\theta}'_{\pi}$ be the maximal ideal of $\bh_\lambda$ 
associated
to
$\pi$.
The localization
$\cW^q_\lambda (\m')$ of
$\cW^q_\lambda$, resp.$\cV^q_\lambda$ at
$\m'$ is contained in the ordinary part
$e\cdot\cW^q_\lambda$ and is therefore a localization of this ordinary
part. Note that
$T(\Z_p)\subset D_p$; by action on $\cW^q_{\lambda,n}$, we obtain (compatible) group homomorphisms
$$\langle\,\rangle_\lambda: T(\Z_p)\rightarrow \bh_\lambda(U_1(p^n);\cO).$$
By linearization, we obtain a continuous $\cO$-algebra homomorphism from the completed group algebra $\cO[[T(\Z_p)]]$ to
$\bh_\lambda$. For any discrete
 $\cO[[T(\Z_p)]]$-module $\cW$, the Pontryagin dual $\cW^\star=Hom(\cW,K/\cO)$ is a compact 
topological $\cO[[T(\Z_p)]]$-module. Let
$$T_1=Ker(T(\Z_p)\rightarrow T(\F_p))\quad {\rm and}\, \Lambda=\cO[[T_1]]$$
$\Lambda$ is an Iwasawa algebra in $(g+1)$-variables.
 Recall that an arithmetic character
$\chi:T(\Z_p)\rightarrow \cO^\times$ is a product $\chi=\varepsilon\mu$ where
$\varepsilon$ is of finite order, factoring through, say, $T(\Z/p^n\Z)$ and $\mu\in X^\star(T)$ is
algebraic. If  $\chi\equiv 1\,{\rm mod.}\,{\varpi}$, it can be identified to a
character of $T_1$. It induces canonically an $\cO$-algebra homomorphism $\chi:\Lambda\rightarrow
\cO$. Its kernel $P_\chi$ is a prime ideal of
$\Lambda$ called an arithmetic prime. We say that $\chi=\mu\varepsilon$ is dominant regular if
$\mu$ is.      

\begin{thm} Given a $\pi$ cuspidal of level $N$; let
$p$ be a prime not dividing $N$ such that the conditions {\bf (Gal)}, {\bf (RLI)}, {\bf (AO)} and {\bf (GO)} hold, 
and that $
p-1>max(a_1+\ldots+a_g+d,4)$;  then 

(i) $\cW^\bullet_\lambda(\m')=\cW^d_\lambda(\m')$ and
$\cW^d_\lambda(\m')^\star$ satisfies the exact control theorem:
 for any regular dominant arithmetic $\chi$, there is a canonical isomorphism  
$$H^d(S_0(p^n),V'_{\lambda\otimes\chi}(K/\cO))_{\m'}\rightarrow \cW^d_\lambda
(\m')[\chi]$$   
Same result for the compactly supported version  $\cC\cW_\lambda (\m')$ of $\cW_\lambda
(\m')$ and for its image $\cW^d_{!,\lambda}(\m')$ in $\cW_\lambda (\m')$.

(ii) The inclusion $\cW^d_{!,\lambda}(\m')\subset \cW^d_\lambda(\m')$ is an equality.

(ii) $\cW^d_\lambda(\m')^\star$ is free of finite rank over
$\Lambda$.
\end{thm}

\noindent{\bf Proof:} 

(i) The proof makes use of Hida's Exact Control criterion (Lemma 7.1 of \cite{H1}) together with the
calculations of Section 3 of \cite{TU} which generalize readily to $GSp(2g)_\Q$. We prove
$\cW^q_\lambda(\m')=0$ and $\cC\cW^q_\lambda(\m')=0$  by induction on $q<d$. For that, by Theorem 3.2 (ii) and isomorphism
(3.16) of \cite{TU}, it is enough to show that
$H^q(S_0(p),V'_\lambda(K/\cO))_{\m'}=0$. By Proposition 3.2 of \cite{TU} and its proof (relating
$\m'$ and $\m$), this amounts to see $H^q(S_U,V_\lambda(K/\cO))_\m=0$. This is precisely what
is stated in Theorem 1 in the introduction, under our assumptions. Thus, exactly as in the proof of
Theorem 3.2 of \cite{TU}, we obtain (i) for $\cW^q$. In an exactly similar manner, we show the control for the
compact support analogue, based on the Exact Control criterion for compactly supported cohomology. 

(ii) Similarly, the degree $d$
boundary cohomology is controlled, and vanishes in weight $\lambda$ (i.e. $\chi=1$) by our Main Th.2.
Therefore, by Nakayama's lemma, it vanishes $\Lambda$-adically, and
$\cW^d_{!,\lambda}(\m')= \cW^d_\lambda(\m')$.

(iii) We use the following criterion: a discrete $\Lambda$-module $\cW$ is $\Lambda$-cofree of corank
$r<\infty$ if and only if there exists an infinite set of arithmetic characters $\chi$ such that
$\bigcap_\chi P_\chi=0$ in $\Lambda$, and for which
$\cW[\chi]$ is $\cO$-divisible, cofree of constant corank $r$. We take the set of algebraic
dominant characters $\chi=\mu\lambda^{-1}$ with $\mu$ regular
dominant and congruent to $\lambda$ mod. $p$, and apply the control formula stated in (i). We need to see that
$H^d(S_0(p),V'_{\mu}(K/\cO))_{\m'}$ is $p$-divisible (and furthermore, of constant
corank). The long exact sequence
$$ H^d(S_0(p),V'_{\mu}(K))_{\m'}\rightarrow
H^d(S_0(p),V'_{\mu}(K/\cO))_{\m'}
\rightarrow $$
$$
\rightarrow H^{d+1}_c(S_0(p),V'_{\mu}(\cO))_{\m'}$$  

shows it is enough to verify
that the
$H^{d+1}$ is torsion-free. By Poincar\'e-duality (Th.6.4 of \cite{TU}), it amounts to see that
$H^{d-1}_c(S_0(p),V'_{\hat{\mu}}(K/\cO))_{\m'}$ is divisible; in fact it is
null because by (i), since $\hat{\mu}$ is regular dominant, one knows that 
$\cC\cW^{d-1}_{\hat{\lambda}}(\m')$ is zero  and that it is controlled:

$$H^{d-1}_c(S_0(p),V'_{\hat{\mu}}(K/\cO))_{\m'}=\cC\cW^{d-1}_{\hat{\lambda}}(\m')
[\hat{\chi}]=0.$$ 
This shows the divisibility of $\cW^d_\lambda(\m')[\chi]$ for all 
$\mu$'s as above. The corank $r(\chi)$ can be read off from the dimension over the residue field $k$ of the
${\varpi}$-torsion. Note that in $\Lambda$, $P_{\chi}+({\varpi})$ is the maximal ideal, hence does not
depend on $\chi$. Thus $r(\chi)={\rm dim}_k \cW^d_\lambda(\m')[\m_{\Lambda}]$ is
independent of $\chi$.QED.  
    
Let $\bh_\m=\bh_\lambda (U;\cO)(\m')$ be the localization of $\bh_\lambda$ at $\m'$. It acts faithfully on
$\cW^\bullet_\lambda(\m')=\cW^d_\lambda(\m')$.

\begin{thm} Under the same assumptions,

(i) $\bh_\m$ is a finite torsion-free $\Lambda$-algebra,

(ii) there exists a finite integrally closed extension $\bI$ of $\Lambda$ and a $\Lambda$-algebra
homomorphism $\bh_m\rightarrow \bI$ such that for any
$\mu\in X$ such that $\mu\equiv\lambda\, mod.p$ and
$\phi=\mu\lambda^{-1}$ is dominant regular, for $P$ a prime in $\bI$ above $P_\phi$ and
$\cO'=\bI/P$, there is a commutative diagram

$$ \begin{array}{ccc}\bh_\m/P_\phi\bh_\m&\rightarrow
&\cO'\\\downarrow&\nearrow&\\\bh_\mu (U;\cO)_\m&&\end{array}$$
where the horizontal arrow is $\Theta\otimes Id_{\bI/P}$ and the oblique arrow is $\theta_{\pi_P}$
for some cuspidal automorphic representation $\pi_P$ occuring in
$H^d(S_U,V_\mu(\C))$. For $\mu=\lambda$, one has $\theta_{\pi_P}= \theta_\pi$ on $\cH^N$.

(iii) If $\pi'$ is another cuspidal representation occuring in $H^d(S_U,V_\lambda(\C))$, if
$\theta_\pi\equiv\theta_{\pi'}\,mod.max(\overline{\Z}_p)$, there exists another finite integrally
closed extension $\bI'$ of $\Lambda$ and a $\Lambda$-algebra homomorphism
$\Theta':\bh_\m\rightarrow \bI'$ lifting $\theta_{\pi'}$ and for any $\mu$ and any arithmetic
ideal $P''$ in the compositum $\bI\cdot\bI'$; let 
$P=P"\cap\bI$ and $P'=P"\cap \bI'$; we have 
$$\theta_{\pi_P}\equiv \theta_{\pi'_{P'}}\,mod.max(\overline{\Z}_p)$$

\end{thm}

{\bf Comment:} 1) We call $\Theta$ a Hida family in $(g+1)$-variables lifting $\theta_\pi$. Statement
(iii) means that congruences to $\theta_\pi$ (outside $N$) can be lifted to families of congruences. 

2) Statement (i) implies that $\bh_\m$ is flat of relative dimension $(g+1)$ over $\cO$; this was
predicted by calculations in Sect.9, example 2, and Sect. 10.5.3, Conjecture I, of \cite{T}; it was
already proven
$g=2$ in
\cite{TU} under stronger assumptions on $p$.

3) The representations $\pi_P$ occuring in the family whose existence is stated in (iii) are
cuspidal because $\bh_\m$ is
cuspidal: by Th.9 (ii), $\cW^d_{!,\lambda}(\m')=\cW^d_\lambda(\m')$ for any $\mu$ as in the theorem,
$H^d(S_U,V_\mu(\cO))_\m\subset H^d_{cusp}(S_U,V_\mu(\C))$ by our Th.2 and the considerations at the end of Sect. 2.1.  

\noindent {\bf Proof:} It results from the previous one as in Corollary 7.5-7.7 of \cite{TU}.

\vfill\eject

\section{Application to Taylor-Wiles' systems}

 In this section, we apply Theorem 1 to show that some
cohomology group $M_Q$ is free over a finite group algebra
$\cO[\Delta_Q]$ (this is the non-trivial condition to be verified for having a Taylor-Wiles' system:
Condition (TW3) of Definition 1.1 in \cite{Fu}, see also Proposition 1 of \cite{TW}. More precisely,
let us fix as above a cuspidal stable representation
$\pi$ whose finite part $\pi_f$ occurs in
$H^d(S_U,V_\lambda (\C))$,  for a regular dominant weight $\lambda$. Let $p$ be a prime at which
the level group
$K$ is unramified. Let $r\geq  1$. We consider sets
$Q=\{q_1,...,q_r\}$ consisting of primes $q$ which are congruent to $1$ mod.$p$  and such that the
four roots of
$\overline{\theta}_\pi(P_{q}(X)$ are distinct and belong to $k$. For each $q\in Q$, we fix one of these
roots and denote it by $\alpha_q$. Let $(\Z/q\Z)^\times=\Delta_{q}\times (\Z/q\Z)^{(p)}$ where
$\Delta_{q}$ is the $p$-Sylow subgroup and
$(\Z/q\Z)^{(p)}$ the non-$p$-part of $(\Z/q\Z)^\times$. Let
$\Delta_Q=\prod_{q\in Q}\Delta_{q}$. We put
$$U_Q=\{g\in U; {\rm for\,any\,}q\in Q,\,
g\equiv\left(\begin{array}{cccc}u&\ddots&*&\ddots\\0&*&*&*\\0&0&u^{-1}&*\\0&0&0&*\end{array}\right)
\,{\rm mod.}\,q\},\quad u\in (\Z/q\Z)^{(p)}\}$$ 
and
$$U_0(Q)=\{g\in U_Q; {\rm for\,any\,}q\in Q,\, g\,{\rm mod.}q\in B(\Z/q\Z)\}$$ Let ${\cal H}_Q$ be
the abstract Hecke algebra for $U_Q$ generated over $\cO$ by 
\begin{itemize}
\item  Hecke operators $T$'s
 outside 
$$S_Q=Ram(U)\cup\{p\}\cup Q$$
\item  the $U_q$'s for each $q\in Q$:
$$ U_q=U_Q\cdot diag(1,\ldots,1,q,\ldots,q)\cdot U_Q
$$
\item and by the normal action of $\Delta_Q=K_0(Q)/K_Q$.
\end{itemize} 
$\theta_\pi:{\cal H}_Q\rightarrow \cO$ resp. $\overline{\theta}_\pi:{\cal H}_Q\rightarrow k$ define
$\cO$-algebra homomorphisms. Let 
$$\m_Q=<{\varpi},T-\theta_\pi(T), (T\,{\rm outside}\, S_Q), U_q-\alpha_q, (q\in Q)>.$$  It is a maximal
ideal of ${\cal H}_Q$. Consider the following ``{$d$-th homology module}'':   

$$M_Q=H^d( S_{U_Q},V_\lambda(K/\cO))_{\m_Q}^\star$$

It has a natural action of $\cO[\Delta_Q]$. It
is a complete intersection noetherian local ring. 
\begin{thm} Assume that {\bf (Gal)}, {\bf (RLI)} and {\bf (GO)} hold, and
$p-1>max(\vert\lambda+\rho\vert,4)$; then, for any $Q$ as above
 $M_Q$ is free over $\cO[\Delta_Q]$.
\end{thm}

{\bf Proof:} By Theorem 1, we know that $M_Q$ is free as $\cO$-module. Hence, it is enough to show
that
$\overline{M}_Q=M_Q/{\varpi}\cdot M_Q$ is free over $\Lambda_Q=k[\Delta_Q]$. By Pontryagin duality, 
$\overline{M}_Q$ is the $k$-dual of the ${\varpi}$-torsion submdodule $N_Q$ of 
$H^d(S_{U_Q},
V_\lambda (K/\cO))_{\m_Q}$. By the long exact
sequence for
$$0\rightarrow  V_\lambda ({\varpi}^{-1}\cO/\cO)\rightarrow  V_\lambda (K/\cO)\rightarrow V_\lambda
(K/\cO)\rightarrow  0$$ 
and the vanishing
 of $H^{d-1}({}_{U_Q}S, V_\lambda (K/\cO))_{\m_Q}$, we see that
$$N_Q=H^d(S_{U_Q},
V_\lambda (k))_{\m_Q}.$$  

Moreover,
$\Lambda_Q$ is complete intersection, hence is a Frobenius algebra: the freeness of $\overline{M}_Q$ is
equivalent to that of $N_Q$. 

To show that $N_Q$ is free, we follow Fujiwara's approach (Sect.3 of \cite{Fu1}). Since $\Lambda_Q$
is artinian local, freeness is equivalent to flatness:
$Tor_j^{\Lambda_Q}(N_Q,k)=0$ for $j>0$.
 For any $\ell$ prime to $N$, consider the sub-semigroup $D'_{Q,\ell}$ of $T(\Q_\ell)\cap M_{2g}(\Z_\ell)_{prim}$
consisting in $t$'s such that $ord_\ell(\alpha (t)\leq 0$ for any positive root $\alpha$ of $(G,B,T)$.
 Let $D_{Q,\ell}=U_{Q,\ell}\cdot D'_{Q,\ell}\cdot
U_{Q,\ell}$. For 
$q\in Q$, the local Hecke algebra $\cH_{Q,q}=\Z[U_{Q,q}\backslash D_{Q,q}/U_{Q,q}]$ is generated by
$$\Delta_q \, {\rm and}\, diag(1,q^{a_2},\ldots,q^{a_g},q^{c-a_g},\ldots,
q^{c-a_2},q^c),\, {\rm for} \, 0\leq a_2\leq\ldots\leq a_g\leq
c/2\} .$$ 

Note that 
$$\cH_Q=\bigoplus_{\ell\notin S_Q} \cH_\ell^{unr}\oplus (\bigoplus_{q\in Q} \cH_{Q,q})$$

We view $V_\lambda(k)$
as an \'etale sheaf over $X_Q=S_{U_Q}\otimes \Q$. For $t\in T(\A^N)$ and $t_\ell\in D'_{Q,\ell}$, the Hecke
correspondence $[U_Q t U_Q]$  acts on $(X_Q,V_\lambda(k))$ via the diagram
$$(10.1) \quad \begin{array}{ccccccc}&&S_{U_Q\cap t^{-1}U_Qt}&\cong&S_{U_Q\cap tU_Qt^{-1}}&&\\&\pi_1\swarrow&&&&\searrow\pi_2&\\
S_{U_Q}&&&&&& S_{U_Q}
\end{array}$$
where $\pi_1$ and $\pi_2$ are the canonical coverings induced by the inclusions of the level groups,
 the horizontal isomorphism is induced by right
multiplication by
$t^{-1}$. The action on the sheaf $V_\lambda(k)$ is via $\pi_{1,*}\circ [t^{-1}]\circ \pi_2^*$, where
 $[t^{-1}]:\pi_2^*V_\lambda(k)\rightarrow\pi_1^* V_\lambda(k)$ is induced by a right action of the $p$-component $t^{-1}$ 
on the representation $V_\lambda$ which preserves integrality: see for instance \cite{TU} Section 3.5. 

We can form a complex $C^\bullet$ representing $R\Gamma(X,V_\lambda(k))$
endowed with an action of $Gal(\overline{\Q}/\Q)\times \cH_Q$. One can take for instance the global sections 
$C^\bullet(X_Q,V_\lambda(k))$ of
the Žtale Godement
resolution 
$$\cC^\bullet(X,V_\lambda(k))$$
 of $V_\lambda(k)$ (see \cite{FK} Sect.12,
p.129, and Section 3.4 \cite{Fu}) whose terms are acyclic. More precisely, by functoriality of the construction,
the diagrams (10.1) still operate on 
$(X_Q,\cC^\bullet )$ and induce
endomorphisms
$[U_QtU_Q]$ of $C^\bullet$.
 The diagrams (10.1) are defined over $\Q$, hence the action of Galois by transport of structure commutes to these
endomorphisms. The main property that we shall use for the Godement resolution is the following. Let
$f:X\rightarrow Y$
be a finite Žtale
 Galois covering
 with Galois group $G$, let $\cG$ be an Žtale sheaf  on $Y$, let $\cC^\bullet(Y,\cG))$, resp. $\cC^\bullet
(X,f^*(\cG)$
be the Godement resolution 
of $\cG$ resp. $f^*\cG$ on $Y$ resp.$X$. $G$ acts on $f_*\cC^\bullet
(X,f^*(\cG)$ and the adjunction map $a:\cG\rightarrow f_*f^*\cG$ induces an isomorphism

$$(f_*\cC^\bullet(X,f^*\cG))^G=\cC^\bullet(Y,\cG).$$

In particular for $q\in Q$ and
$G=\Delta_q$, we shall make use of the formula

$$(10.2)\quad (C^\bullet(X_Q,V_\lambda(k)))^{\Delta_q}=C^\bullet(X_Q/\Delta_q,V_\lambda(k)).$$

 The hypercohomology spectral sequence applied to
$C^\bullet\otimes_{\Lambda_Q}k$ gives rise to the
$Tor$-spectral sequence:
$$E_2^{i,j}=Tor^{\Lambda_Q}_{-i}(H^j(C^\bullet),k)\rightarrow H^{i+j}(C^\bullet\otimes k)$$ All the
maps involved are $k[Gal(\overline{\Q}/\Q)]\times {\cal H}_Q$-linear. Let us tensor this spectral
sequence with the localized Hecke algebra ${\cal H}_{Q,\m_Q}$. We get
$$E_2^{i,j}({\m_Q})=Tor^{\Lambda_Q}_{-i}(H^j(C^\bullet)_{\m_Q},k)\rightarrow
H^{i+j}(C^\bullet\otimes k)_{\m_Q}$$

{\bf Fact:} $H^j(C^\bullet)_{\m_Q}=0$ for any $j\neq d$.

{\bf Proof:} By Theorem 1, we know that
 
$$H^j(S_{U_Q}\otimes \overline{\Q}, V_{\lambda}(k))_{\m_Q}=0\quad for\, j>d.$$

This fact implies that the spectral sequence is concentrated on
$E_2^{i,d}({\m_Q})=Tor^{\Lambda_Q}_{-i}(N_Q,k)$ and therefore degenerates: 

$$H^{i+d}(C^\bullet\otimes_{\Lambda_Q}k)_{\m_Q}=E_2^{i,d}({\m_Q}).$$ 

It remains to see that
$H^{i+d}(C^\bullet\otimes_{\Lambda_Q}k)_{\m_Q}=0$ unless $i=0$.

For this purpose, we consider the exact
sequence of complexes
$$(10.3)\quad 0\rightarrow \prod_{q\in Q}(C^\bullet)^{\Delta_q}\rightarrow (C^\bullet)^{\oplus Q}\rightarrow
(C^\bullet)^{\oplus Q}\rightarrow C^\bullet\rightarrow C^\bullet\otimes_{\Lambda_Q}k\rightarrow 0$$

 where for each $q\in Q$, the $q$-th component of the middle arrow is the multiplication by $\delta_q-1$
on $C^\bullet$, for $\delta_q$ a generator of $\Delta_q$.
 By Theorem 1 of this paper and by (10.2),
 we see that the first
four complexes of (10.3) have no $\m_Q$-localized cohomology in degree $>d$. By considering long exact sequences,
and by exactness of $\m_Q$-localization, this implies that the same holds for the complex of
$\Delta_Q$-coinvariants
$C^\bullet_{\Delta_Q}=C^\bullet\otimes_{\Lambda_Q}k$.  This concludes the proof.

\vfill\eject

\section{Appendix I: On the constructibility of certain Žtale sheaves}
Let $X^*$ be the minimal compactification over $\Q$ of the Siegel variety $X$ over the rationals. Let
$\Sigma$ be the standard stratification on $X^*$; the strata have dimension $c_r=r(r+1)/2$, $r=g,g-1,\ldots, 0$. Let
$r\geq 0$ and $U_r$ be the union of the strata of dimension greater than $c_r$; we write $\Sigma_r$ the
stratification on $U_r$ induced by $\Sigma$. Let
$j_r:U_r\hookrightarrow X^*$ be the natural open immersion. The goal of this appendix is to provide a proof for the following
proposition which is used in Sect.8.7 for proving Lemma 18.

\begin{pro}For any $\Sigma_r$-constructible torsion
Žtale sheaf
$V$ on
$U_r$, for any $i\geq 0$, $R^i j_{r,*}V$ is $\Sigma$-constructible.  
\end{pro}

\noindent{\bf Proof:} Since $r$ is fixed, we abbreviate $j_r=j$. We use a smooth toroidal compactification
of
$X$.  Let $U$ be the level group of our Siegel variety. Let
${\bf S}=({\bf S}_\xi)_\xi$ be a
$U$-admissible regular rational polyhedral cone decomposition of $S^2(\Z^g)$ (see \cite{CF} Chap.IV, Th.6.7 and \cite{PinkD}
Sect.12.4); in the above notation, $\xi$ runs over the set of rational boundary components in the minimal compactification
$X^*$ and ${\bf S}_\xi$ is a polyhedral cone decomposition of $S^2(N_\xi)$ for a quotient $N_\xi$ of $\Z^g$ of rank $r_\xi$,
depending only on $\xi$ (here, $r_\xi$ is the genus of the Siegel variety $\xi$). Let
$X_{{\bf  S}}$ be the corresponding toroidal compactification of $X$ over $\Q$. It is smooth and $X_{{\bf  S}}-X$ is a divisor
with normal crossings, whose irreducible components are  smooth; it is endowed with a proper morphism
$\pi:X_{{\bf  S}}\rightarrow X^*$ defined over $\Q$, inducing the identity on $X$. The toroidal stratification
$\{Z(\sigma)\}_{\sigma\in {\bf S}/GL(X)}$ is compatible to (and finer than) the inverse image $\pi^{-1}(\Sigma)$ of the
stratification $\Sigma$ (see Th.6.7 of \cite{CF}).  By
\cite{CF} Chap. IV.3 or
\cite{Pink} 3.10, the restriction $\pi_\xi$ of
$\pi$ above any rational boundary component $\xi$ of $X^*$ is a proper morphism with singularities of smooth dnc type: let
$F_\xi=X_{\bf S}\times_{X^*}\xi$, then, locally for the Žtale topology, we have $\cO_{F_\xi}\cong
\cO_{\xi}[T_1,\ldots,T_m]/(T_1\cdot\ldots\cdot T_n)$. 
More precisely, $F_\xi$ is a disjoint union 

$$F_\xi=\bigcup_{\sigma\in {\bf T}_\xi}Z(\sigma)$$ 
where 
\begin{itemize}
\item ${\bf T}_\xi$ is the set of cones $\sigma\in {\bf S}_\xi$ whose elements are all definite positive on $N_\xi$,

\item $Z(\sigma)=\Xi_\xi\times^{E_\xi}Z_\xi(\sigma)$ (in the notations of \cite{CF} p.106) are the toroidal strata.
\end{itemize}

Note that ${\bf T}_\xi$ has the property that any cone of ${\bf S}_\xi$ containing a cone in ${\bf T}_\xi$ is in ${\bf
T}_\xi$; therefore,
$F_\xi$ is closed in the toric immersion $\Xi_{\xi,{\bf S}_\xi}$.  
Moreover, the $Z(\sigma)$ are
smooth as well as their closures; thus, $F_\xi$ is Žtale-locally the boundary of a toric immersion of $E_\xi$ for $T_\xi$,
of smooth dnc type, as desired.

 Let
$U_{r,{\bf  S}}$ be the inverse image of $U_r$ by $\pi$, and $j_{{\bf  S}}:U_{r,{\bf  S}}\hookrightarrow X_{{\bf  S}}$ the
corresponding open immersion. We have $\pi\circ j_{{\bf  S}}=j\circ \pi$. Similarly, let $k:X\hookrightarrow U_r$ resp
$k_{{\bf  S}}:X\hookrightarrow U_{r,{\bf  S}}$. By a simple dŽvissage, one can assume that our Žtale sheaf is
of the form $V=k_!W$ for a locally constant sheaf $W$ on $X$. Then, we have
$$k_!W=\pi_*\circ k_{{{\bf  S}},!}W$$
Let $V_{\bf S}=k_{{{\bf  S}},!}W$.
We have $R^q\pi_*V_{\bf  S}=0$ if $q>0$, by proper base change. Hence,
$R^ij_*\circ \pi_*V_{\bf  S}=R^i(j_*\circ\pi_*)V_{\bf  S}=R^i(\pi_*\circ j_{{\bf  S},*})V_{\bf  S}$ which is the abutment of a
spectral sequence whose $E_2$-term is $R^p\pi_*\circ R^q j_{{\bf  S},*} V_{\bf  S}$.

We show now that the sheaves 
$R^q
j_{{\bf  S},*} V_{\bf  S}$ are constructible for the natural toroidal stratification. By compatibility of the
toroidal stratification of $X_{\bf S}$ with that of the toric immersion of $E=Hom(S^2(\Z^g),\G_m)$,  we can view
$X\hookrightarrow U_{r,{\bf S}}\hookrightarrow X_{\bf S}$, local-etally as $E\hookrightarrow E_r(\sigma)\hookrightarrow
E(\sigma)$ where $E=\G_m^N$, $E_r(\sigma)=\G_m^{(N-n)}\times \A^n$ and $E(\sigma)=\A^N$. We are now in a cartesian product
situation, and therefore, by KŸnneth formula, we are left with the one-dimensional case $\G_m{\stackrel{k'}\hookrightarrow}
\G_m{\stackrel{j'}\hookrightarrow}
\A^1$ or
$\G_m{\stackrel{k'}\hookrightarrow} \A^1{\stackrel{j'}\hookrightarrow} \A^1$. It is easy then to see that $R^ij'_*$ of $k'_!V$
is constructible.

By Lemma 21 below, the higher direct images $R^p\pi_*(R^q j_{{\bf  S},*} V_{\bf  S})$ are
$\Sigma$-constructible.

 This property being preserved by subquotients, we conclude that $R^ij_*\circ \pi_*V_{\bf  S}$
is $\Sigma$-constructible.

\vskip 1cm

\begin{lem} Let $Y$ be an integral scheme over $\Q$ and $f:X\rightarrow Y$ be a proper morphism of smooth dnc type. Let
$T=(X_0,X_1,\ldots ,X_n)$ be the stratification of $X$ defined by $X_0=X^{smooth}$, $X_{i+1}=(\overline{X_i}-X_i)^{smooth}$.
Let $\cF$ be a
$T$-constructible torsion Žtale sheaf on $X$. Then $R^i f_*\cF$ is locally constant.
\end{lem}

\noindent{\bf Proof:} By properness of $f$, we know that  $R^i f_*\cF$ is constructible on $Y$ with finite fibers. To check it
is locally constant we proceed by induction on dimension of $X$; the maps 
$$X_0{\stackrel{j}\hookrightarrow}
X{\stackrel{i}\hookleftarrow} X_1$$

 provide a dŽvissage:

$$0\rightarrow j_!\cF\vert_{X_0}\rightarrow \cF\rightarrow i_*i^*\cF\rightarrow 0$$

By stability of locally constant sheaves by kernels and extensions, we are left with the case of
$$R^i f_*j_!\cF\vert_{X_0}.$$
By a theorem of M.Artin (exposŽ XII \cite{SGA4}, see also Illusie's Appendix, p. 252-261 in \cite{SGA4etdemi}) this sheaf is
locally constant (in general, we would need that
$\cF\vert_{X_0}$ is tamely ramified along the divisor with normal crossings $\overline{X_0}-X_0$ for a smooth compactification
$X_0\hookrightarrow
\overline{X_0}$ over $Y$, but it is automatic here, since we are in characteristic $0$).

\vfill\eject

\section{Appendix II: An explicit construction of the log crystal $\overline{\cV}_\lambda$}

In this appendix, we use Weyl's invariant theory to construct automorphic vector bundles over $\Z_p$, associated to  dominant
weights of the symplectic group $G=GSP_{2g}$ and of the Levi $M$ of the Siegel parabolic of $G$. The defect of this method
(comparing with that of section 5.2) is the lack of functoriality. The advantage is to show clearly how the Hodge structure is
obtained by plethysms from that of $R^1f_*\Omega^\bullet_{A/X}$.

  As before, $X$ is the natural smooth model of
$S_U$ over
$\Z_{(p)}$, $\overline{X}$ is a toroidal compactification over $\Z_{(p)}$. It is
projective smooth and its divisor at infinity $D$ has normal crossings. Let $f:A\rightarrow X$ be  the universal principally
polarized $g$-dimensional abelian variety over $X$; let 
$Y=A\times_X\ldots\times_X A$ be the fiber product of
$A$ by itself $s$-times above $X$ and $f_s:Y=A^s\rightarrow X$ its structural map. Let us recall some
facts on algebraic correspondences.

\subsection*{II.1 Correspondences over $\Z_{(p)}$} We view $f:A\rightarrow X$ over $\Z_{(p)}$ for a
prime $p$ not dividing $N$. Let
$s\geq 1$. Let $Z^\bullet (Y/X)$ be the free abelian group generated by irreducible closed
$X$-subschemes
$Z\subset Y\times_X Y$, flat over
$X$. It is graded by the relative codimension of cycles. Its quotient $A^\bullet(Y\times_X Y/X)$ by
the submodule of cycles on
$Y\times_X Y$ rationally equivalent to zero is denoted by
$Corr^\bullet(Y/X)$ and is called the group of correspondences on $Y$ relative to $X$ (\cite{Ful}
Section 20.1). By smoothness of
$f_s:Y\rightarrow X$ and of
$X$ over
$\Z_{(p)}$, the group $Corr^\bullet(Y/X)$ carries a natural structure of graded ring (see Ex.20.1.1
(c) and Ex. 20.2.3 of
\cite{Ful}). 

Let $C^\bullet (Y/X)_{(p)}=C^\bullet (Y/X)\otimes\Z_{(p)}$. 

A correspondence $Z\in Corr^r (Y/X)_{(p)}$ gives rise (because of the smoothness of the base $X$
over $\Z_{(p)}$)
 to
a cohomology class

$$Cl(Z)\in R^{2r}(f_s\times f_s)_*\Omega^\bullet_{Y\otimes Y/X}$$
 defined by the relative cycle
map (See \cite{EZ} Chap.IV). Let $\delta=g\cdot s= dim\, Y$.

We follow \cite{Kl}, Sect 3 in a relative setting:  by KŸnneth formula and Poincar\'e duality, we
have
$$ R^{2r}(f_s\times f_s)_*\Omega^\bullet_{Y\otimes Y/X}= \bigoplus_{0\leq m\leq
2r}Hom_{\cO_X}(R^{m+2\delta-2r}f_{s,*}\Omega^\bullet_{ Y/X},R^mf_{s,*}\Omega^\bullet_{ Y/X})$$

We can therefore view the $m$-th component of $Cl(Z)$ as a degree $2r-2\delta$ endomorphism of
$R^\bullet f_{s,*}\Omega^\bullet_{ Y/X}$. This defines a homomorphism

$$Corr^\bullet (Y/X)_{(p)}\rightarrow End_{\cO_X}R^\bullet f_{s\,*}\Omega^\bullet_{Y/X}$$ 

which corresponds to letting a
cycle $Z$ act by "$pr_{1\,*}\circ pr_2^*$" on the sheaf $R^\bullet f_{s\,*}\Omega^\bullet_{Y/X}$.
More precisely, we have : 

\begin{lem}
Let $u\in R^{*}(f_s\times f_s)_*\Omega^\bullet_{Y\otimes Y/X}$, then $u(x)=pr_{1\,*}( pr_2^*(x)\cup
u)$.
\end{lem}

{\bf Proof :} \cite{Kl} Sect.3.

This
homomorphism sends cycles $Z$ of relative codimension $\delta+r$ ($-\delta\leq r\leq \delta$)
to degree $2r$ endomorphisms. We denote by 
$$\cC(Y/X)=\bigoplus_{-\delta\leq r\leq \delta} \cC^{2r}(Y/X)  $$
 the graded algebra generated by the
cycle classes of correspondences; it is a finite free $\Z_{(p)}$-algebra. 

In particular, we can view cycles $D$ of $Y$ as cycles in $Y\times_X Y$ via the
diagonal immersion
$Y\hookrightarrow Y\times_X Y$ (the two resulting projections $pr_i: D\rightarrow Y$ are equal).
This yields
$$A^r(Y/X)\rightarrow Corr^{r+\delta} (Y/X)_{(p)}\rightarrow End_{\cO_X}R^\bullet
f_{s\,*}\Omega^\bullet_{Y/X}.$$
Write $D\mapsto [D]$ for this homomorphism.
On the other hand, the action of the cycle $D$ by $-\cup Cl(D)$ yields another homomorphism
$$A^r(Y/X)\rightarrow  End_{\cO_X}R^\bullet f_{s\,*}\Omega^\bullet_{Y/X}$$

\begin{lem} Let $\iota :Y\rightarrow Y\times_X Y$ be the diagonal immersion and $\Delta$ its image.
Then for any cycle $D$ of $Y$, we have 
$$Cl_{Y\times Y}(\iota_*D)=\iota_*Cl_Y(D)=pr_1^*(Cl_Y(D))\cup Cl_{Y\times
Y}(\Delta)=$$
$$pr_2^*(Cl_Y(D))\cup Cl_{Y\times Y}(\Delta)$$
\end{lem}

\noindent {\bf Proof :} By the functoriality of the cycle class map 
we have the following commutative diagram:
$$\begin{array}{lcr} A^r(Y/X)&\rightarrow &R^{2r}
f_{s\,*}\Omega^\bullet_{Y/X}\\ \iota_*\downarrow&&\iota_*\downarrow\\
Corr^{r+\delta}(Y/X)&\rightarrow&R^{2r+2\delta}
f_{s\,*}\Omega^\bullet_{Y\times_X Y/X}
\end{array}$$
where the horizontal arrows are the cycle maps, the left vertical arrow exists by properness of
$\iota$ and the right vertical one is the Poincar\'e dual of $\iota^*$. It remains to check that the
$\iota_*$ on the right satisfies 
$$\iota_*(x)=pr_1^*(x)\cup Cl_{Y\times Y}(\Delta)=pr_2^*(x)\cup
Cl_{Y\times Y}(\Delta).$$
By definition of the PoincarŽ duality, it amounts to 
$$Tr_{Y\times Y}(x\cup \iota^*(y))=Tr_Y(pr_1^*(x)\cup Cl_{Y\times Y}(\Delta)\cup y)$$
One has $\Delta=\iota_*(Y)$, therefore by using PoincarŽ duality, we can rewrite the right hand
side as $Tr_{Y\times Y}(\iota^*\circ pr_1^*(x)\cup\iota^*(y))$, or $Tr_{Y\times Y}(x\cup
\iota^*(y))$, as desired. same for $pr_2$.

\begin{cor}
We have 
$$[D]=-\cup Cl_Y(D).$$
\end{cor}

\noindent{\bf Proof :} We apply the two previous lemmata, noticing that
$$pr_1^*(pr_2^*(x\cup Cl_Y(D))\cup Cl_{Y\times Y}(\Delta))=pr_1^*(pr_1^*(x\cup Cl_Y(D))\cup
Cl_{Y\times Y}(\Delta))$$
$$=x\cup Cl_Y(D).$$

Another particular correspondences used in the next, are given by cycles of the form
$D\times_X Y$ in $Y\times_X Y$ where $D$ is a relativ cycle in $Y$ of relative codimension
$r$. The action of such correspondence is given by the following diagram:

$$\begin{array}{lcr}R^mf_{s,*}\Omega^\bullet_{Y/X}&\stackrel{[D\times Y]}\longrightarrow
&
  R^{m-2r}f_{s,*}\Omega^\bullet_{Y/X}\\
\hbox{  }\hbox{  }\hbox{  }\downarrow&& \uparrow\hbox{  }\hbox{  }\hbox{  }\hbox{ 
}\hbox{  }\hbox{  }\hbox{  } \\  R^{2\delta-m}f_{s,*}\Omega^\bullet_{Y/X}&
\stackrel{-\cup D}\longrightarrow & R^{2\delta-m+2r}f_{s,*}\Omega^\bullet_{Y/X}
\end{array}$$

where the vertical maps are given by the polarization of the abelian scheme $Y$ wich
identifie each cohomology space with it's dual and by PoincarŽ duality.

\subsection*{II.2 The $\Z_{(p)}$-schematic version of Construction 5.1}

In this section, we consider dominant weights $\lambda$ for $(G,B,T)$ such that $s=\vert\lambda\vert$
satisfies $s+d<p-1$. We attach 
 to such  weights $\lambda$ a vector bundle
$\cV_\lambda$ with connection. Note that because of the need of compatibility with the transcendental
construction over $\C$ (using the restriction of the $G$-representation on $V_\lambda$ to the
Siegel parabolic), the definition will involve duals. We define first the vector bundle
$\cV_1$ associated to the standard representation
$V_1$ of $G$ as
$$\cV_1^\vee=R^1f_*\Omega_{A/X}^\bullet,$$  endowed with the Gauss-Manin connection.

We now use the sheaf-theoretic analogue of Construction 5.1  to define the dual of $\cV_\lambda$
 over $X$ and
$X_n$ as a direct factor in $R^\bullet f_{s\,*}\Omega^\bullet_{Y/X}$ cut out by algebraic
correspondences over
$\Z_{(p)}$. More precisely, we find an idempotent $e_\lambda$  in 
$\cC(Y/X)_{(p)}$ realizing this cut out:
$$\cV_\lambda^\vee =e_\lambda\cdot R^\bullet f_{s\,*}\Omega^\bullet_{Y/X}$$

The construction is in four steps: 
\begin{enumerate}

\item Project $R^\bullet f_{s\,*}\Omega^\bullet_{Y/X}$ to $(\cV_1^{\vee})^{\otimes s}$. This is
realized by the Liebermann trick. By KŸnneth formula, we have
$$R^\bullet f_{s\,*}\Omega^\bullet_{Y/X}=(R^\bullet f_{*}\Omega^\bullet_{A/X})^{\otimes s}$$
Moreover, since $A$ is an abelian scheme, one has

$$R^\bullet f_{*}\Omega^\bullet_{A/X}=\bigwedge^\bullet \cV_1^\vee$$

Therefore,
$$R^\bullet f_{s\,*}\Omega^\bullet_{Y/X}=\bigoplus_{0\leq j_1\leq 2g,\ldots,0\leq j_s\leq
2g}\bigwedge^{j_1}\cV_1^\vee\otimes\ldots\otimes\bigwedge^{j_s}\cV_1^\vee$$

The summand corresponding to $(j_1,\ldots,j_s)$ in the decomposition above is the kernel of the
correspondences on $Y$ given by
$[m_1]^*\times\ldots \times [m_s]^*-m_1^{j_1}\cdot\ldots\cdot m_s^{j_s}$ for all $m_1,\ldots,m_s\in
\Z$. Recall that we assumed also
$p>5$, hence $max(d,4)<p-1$ implies for any $g\geq 1$ that $2g<p-1$. Hence for any
$\alpha=1,\ldots,s$, we have
$j_\alpha<p-1$.Therefore by choosing $(m_1,\ldots,m_s)$ suitably (that is, with coordinates
generating $(\Z/p\Z)^\times$), we can construct an idempotent
$e_1$ in $\cC(Y/X)_{(p)}$ (of degree $0$) such that $e_1\cdot R^\bullet
f_{*}\Omega^\bullet_{A/X}=\cV_1^{\vee\,\otimes s}$.

Then, we realize the contractions $\phi_{i,j}$'s and their duals $\psi_{i,j}$'s defined in
Sect.5.1.1,  as
 algebraic correspondences in $\cC(Y/X)_{(p)}$.

\item The $\psi_{i,j}$'s:

For any $t\geq 1$, let $Y_t=A\times_X\ldots\times_X A$, $t$ times, and $f_t:Y_t\rightarrow X$ the
corresponding structural map. We abbreviate
$Y_s=Y$. Let $p_{i,j}: Y\rightarrow A\times A$ be the projection 
to the
$i$th and
$j$th components.
 Consider the Poincar\'e divisor $P$ in
$A\times_X A$ (corresponding to the PoincarŽ bundle).

\begin{de} The de Rham polarisation $\Psi_P\in \cV_1^{\vee\,\otimes 2}$ is defined as the projection
of
$Cl_{A\times A}(P)\in R^2f_{2,*}\Omega^\bullet_{A^2/X}$ to
$(Rf_*\Omega^\bullet_{A/X})^{\otimes 2}$ given by the KŸnneth formula.  
\end{de}

Consider the pull-back of $P$ by $p_{i,j}$; it is a divisor
$P_{i,j}$ in
$Y$. By 5.2.1, it defines a degree
$2$ endomorphism $[P_{i,j}]$ of $R^\bullet f_{s\,*}\Omega^\bullet_{Y/X}$. 
We have a commutative diagram
$$\begin{array}{ccccl}\cV_1^{\vee\,\otimes s-2}&\hookrightarrow &
R^{s-2}f_{s-2}\Omega^\bullet_{Y_{s-2}/X}&\hookrightarrow & R^{s-2}f_{s,*}\Omega^\bullet_{Y/X}\\
&&&&\\
\downarrow\Psi_{P,i,j}&&&& \downarrow -\cup Cl(P_{i,j})\\&&&&\\ \cV_1^{\vee\,\otimes s}&&
\hookrightarrow &&
R^{s}f_{s,*}\Omega^\bullet_{Y/X}\end{array}.$$ 
where the horizontal arrows are given by KŸnneth formula, and $\Psi_{P,i,j}$ consists in
inserting $\Psi_P$ at $i$th and $j$th indexes. Therefore, the morphism $\Psi_{P,i,j}$ is
induced by the divisor $P_{i,j}$.

\item The $\phi_{i,j}$'s: Consider the self-intersection $2g-1$ times of $P$; it is a
$1$-cycle on $A\times A$. Take its pull-back to $Y$ by the projection $p_{i,j}:Y\rightarrow A\times
A$ and again to $Y\times_X Y$ by the first projection $p_1:Y\times_X Y\rightarrow Y$. Then,
intersect this with the pull-back of the diagonal $\Delta_{s-2}$ in the self-product of the
remaining $s-2$ copies of $A$ in $Y$. The resulting cycle $Z_{P,i,j}$ is codimension
$\delta -1$ in
$Y\times_X Y$; therefore, it gives rise to a degree
$-2$ endomorphism of the cohomology.

\begin{de} Let $\Phi_P:\cV_1^{\vee\,\otimes 2}\rightarrow \cO_X$ be the linear dual of the projection
to
$(R^{2g-1}f_*\Omega^\bullet_{A/X})^{\otimes 2}$ by KŸnneth formula of $cl(P^{2g-1})\in
R^{4g-4}(f\times f)_*\Omega^\bullet_{A\times A/X}$.
\end{de}
Consider the contraction $\Phi_{P,i,j}:\cV_1^{\vee\,\otimes s}\rightarrow \cV_1^{\vee\,\otimes s-2}$
by
$\Phi_P$ at indexes $i$ and
$j$. We have a commutative diagram:

$$\begin{array}{lcr}\cV_1^{\vee\,\otimes s}&\hookrightarrow &
  R^{s}f_{s,*}\Omega^\bullet_{Y/X}\\

\downarrow\Phi_{P,i,j}&& \downarrow Z_{P,i,j}\\ \cV_1^{\vee\,\otimes s-2}&\hookrightarrow &
R^{s-2}f_{s,*}\Omega^\bullet_{Y/X}\end{array}.$$
Thus, $\Phi_{P,i,j}$ is given by the correspondence $Z_{P,i,j}$.   

\item Apply the Young symmetrizer
$c_\lambda$ to
$\cV_1^{\vee\,<s>}$. This projector has $\Z_{(p)}$-coefficients and belongs to a group algebra of
automorphisms of $f_s$, hence defines an element of $\cC(Y/X)$ as in 5.2.1. 

\end{enumerate}

Let us summarize the above constructions. For any dominant weight $\lambda$ of $G$ such that
$\vert\lambda\vert<p$, we associate a coherent locally free
$\cO_X$-module $\cV_{\lambda}$ such that 
\begin {itemize}
\item $\cV_1^\vee = R^1f_*\Omega_{A/X}^\bullet$ is associated to the
standard representation.
\item  $\cV_{\lambda}^\vee\otimes_{\Z_{(p)}}\C$ is the classical complex automorphic bundle associated
to
$\lambda$ (see for example \cite{CF} p.222).

\item Let us consider the additive functor $V\to {\cal V}^\vee$ from the semisimple category of
$G$-representations over $\Z_{(p)}$ of $p$-small weights  to the category 
of coherent locally free $\cO_X$-modules defined as above for simple objects. It is a functor of abelian categories
which commutes with tensor products and duality. This functor sends the $\phi_{i,j}$'s resp. $\psi_{i,j}$ of Sect.5.1.1 to
the $\Phi_{i,j}$'s resp. $\Psi_{i,j}$ of the present section.
\end{itemize}

\subsection*{II.3 The Gauss-Manin connection}

Over $\C$, the automorphic vector bundle $\cV_\lambda(\C)$ over  $S_U$ carries a natural
integrable connection given by the monodromy action $G(\Q)\to Aut(V_\lambda), g\mapsto
(v\mapsto g.v)$, where $V_\lambda$ est the irreducible $G(\C)$-representation of highest
weight $\lambda$. We call this connection the monodromy connection. 
To get an algebraic connection on the algebraic locally free $\cO_X$-module $\cV^\vee_\lambda$,
we first note that the sheaves ${\cal H}^m_{dR}(Y/X)=R^mf_{s,*}\Omega^\bullet_{ Y/X}$ are naturally
endowed with the Gauss-Manin connection (\cite{Katz}). We claim that this connection induces after
analytification, the monodromy connection. Indeed, we have just to verify this compatibility on
${\cal H}^1_{dR}(A/X)=R^1f_*\Omega^\bullet_{A/X}$. This implies in particular that the Gauss-Manin
connection commute to the idempotent used to define $\cV(\C)$. 

\begin{cor} Over $\Z_{(p)}$, the Gauss-Manin connection on
$\cV_1^\vee$ commutes to algebraic correspondences and therefore induces an integrable connection on
$\cV_{\lambda}$ ($\vert\lambda\vert<p$).
\end{cor}

{\bf Proof :}   Note that ${\cal H}^i_{dR}$ is locally free, hence commutes to
base-change: Cor.2 Chapt.2.5 of \cite{Mu}. We may replace $\Z_p$ by  $\C$ and the assertion follows
from the discussion above.

\subsection*{II.4 Canonical extension to toroidal compactification over $\Z_{(p)}$}

In the complex setting, Mumford (\cite{Mu2}, see also \cite{CF}, section VI.4) define a canonical extension
$\overline{\cV}_\lambda(\C)$ over $\overline X(\C)$ of the automorphic vector bundle $\cV_\lambda(\C)$. As explained by
Harris (\cite{Ha2}, (4.2.2)), this canonical extension is the extension provided by Deligne's existence theorem. As the
toroidal extension is defined over $\bQ$, we deduce that the extension is also defined over $\bQ$, we denote by
$\overline{\cV}_{\lambda,\Q}$ this extension over $\Q$, viewed as a coherent locally free module over ${\overline
X}_\Q={\overline X}\otimes_{\Z_p}\Q$. To extend this automorphic sheaves to $\Z_{(p)}$, we proceed as follows.

 First, consider
$$\begin{array}{ccc} A&\hookrightarrow & \overline{A}\\ \downarrow f
&&\downarrow\overline{f}\\X&\hookrightarrow &\overline{X}\end{array}$$ (for the construction of
$\overline{A}$ over $\Z[1/N]$, see Th.1.1 of IV.1
\cite{CF}) then, the canonical extension 
$\overline{\cV_1}^\vee$ of the standard sheaf $\cV_1=R^1f_*\Omega_{A/X}^\bullet$ to $\overline{X}$ is 
$$\overline{\cV}_1^\vee=R^1\overline{f}_*\Omega^\bullet_{\overline{A}/\overline{X}}(log\,\infty_{\overline{A}/\overline{X}})$$
 (where $\Omega^\bullet_{\overline{A}/\overline{X}}(log\,\infty_{\overline{A}/\overline{X}})$ denotes
the complex of relative differentials with relative logarithmic poles as defined in section 4.3). 

For $s<p$, let ${\overline
f}_s:\overline{Y}\to {\overline X}$ be a toroidal compactification of $f_s:Y\rightarrow X$. The canonical extension of
$R^s f_{s\,*}\Omega^\bullet_{Y/ X}$ to $\overline{X}$ can be defined in a similar way as before. That is 
$R^s {\overline f}_{s\,*}\Omega^\bullet_{{\overline
Y}/{\overline X}}(log\infty)$. Then, for a dominant weight $\lambda$ such that $\vert \lambda\vert =s<p$, the
canonical extension
$\overline{\cV}^\vee_\lambda$ of
${\cV}^\vee_\lambda$ is defined by 
$$\overline{\cV}^\vee_\lambda=j_*\overline{\cV}^\vee_{\lambda,\Q}\cap R^s {\overline
f}_{s\,*}\Omega^\bullet_{{\overline Y}/{\overline X}}(log\infty)$$

where $j : {\overline X}_\Q\to {\overline X}$ is the open immersion of the generic fiber 
${\overline X}_\Q$ in $\overline X$. $\overline{\cV}^\vee_\lambda$ is a coherent locally free ${\cal
O}_{\overline X}$-module,  direct factor of $R^s {\overline f}_{s\,*}\Omega^\bullet_{{\overline
Y}/{\overline X}}(log\infty)$ and $ \overline{\cV}^\vee_\lambda\otimes
_{\Z_p}\Q=\overline{\cV}^\vee_{\lambda,\Q}$. Moreover the Gauss-Manin connexion induces an integrable
connection on $\overline{\cV}^\vee_\lambda$. Note that this definition is legitimate by the
semisimplicity of the category of $G$-representions over $\Z_{(p)}$ with $p$-small weight (Lemma 7 of
sect. 5.1.1 with $G$ instead of $M$).

\vskip 5mm

{\bf Remark :} A better way to extend this automorphic sheaves is to extend the idempotents
$e_\lambda$ to the toroidal compactification: if
$\overline Y$ is a scheme and
$Y$ is an open subscheme, then there is an exact sequence (\cite{Ful} I.1.8):
$$ A_\bullet({\overline Y}-Y) \to A_\bullet({\overline Y})\to A_\bullet(Y)\to 0 $$

The natural way to extend a cycle of $Y$ to $\overline Y$ is to take it's closure. In the case of a
toroidal imbedding, Lemma 3.1. of \cite{Ha2} suggest to consider the normalization of the
closure. So we obtain correspondances ${\overline e}_\lambda$ over $\overline Y$. Unfortunatly, 
we can not see that
${\overline e}_\lambda$ is an idempotent. The problem is that the closure of the intersection of two cycles is not equal, in general, to
the intersection of the closure of this cycles. 

\subsection* {II.5 Automorphic bundles for the Levi $M$ }   

To every $B_M$-dominant weight $\mu$, one can also associate
${\cal W}_{\mu,n}$, a locally free ${\cal O}_{X_n}$-module; it is called the automorphic bundle
attached to
$\mu$. The construction is similar to the one sketched above. Consider the semiabelian scheme $f_{\cal
G}:{\cal G}\rightarrow
\overline{X}$ associated to our fixed toroidal compactification (see Th. 5.7, Chapt.IV of \cite{CF}),
which extends the universal abelian surface $f:A\rightarrow X$. Then, 
  the automorphic bundle on $X_n$ associated to the standard representation $W_1$ is
${\bf Lie}(A/X_n)^\vee$, and by part (3) of Theorem 5.7 of
\cite{CF} mentioned above, its canonical extension 
$\overline{\cal W}_{\mu,n}$ is ${\bf Lie}({\cal G}/\overline{X}_n)^\vee$. 
Then one uses the same trick as above to
construct $\overline{\cal W}_{\mu,n}$  
from the tensor product of ${\bf Lie}({\cal G}/\overline{X}_n)^\vee$ by itself
$s$-times. We note here that we can use the deep result of Harris (\cite{Ha2}, Th.4.2) 
to recover the rationality of the canonical extention of such automorphic vector bundles.

\vfill\eject

\vskip 5mm

{\small A. Mokrane, D\'epartement de Math\'ematiques, UMR 7539, Institut Galil\'ee, Universit\'e de
Paris 13, 93430 Villetaneuse. mokrane@math.univ-paris13.fr

\vskip 5mm 

J. Tilouine, Institut Universitaire de France et  D\'epartement de Math\'ematiques, UMR 7539,
Institut Galil\'ee, Universit\'e de Paris 13, 93430 Villetaneuse.  tilouine@math.univ-paris13.fr}

\end{document}